\documentclass[11pt]{article}
\usepackage{ulem}
\usepackage[makeroom]{cancel}
\usepackage{color}
\usepackage{xcolor}
\usepackage{amsfonts}
\usepackage{amssymb}
\usepackage{amscd}
\usepackage{amsmath}
\usepackage{epsfig}
\usepackage{graphicx, subfigure}
\usepackage{latexsym}
\usepackage{pst-3dplot}
\usepackage{verbatim}
\usepackage{tikz}
\usepackage{pdfsync}
\usepackage{hyperref}
\usepackage{multirow}
\usepackage[utf8]{inputenc}

\usepackage{color}
\usepackage{graphicx}

\usepackage{caption}
\usepackage{lipsum}
\newlength{\dinwidth}
\newlength{\dinmargin}
\newtheorem{definition}{Definition}
\newtheorem{theorem}{Theorem}

\newtheorem{proposition}{Proposition}
\newtheorem{corollary}{Corollary}
\newtheorem{remark}{Remark}
\newtheorem{lemma}{Lemma}
\newtheorem{example}{Example}

\def \surf{{\cal L}}

\def\Jac{\mathrm{Jacobian}}
\def \i{{\rm i}}

\def\px1{p_{x_1}}
\def\px2{p_{x_2}}
\def\pu1{p_{u_1}}

\def\Axi{A_{x_i}}
\def\Aum{A_{u_m}}
\def\Aualpha{A_{u_\alpha}}
\def\Auk{A_{u_k}}
\def\Aaj{A_{a_j}}

\def\Aai{A_{a_i}}
\def\Ainfty{A_\infty}

\def\Pxi{P_{x_i}}

\def\Pum{P_{u_m}}
\def\Pualpha{P_{u_\alpha}}
\def\Puk{P_{u_k}}
\def\Paj{P_{a_j}}
\def\Pak{P_{a_k}}

\def\Pinfty{{P_\infty}}

\def\Pualpha{P_{u_\alpha}}

\def\betaaj{\beta_{a_j}}
\def\betaai{\beta_{a_i}}
\def\betaxi{\beta_{x_i}}
\def\betaum{\beta_{u_m}}

\def\betaualpha{\beta_{u_\alpha}}

\def\c{\hat c}

\def\a{{\bf a}}
\def\b{{\bf b}}
\def\x{{\bf x}}

\def\B{\mathcal B}

\makeatletter
\@addtoreset{equation}{section}
\makeatother

\begin{document}
\title{Isoharmonic deformations and constrained Schlesinger systems}

\author{Vladimir Dragovi\'c$^1$ and Vasilisa Shramchenko$^2$}
\date{}

\maketitle

\footnotetext[1]{Department of Mathematical Sciences, University
	of Texas at Dallas, 800 West Campbell Road, Richardson TX 75080,
	USA. Mathematical Institute SANU, Kneza Mihaila 36, 11000
	Belgrade, Serbia.  E-mail: {\tt
		Vladimir.Dragovic@utdallas.edu}--the corresponding author}

\footnotetext[2]{Department of mathematics, University of
	Sherbrooke, 2500, boul. de l'Universit\'e,  J1K 2R1 Sherbrooke, Quebec, Canada. E-mail: {\tt Vasilisa.Shramchenko@Usherbrooke.ca}}

\

{\it Dedicated to the bicentennial of birth of P. L. Chebyshev (1821-1894), the founding father of the modern St. Petersburg school of mathematics.}

\

\begin{abstract}
	We introduce and study the dynamics of Chebyshev polynomials on $d>2$ real intervals. We define isoharmonic deformations as a natural generalization of the Chebyshev dynamics. This dynamics  is associated with  a novel class of constrained isomonodromic deformations for which we derive the constrained Schlesinger equations. We provide explicit solutions to these equations in terms of differentials on  an appropriate family of  hyperelliptic curves of any genus $g=d-1\ge 2$. The verification of the obtained solutions relies on the combinatorial properties of the Bell polynomials and on the analysis on the Hurwitz spaces. From the point of view of the classical algebraic geometry we formulate and solve the problem of constrained Jacobi inversion for hyperelliptic curves.  We discuss applications of the obtained results in integrable systems, e.g. billiards within ellipsoids in $\mathbb R^d$.
	
	\vskip 1cm
	
	
	
\end{abstract}

\

\newpage

\tableofcontents

\section{Introduction}
\label{sect_introduction}

Let us consider $d$ real intervals with $d>2$ and a point $y_0$ in the extended complex plane outside the intervals.   A harmonic measure  with a pole at $y_0$ is assigned to each of these $d$ intervals (see e.g. \cite{Si2015a}).  In this paper we introduce and study deformations of the endpoints of the intervals and of the  position of the  pole $y_0$ in a way that the harmonic measures of these deformed $d$ intervals with respect to the deformed pole $y_0$ remain unchanged under the deformation.
We will call such deformations {\it isoharmonic}  and understand them as deformations of a pair formed by a region and a point in it. Here, the region is the complement of the union of $d$ intervals  in the  extended complex plane and the point is the position of the pole $y_0$.

\smallskip

Since we assume that $d>2$, there is a hyperelliptic curve of genus $g=d-1\ge 2$ obtained as a double covering of the Riemann sphere  ramified over the $2d$ endpoints of the intervals.   The position $y_0$ of the pole of the harmonic measures defines  two points $Q_0, Q_0^*$ on the curve, paired by the hyperelliptic involution,  lying above $y_0$ on the double covering.   The  isoharmonic deformations     along with their   interpretation in potential theory  also  carry  a transparent algebro-geometric meaning  by inducing special deformations of the hyperelliptic curve, and   are related to a class of isomonodromic deformations of linear systems.

\smallskip

Let us first discuss the algebro-geometric context.  Assuming one of the endpoints of the intervals to be the point at infinity,
 the  isoharmonic deformations determine  a smooth family of hyperelliptic curves together with a section consisting of the points $Q_0$. More precisely, we will define a smooth family of hyperelliptic curves $\hat T: \hat{\mathcal  H}\to \hat X $ parameterized by the parameter set $\hat X$  introduced  below (see Sections \ref{sect_Chebyshev}, \ref{sect_isoharmonic},  and Definition \ref{def:Tfamily} in Section \ref{sect_Tfamily}). The fiber  $\mathcal H_{\bf x}$  over ${\bf x}\in \hat X$ is the projective closure of the algebraic curve of the equation
\begin{equation}
\label{T0}
\mu^2=\Delta^{{\bf  x}}_{2d -1}(z),
\end{equation}
where $\Delta^{{\bf  x}}_{2d-1}$ is  a polynomial of degree  $2d-1$. Taking a small enough neighbourhood $\hat X^0$ of ${\bf  x^0}\in \hat X$ we can choose a canonical homology basis in a consistent way for all curves $\mathcal H_{\bf x}$  with ${\bf x}\in \hat X^0$ by requiring that the projections of the basis cycles to the $z$-sphere be independent of $\bf x$. There is also an induced vector bundle $ V\to \hat X^0$ whose fiber over $\bf x$ is the vector space of holomorphic differentials on the curve $\mathcal H_{\bf x}$. Fix a basis  of  sections of this bundle normalized with respect to the $\bf a$-cycles of the chosen canonical homology basis.
Let us base the Abel map at the point at infinity $P_\infty$ of each curve. Then, there exists a section $s_{Q_0}$ of the family $\hat T$ of the hyperelliptic curves, such that
\begin{equation}
\label{section_cond}
\mathcal A_{P_\infty}(s_{Q_0}({\bf  x})) =\c_1+\mathbb B_{\bf x}\c_2.
\end{equation}
Here $\mathbb B_{\bf x}$ is the Riemann matrix of $\mathcal H_{\bf x}$ for the given choice of the homology basis and $\c_1, \c_2\in\mathbb R^g$. Then the isoharmonicity of the deformations implies that  $\hat c_1, \hat c_2$ are constant vectors independent of ${\bf x}\in \hat X^0$. By applying to $s_{Q_0}$ the hyperelliptic involution in the fibers, we get the associated section $s_{Q_0^*}$ of $\hat T\,.$

\smallskip

 The isoharmonic deformations are closely related to the theory of isomonodromic deformations of linear differential systems. In order to articulate and explore this relationship, we  introduce a novel class of isomonodromic deformations of Fuchsian systems, the so called {\it constrained } isomonodromic deformations.  Some  constrained isomonodromic deformations are  described by solutions to the {\it constrained Schlesinger system}.  We construct explicit families of such deformations in terms of differentials of the third kind on the fibers of the family $\hat T$ of hyperelliptic curves with simple poles along sections $s_{Q_0}$  and $s_{Q_0^*}$. In each fiber, the periods of the differential with respect to the chosen homology basis are constant multiples of $\c_1, \c_2$. Thus, the periods are constants {\it independent of the fiber}. The bridge connecting the potential theory side with the algebro-geometric framework is the fact that the above differentials of the third kind correspond to the differentials of the Green functions of the complement of the union of $d$ intervals with the pole at $y_0$.  The independence on the fibre of the periods of the differentials  then relates to variations of the region and the pole in an isoharmonic way.
\smallskip

A section of a family of {\it elliptic} curves defined by the condition \eqref{section_cond} with fixed constants $\c_1, \c_2$  arose in the context of isomonodromic deformations in the classical work of Picard \cite{Picard}. Picard provided a general solution of one of the Painlev\'e VI equations more than a decade prior to works of Gambier \cite{Gambier} and of R. Fuchs \cite{Fuchs}, who
derived the general form of such equations. The Painlev\'e VI equation is a second order ordinary differential equation
with parameters $\alpha, \beta, \gamma, \delta \in \mathbb C$, denoted PVI$(\alpha, \beta, \gamma, \delta)$, see eg. \cite{Okamoto}, \cite{Japan}, \cite{Ma1998}.   R. Fuchs derived it in 1907 \cite{Fuchs} in his
study of isomonodromic deformations of a Fuchsian linear system with four singularities at
 $\{u_1=0, u_2=1, u_3=x, u_4=\infty\} \subset {\mathbb {CP}}^1$

\begin{equation}
\label{linsys_intro}
\frac{d\Phi}{du} = A(u) \Phi, \qquad u \in {\mathbb C},
\end{equation}
for a $2\times 2$ matrix function $\Phi(u)$ defined in the complex
plane, where the matrix $A\in sl(2,{\mathbb C})$ is of the form:
\begin{equation*}
A(u) = \frac{A_1}{u} + \frac{A_2}{u-1} +
\frac{A_3}{u-x}.
\end{equation*}
As was already known to Painlev\'e in 1906 and  rediscovered  by Manin in \cite{Ma1998}, the Painlev\'e VI equations have an equivalent elliptic form:
$$
	\frac{d^2\, z}{d\,\tau^2}=\frac{1}{2\pi i}\sum_{j=0}^3\alpha_j\wp_z\left(z+\frac{T_j}{2},\tau\right),
	$$
where the transformed Weierstrass $\wp$-function  satisfies: $(\wp'(z))^2 =  \wp(z)\;(\wp(z) -1)\;(\wp(z) - x)$
and 	the parameters are related by:
	$$
	(\alpha_0, \alpha_1, \alpha_2, \alpha_3)=(\alpha, -\beta, \gamma, 1/2-\delta).
	$$

Let $\mathcal E_x$ be an elliptic curve represented as a two-fold covering of the Riemann sphere ramified over the set $\{0,1,x,\infty\}$ with $x\in\mathbb C\setminus\{0,1\}$. For some chosen canonical homology basis, let the Abel map be based at the point at infinity and let the Jacobian be generated by the Weierstrass vectors $2w_1, \,2w_2\in\mathbb C$. Then $T_j$ are the periods $(0, 2w_1, 2w_2, 2w_1+2w_2)$.

Now we can say more about the general solution of the Picard equation,  which is the Painlev\'e VI equation with parameters  $(\alpha_0, \alpha_1, \alpha_2, \alpha_3)=(0,0,0,0)$, that is  PVI$(0, 0, 0, 1/2)$.  Let us fix an arbitrary point $z_0$ in the Jacobian of $\mathcal E_x$. Then $z_0=2w_1\c_1+2w_2\c_2$ for some $\c_1, \,\c_2\in \mathbb R.$ The Jacobi inversion of $z_0$ gives a point $Q_0\in\mathcal E_x$ for which the projection on the base of the two-fold ramified covering is given by the image of $z_0$ under the Weierstrass $\wp$-function,  $y_0=\wp(z_0).$ If we now allow $x$ to vary, assuming that the projection of the homology basis on the base of the covering stays fixed, we need to define the corresponding variation of $z_0.$ The natural way to define this variation is the one leading to the Picard solution: $z_0$ is defined to be the point given by the above equality with $\c_1$ and $\c_2$ being fixed for all curves $\mathcal E_x$, that is $z_0(x)=2w_1(x)\c_1+2w_2(x)\c_2$ with $\c_1, \c_2$ independent of $x$. This gives the Picard solution $y_0=\wp(z_0(x))$ of PVI$(0,0, 0, 1/2).$

\smallskip

In other words, we have a family of elliptic curves $\pi:\mathcal E \to \mathcal B_{0,1}$ over the set $x\in\mathcal B_{0,1}:=\mathbb C\setminus\{0,1\}$ with a fiber given by an elliptic curve $\mathcal E_x$ as above. Points at infinity of the elliptic curves form a section of $\mathcal E$ which is taken as the zero for the  group law on the fibers. The points $Q_0$ defined as above form a local section $s_{Q_0}$ of $\mathcal E$.

\smallskip

We can then define a differential $\Omega$ of the third kind on $\mathcal E_x$ with a simple pole of residue $+1$ along the section $s_{Q_0}$ and a simple pole of residue $-1$ along the section $s_{Q_0^*}$ obtained from $s_{Q_0}$ by the elliptic involution of $\mathcal E_x$ interchanging the sheets of the covering; let the differential be normalized by the condition that its $a$-period with respect to the chosen homology basis is equal to $-4 \pi \i \c_2\,.$ It turns out, see \cite{DS2019}  and \cite{Hitchin1}, that the two zeros of $\Omega$, paired by the elliptic involution, form a local section of the family $\pi$ of elliptic curves whose projection on the base of the two-fold coverings representing the curves, seen as a function of $x\in \mathcal B_{0,1}$ coincides with the solution of PVI$(\frac{1}{8}, -\frac{1}{8}, \frac{1}{8}, \frac{3}{8})\,$. The constants $\c_1, \c_2$ play the role of the initial condition for the differential equation. Moreover, the relationship between the poles and zeros of $\Omega$ is given by the Okamoto transformation between solutions of PVI$(0,0,0, \frac{1}{2})$ and PVI$(\frac{1}{8}, -\frac{1}{8}, \frac{1}{8}, \frac{3}{8})\,$, see also \cite{Lor}.

\smallskip

Our objective now is to find a counterpart of solutions to PVI$(0,0,0, \frac{1}{2})$ and PVI$(\frac{1}{8}, -\frac{1}{8}, \frac{1}{8}, \frac{3}{8})\,$ in the case of a family of hyperelliptic curves. Thus, our starting point is to  find an analogue of the Picard solution for hyperelliptic curves. Given that the Painlev\'e VI equations describe isomonodromic deformations of linear systems \eqref{linsys_intro}, where the parameters $(\alpha, \beta, \gamma, \delta)$ are related to the eigenvalues of the residues of the matrix $A(u)\,,$ it is natural to expect that, for the sake of possible generalizations, the role of the Painlev\'e VI equation is played by the Schlesinger systems. Recall that the Schlesinger systems give the dependence of the residue matrices of $A(u) du$ on the positions of Fuchsian singularities of  linear system \eqref{linsys_intro}, with an arbitrary number of singularities, such that the linear system be isomonodromic.

\smallskip

Suppose now that we have a family  $T:\mathcal H \to X$ of hyperelliptic curves of genus $g>1$ over some parameter space $X$ where we can choose the canonical homology bases consistently for all the curves. Following the logic of the Picard solution, let us fix $\c_1, \c_2\in\mathbb R^g$ and consider a point $z_0=\c_1+\mathbb B_{\bf x} \c_2$ in the Jacobian of the fiber $\mathcal H_{\bf x}$ where ${\bf x} \in X$ and $\mathbb B_{\bf x}$ is the Riemann matrix of $\mathcal H_{\bf x}$ with respect to the chosen homology basis.
We then run into problems since a local section as in \eqref{section_cond} cannot be defined over an arbitrary family of hyperelliptic curves. This is due to the fact that the Abel map is only surjective on a surface of genus one. The Jacobi inversion of a  given generic point in the Jacobian of a hyperelliptic curve of genus $g$ gives a positive divisor of degree $g$. A generalization based on such a divisor is possible and was done in \cite{DS2017a}; it leads to isomonodromic deformations described by classical non-constrained Schlesinger systems.

\smallskip

In this paper, we focus on the situation completely opposite to a generic one, when a given point in the Jacobian of a hyperelliptic curve of genus $g$ is the image by the Abel map of a {\it single point} of the curve. Moreover, we want
 there to be a continuous variation of the point in the Jacobian such  that
 this condition be satisfied under a deformation of the curve and the Jacobian.
 We thus search for a special, non-generic family of compactified hyperelliptic curves of genus $g>1$, for which a section defined by \eqref{section_cond} exists for some fixed vectors $\c_1, \c_2\in\mathbb R^g.$
We will refer to this as {\it the problem of the constrained variations of the Jacobi inversion for curves of genus $g>1$.} It is shown in Section \ref{sect_isoharmonic} that the isoharmonic deformations generate such constrained variations of the Jacobi problem. In its turn, a solution to the constrained Jacobi problem allows us to construct solutions to the constrained Schlesinger system, see Theorem \ref{thm_main_g}.

\smallskip

The motivation for our study comes from the theory of extremal polynomials on $d$ real intervals.  We consider polynomials $\mathcal P_n$ of degree $n$ satisfying the Pell equation
\begin{equation}
\label{Pell0}
{\mathcal P}_n^2(z)-\Delta_{2d}(z){\mathcal Q}_{n-d}^2(z)=1.
\end{equation}
Here $\Delta_{2d}$ is a monic polynomial vanishing at the ends of the $d$  finite  real intervals and $\mathcal Q_{n-d}$  is a polynomial of degree $n-d$. In other words, $\mathcal P_n$ is a solution of the Pell equation, if there exist polynomials $\Delta_{2d}$ and  $\mathcal Q_{n-d}$ such that \eqref{Pell0} holds.  Solutions $\mathcal P_n$ of the Pell equation are also called the (generalized) Chebyshev polynomials;  they satisfy the  following extremality  conditions    \cite{Bogatyrev2012}.  By dividing out the leading coefficient of $\mathcal P_n$, we obtain a monic polynomial which has the least possible deviation from zero over the given set of intervals among all monic polynomials of the same degree and of the same signature. The signature takes into account the number of oscillations of the polynomial on each interval, i.e. this is a $d$-vector of integers, where the $j$-th component counts the number of zeros of ${\mathcal Q}_{n-d}$ in the $j$-th interval. The corresponding set of intervals is then the maximal subset of the real line on which $\mathcal P_n$ has the above minimal deviation property. The study of such polynomials was initiated by Chebyshev for the case of one interval. His student Zolotarev initiated the study of such polynomials on two intervals. Important results were obtained by Markov, Bernstein, Borel, Akhiezer in the early 20th century. The study of such problems significantly expanded with varied applications and continues nowadays \footnote{Lebesgue wrote on June 17 1926: ``{\it I assume I am not the only one who does not understand the interest  in and significance of these strange problems of maxima and minima studied by Chebyshev in memoirs whose titles often begin with, 'On functions deviating least from zero...' Could it be that one must have a Slavic soul to understand the great Russian Scholar?}\,'' (See \cite{SY1992}, \cite{Si2015}.)}, see eg. \cite{Akh4}, \cite{AhiezerAPPROX}, \cite{Widom}, \cite{SY1992}, \cite{PS1999}, \cite{Bogatyrev2012}, \cite{Si2011}, \cite{Si2015} and references therein.

\smallskip

We show in Section \ref{sect_Chebyshev} that constrained variations of the Jacobi inversion problem are naturally induced by the dynamics of Chebyshev polynomials for $d>2$ intervals, which we introduce there as well. This dynamics provides an important class of isoharmonic deformations.

\smallskip

The Chebyshev dynamics for $d=2$  was linked in \cite{DS2021} to Hitchin's discovery \cite{Hitchin1} of a direct connection between PVI$(\frac{1}{8}, -\frac{1}{8}, \frac{1}{8}, \frac{3}{8})$ and the Poncelet Closure Theorem. Griffiths and Harris in \cite{GriffithsHarris} presented classical results about the Poncelet Theorem in a modern language and attracted attention of contemporary mathematicians to this gem of classical projective and algebraic geometry.  Poncelet trajectory is a polygonal line whose vertices lie on a conic, called {\it boundary}, and whose segments are tangent to another conic, called {\it caustic}. The Poncelet theorem says that if, for the given pair of conics, there exists a closed Poncelet trajectory consisting of $n$ segments, then there is such a trajectory passing through every point of the boundary conic.
 Any pair of conics generates a {\it tangential pencil of conics}, that is a one-parameter family of conics inscribed in the four common tangents of the two given conics.  There exists an elliptic curve constructed using this pencil and one distinguished conic from the pencil. If we choose the caustic to be this distinguished conic, then the boundary conic determines a point on the elliptic curve.  The Poncelet trajectories are closed of length $n$ if the pair caustic-boundary corresponds to a point of order $n$ of the elliptic curve.  A  criterion which tells  if a pair of conics determines a point of order $n$ was derived by Cayley \cite{Cay}. This criterion was effectively used by Hitchin in his construction of algebraic solutions to PVI$(\frac{1}{8}, -\frac{1}{8}, \frac{1}{8}, \frac{3}{8})$.

\smallskip

Hitchin's observation has been one of many significant appearances of the Poncelet Theorem  in various mathematical contexts, from approximation theory, integrable systems \cite{MV1991}, to algebraic geometry of stable bundles on projective spaces \cite{NT}, billiards, Jacobians of elliptic and hyperelliptic curves, numerical ranges, (para)orthogonal polynomials, including works of several outstanding mathematicians, like Jacobi, Cayley, Darboux, Trudi, Lebesgue, Griffiths, Harris, Berger, Narasimhan, Kozlov, Duistermaat, Simon and many others. Poncelet Theorem
also has an important mechanical interpretation. The elliptical billiard \cite{KT1991, Tab2005book} is a dynamical system where a material
point of unit mass is moving frictionlessly, with constant velocity inside an ellipse, obeying the reflection law at the boundary, that is having congruent angles of incidence and reflection. Any segment of a given elliptical billiard trajectory is tangent to the same conic, confocal with the
boundary, called {\it the caustic} of a given trajectory \cite{DR2011}. If a trajectory becomes closed after $n$ reflections, then the
Poncelet Theorem implies that any trajectory of the billiard system, which shares the same caustic curve, is also periodic with
the period $n$. This setting generalizes naturally to the Poncelet Theorem and billiards within quadrics in any dimension $d>2$.
Our Chebyshev dynamics from this paper with $d>2$ intervals corresponds to such billiards in $d$-dimensions exactly in the same way the Hitchin construction corresponded to the Poncelet theorem in $d=2$, see Section \ref{sect_Poncelet}.   Let us note that the first consideration of billiards within quadrics in $\mathbb R^d$ for $d>2$ goes back to an isolated short note \cite{Dar1870} of Darboux from 1870, where the case $d=3$ was studied.
The area of integrable billiards in higher dimension started to develop intensively in 1990's, in works of Moser, Veselov \cite{MV1991}, Chang, Crespi, and Shi \cite{CCS} and many others (see \cite{DR2011, DR2019b} and references therein).

\smallskip

 Exploiting further the connection with potential theory and conformal geometry, in Section \ref{sec_SC}, we apply  a Schwarz--Christoffel mapping generated by the Green function, and map the upper-half plane to a semi-strip with vertical slits, the {\it comb region}. The idea of conformal maps to the comb regions goes back to Marchenko-Ostrovsky \cite{MO1975}, see also e.g. \cite{SY1992, EY}. It appears to be very effective in the study of isoharmonic deformations. It provides an explicit rectification of
 such deformations, see Theorem \ref{thm:SC}:
Any isoharmonic deformation  preserves the horizontal base of the comb, together with the foot-points of vertical slits. Thus the deformation manifests only in changing of the heights of the vertical slits.  In other words,  the
Schwarz--Christoffel mapping rectifies the isoharmonic deformations, transforming them into a vertical dynamics along the slits, which belong to fixed vertical rays, while the rays do not change under the deformations.

\smallskip

Let us conclude  this part of  the introduction by mentioning various approaches to non-constrained Schlesinger systems based on algebraic geometry, see e.g. \cite{DIKZ}, \cite{KiKo}, \cite{Lor}, \cite{DGS2021}, or conformal geometry, see e.g. \cite{CKLW, EG}.

\medskip

{\bf Organization of the paper.}
\smallskip

In Section \ref{sect_Chebyshev} we introduce the dynamics of Chebyshev polynomials, the {\it Chebyshev dynamics}, and link it to a solution to the constrained variation of the Jacobi inversion problem.

In Section \ref{sect_isoharmonic}, we provide the necessary information from potential theory and introduce in
 Definition \ref{def:isoharmonic} the isoharmonic deformations. Lemma  \ref{lemma:z0fixed} relates the isoharmonic deformations with the constrained variation of the Jacobi inversion.
In Section \ref{sect_constrainedSS}, we introduce a novel class of isomonodromic deformations and derive the constrained Schlesinger system, see Theorem \ref{th:constrainedSS} and the system of equations \eqref{constrained}.
Section \ref{sect_surfaces} fixes the notation from the theory of hyperelliptic curves, provides a definition of $\mathcal T$-families of hyperelliptic curves naturally associated with the Chebyshev dynamics, see Definition \ref{def:Tfamily},
and sets the stage for the formulation of Theorem \ref{thm_main_g}. Theorem  \ref{thm_main_g} is stated in Section \ref{sect_solution} and gives  explicit solutions to the constrained Schlesinger system \eqref{constrained} in terms of differentials on underlying hyperelliptic curves.

\smallskip

A significant part of the paper is devoted to the proof of Theorem \ref{thm_main_g}. The proof is quite involved. There are two main technological ingredients of the proof. The first is related to combinatorial  properties of the Bell polynomials. The basics on Bell polynomials can be found for example in \cite{Andrews}. The results we  derived for Bell polynomials for the purpose of proving Theorem \ref{thm_main_g} are collected in Section \ref{sect_Bell}. The second main technological tool in the proof is the calculus on Hurwitz spaces of hyperelliptic surfaces. We refer to \cite{Fay92} and \cite{KokoKoro} for the background  on the Rauch variational formula and related material. We assembled the results  extending the Rauch variation to $\mathcal T$-families of curves   in Section \ref{sect_variational} in the form in which they are used in Section \ref{sect_proof} devoted to the proof of Theorem \ref{thm_main_g}.

\smallskip

Section \ref{sect_variational} also shows, see Remark \ref{rmk_dynamics}, how  $\mathcal T$-families  for rational values of the parameters $\c_1, \c_2\in \mathbb Q^{d-1}$ provide  solutions to the Chebyshev dynamics,  introduced in Section \ref{sect_Chebyshev}.  For the general value of the parameter $\c_1 \in \mathbb R^{d-1}$ with $\c_2=0$   we get a solution  of the isoharmonic deformations and thus of a more general problem of constrained variation of Jacobi inversion for hyperelliptic curves of any genus.
Section \ref{sect_variation_pts} contains Theorem \ref{th:jacobiconst} which describes solutions to the problem of constrained variation of Jacobi inversion.
\smallskip

In Section \ref{sect_Poncelet}, we explain the relationship between the Chebyshev dynamics on $d>2$ intervals and billiards in $d$-dimensions.  This generalizes Hitchin's work \cite{Hitchin1} and interrelates the results obtained in previous sections with the theory of integrable systems,  in particular in terms of integrable billiards in $d$-dimensional space.  It also employs further  conformal geometry and rectifies the isoharmonic deformations in Theorem \ref{thm:SC}.  We conclude with a brief discussion about the injectivity of the frequency map in the context of the inheritance problem in Section \ref{sec:inher}.

\section{Dynamics of Chebyshev polynomials over $d>2$ intervals}
\label{sect_Chebyshev}

It is well known that an arbitrary choice of $d$ real intervals does not guarantee the existence of a solution to the corresponding Pell equation \eqref{Pell0}. A set of intervals for which a solution does exist is called the {\it support} of a Chebyshev polynomial. Supports of Chebyshev polynomials of degree $n$ are called {\it $n$-regular} in \cite{SY1992}. The supports which can be obtained from one another by an affine change of variables are considered equivalent. To remove the freedom of such a change of variables, we assume one of the intervals to be fixed at $[0,1]\,.$ Two  important theorems of Peherstorfer and Schiefermayr (\cite{PS1999}, Th. 2.7 and Th. 2.12) imply  that if we start with a support of a Chebyshev polynomial and vary one endpoint of each of the remaining $d-1$ intervals, the positions of the other ends of those intervals are uniquely determined by the condition of solvability of the Pell equation in the class of polynomials  of the original Chebyshev polynomial,  that is  with the  same degree and signature.
This brings us to the following question.

\begin{itemize}
\item {\it Given a set of $d\ge 3$ real intervals which support polynomial solutions of the Pell equation \eqref{Pell0}, while keeping one interval fixed and   varying  one endpoint of each of the remaining $d-1$ intervals, how to describe the dynamics of the remaining $d-1$ endpoints, under the condition that the Pell equation remains solvable during the entire process with the same degree of the Chebyshev polynomial and the same signature? We call this variation of support of Chebyshev polynomials} {\rm the  Chebyshev dynamics}.
\end{itemize}

From the point of view of potential theory, the supports of the Chebyshev polynomials are characterized as unions of intervals each of which has a rational equilibrium measure (see \cite{SY1992}  and Section \ref{sect_isoharmonic}). Taking into account that the equilibrium measure is the harmonic measure with respect to the point at infinity, and that rational numbers which deform continuously remain fixed, we see that the above dynamics provides one instance of the isoharmonic deformations.

As it is natural to associate a hyperelliptic curve with a set of $d\geq 3$ intervals, the above Chebyshev dynamics defines a very special family of hyperelliptic curves. Under our assumption, all curves of this family are ramified over the points $0$ and $1$. Thus this family is parameterized by the positions of $d-1$ branch points; the remaining $d-1$ branch points being functions of the $d-1$ independently varying ones. The genus of our curves is thus $g=d-1.$ It is by studying this family of curves that we are able to answer the above question concerning the dynamics of supports of Chebyshev polynomials,  see Remark \ref{rmk_dynamics} in Section \ref{sect_variation_pts}.

More precisely,  aforementioned  Peherstorfer-Schiefermayr theorems \cite{PS1999}  allow  us to define a family of hyperelliptic curves $\hat T: \hat{\mathcal  H}\to \hat X $ parameterized by the set $\hat X:=\{\hat x_1, \dots, \hat x_{d-1}\,|\, \hat x_1<\hat x_2<\dots<\hat x_{d-1}\}\subset (\mathbb R\setminus [0, 1])^{d-1}$ such that the fiber $\hat{\mathcal H}_{\bf x}$ over ${\bf x}=(\hat x_1, \dots, \hat x_{d-1})$ is the projective closure of the algebraic curve of the equation
\begin{equation}
\label{T}
\mu^2=z(z-1)\prod_{j=1}^{d-1}(z-\hat x_j)(z-\hat u_j)=:\Delta^{\bf x}_{2d}(z),
\end{equation}

where  $\hat x_j$ and $\hat u_j$ are real numbers smaller than $\hat x_{j+1}$ and $\hat u_{j+1}$ for $1\leq j\leq d-2$ with the $\{\hat u_j\}$ being functions of $\{\hat x_j\}$ such that
\begin{equation}
\label{support}
\left(\underset{\substack{1\leq j\leq d-1\\ \hat x_j< \hat u_j}}{\bigcup}[\hat x_j, \hat u_j]\right)\cup \left(\underset{\substack{1\leq j\leq d-1\\ \hat x_j> \hat u_j}}{\bigcup}[\hat u_j, \hat x_j]\right)
\end{equation}
 is the support of a Chebyshev polynomial.  In other words, according to the Peherstorfer-Schiefermayr Theorem 2.12 from \cite{PS1999}, the $x_j$ may be either the left endpoint or the right endpoint of the interval it belongs to. In order to include all such possible options of orders between points $x_j$ and $u_j$, we may replace the parameter space $\hat X$ by the pair $\hat X=(\hat X, \sigma)$ where $\sigma: \{1, 2, \dots, d-1\}\rightarrow \{\ell, r\}$. Now, the points $\hat x_j, \, j\in {\sigma^{-1}(\ell)}$ occupy the left endpoints of the intervals $[\hat x_j, \hat u_j]$ while $\hat x_k, \, k\in{\sigma^{-1}(r)}$ occupy the right endpoints of the intervals $[\hat u_k, \hat x_k]$.

 \smallskip

 The family of curves \eqref{T}  admits two sections at infinity, we denote them $s_{\infty^+}$ and $s_{\infty^-}$.   We have $\mu\sim z^{d}$ locally near $s_{\infty^+}$ and $\mu\sim -z^{d}$ at $s_{\infty^-}$.
 The section $s_{\infty^+}$ consists  of points of the same finite order, while $s_{\infty^-}$ consists of the base points of the Abel maps.
 This is the result of the following Lemma, see, for example, \cite{Bogatyrev2012}.
\begin{lemma}
\label{lemma_Abel_Chebyshev}
Let us assume that the family $\hat T$ of compactified hyperelliptic curves is restricted to some subset $\hat X^0\subset \hat X$ which is small enough to allow for a canonical homology basis to be chosen consistently for all curves of the family in such a way that the projections of the cycles onto the $z$-sphere are the same for all ${\bf x} \in \hat X^0$.
Let $s_{\infty^-}(\x)$ be the base for the Abel map $\mathcal A_{s_{\infty^-}(\x)}$ on the fiber $\hat{\mathcal H}_{{\bf x}}$ of the family $\hat T:\hat {\mathcal H}\to \hat X^0$.  A solution to the Pell equation \eqref{Pell0}  with $\Delta_{2d}$ given by \eqref{T}  exists if and only if  $s_{\infty^+}$ is a section of $\hat T$ of order $n$ where $n$ is the degree of the corresponding Chebyshev polynomial, that is
\begin{equation}
\label{Abel_Chebyshev}
n\mathcal A_{s_{\infty^-}(\x)}(s_{\infty^+}(\x))\equiv 0  \;\;(\Jac(\hat{\mathcal H}_{\bf x}))
\end{equation}
for all  ${\bf x} \in \hat X^0\,.$
\end{lemma}

Before we provide a proof of the lemma, let  us define the {\it Akhiezer function} and the related differential
\begin{equation}
\label{eq:Akhiezer}
\mathbb A(P):=\mathcal P_n(z) +\mu \mathcal Q_{n-g-1}(z); \qquad \Omega_{\mathbb A}= \frac{1}{n}d\log \mathbb A
\end{equation}
on each fiber $\hat{\mathcal H}_{\bf x}$ of the family $\hat T$. Here $P=(z,\mu)$ is a point on $\hat{\mathcal H}_{\bf x}$  and $\mathcal P_n\,,\;\mathcal Q_{n-g-1}$ are the polynomials satisfying the Pell equation \eqref{Pell0}  with $\Delta_{2d}$ given by \eqref{T}. By this definition,  $\mathbb A(P)$ is a meromorphic function on $\hat{\mathcal H}_{\bf x}$ with a pole of order $n$ at $s_{\infty^+}(\x)\,.$

{\it Proof of Lemma \ref{lemma_Abel_Chebyshev}.}
Due to the Pell equation, we have
$\mathcal P^2_n(z)-\mu^2(z)\mathcal Q^2_{n-g-1}(z)=1$
for each ${\bf x} \in \hat X^0\,.$  Applying the hyperelliptic involution $(z, \mu)^*=(z, -\mu)$, we obtain for the Akhiezer function $\mathbb A(P)\,:$
\begin{equation*}
\mathbb A(P^*)\mathbb A(P) = (\mathcal P_n -\mu \mathcal Q_{n-g-1})(\mathcal P_n +\mu \mathcal Q_{n-g-1})=1
\end{equation*}
and thus $\mathbb A(P^*)=\frac{1}{\mathbb A(P)}\,.$ Therefore we conclude that  $\mathbb A(P)$ has a zero of order $n$ at $s_{\infty^-}(\x)$. By the Abel theorem, the existence of a function with a pole of order $n$ at $s_{\infty^+}(\x)$ and a zero of order $n$ at $s_{\infty^-}(\x)$  implies the statement of the lemma. $\Box$

Let us call {\it {\bf T}-curves} the fibers of the family $\hat T:\hat {\mathcal H}\to \hat X^0$, the compactified hyperelliptic curves defined by the projectivization of equations \eqref{T} subject to condition \eqref{support} and such that a canonical homology basis can be chosen for all curves as in Lemma \eqref{lemma_Abel_Chebyshev}. The notion of {\bf T}-curve coincides with the notion of {\it Toda curve}, see for example the classical McKean's survey \cite{McK1979} for their   important role in  spectral theory and integrable systems. We use the letter T
 as a common initial of Toda and Tchebysheff, recalling the traditional Western (French) transliteration of Chebyshev.

\smallskip

Note that in genus one, which corresponds to $d=2$, the $\bf T$-curves have one dependent branch  point $u_1$ and one independent branch point $x_1$. In this case, the Akhiezer function \eqref{eq:Akhiezer} is related by a M\"obius transformation in the $z$-sphere fixing $0,1$ and sending $u_1$ to $\infty$  to the function $s$ used by Hitchin in \cite{Hitchin1} to establish the direct link between PVI$(\frac{1}{8}, -\frac{1}{8}, \frac{1}{8}, \frac{3}{8})$ and the Poncelet theorem. In Section \ref{sect_Tfamily}, we perform this M\"obius transformation  to move the section $s_{\infty^+}$ of a finite order into the affine part of the $\bf T$-curves. The resulting family of curves is  an example of a $\mathcal T$-family, according to Definition \ref{def:Tfamily} from Section \ref{sect_surfaces}.

Let us mention also that in the case of two intervals, $d=2,$ exactly one critical point of the Chebyshev polynomial falls in the gap between the two intervals and is called {\it gap critical point}. The image of this point under the mentioned M\"obius transformation as a function of the image of $x_1$  satisfies PVI$(\frac{1}{8}, -\frac{1}{8}, \frac{1}{8}, \frac{3}{8})\,,$ see \cite{DS2021}.
Since the Painlev\'e-VI equation is equivalent to the Schlesinger system in the matrix dimension $2\times 2$, the dynamics of Chebyshev polynomials defined on two real intervals is related to the Schlesinger isomonodromic deformations of a linear  Fuchsian system for a $2\times 2$ matrix. Motivated by this relationship, we consider the isomonodromic deformations of a $2\times 2$ Fuchsian system naturally associated with our $\mathcal T$-family of hyperelliptic curves (Definition \ref{def:Tfamily}, Section \ref{sect_surfaces})  in the case of $d\geq 3$ intervals and show that this leads to a generalization of the Schlesinger system, which we call the  constrained Schlesinger system, see Section \ref{sect_constrainedSS}.

\section{Isoharmonic deformations}
\label{sect_isoharmonic}

In the case $g=1$ and $\c_1, \c_2\in \mathbb Q$  in \eqref{section_cond},  the important roles were played  in \cite{Hitchin1}, \cite{DS2019}, and \cite{DS2021}   by  Hitchin's  function $s$ and  an associated  differential $(d\log s)/n$ for the constructions of algebraic solutions of Painlev\'e VI equation  PVI$(\frac{1}{8}, -\frac{1}{8}, \frac{1}{8}, \frac{3}{8})\,$. Hitchin's function  and the differential are related by a M\"obius transformation in $z$ to the Akhiezer function and the differential from \eqref{eq:Akhiezer}.  Without the assumption that $\c_1, \c_2$ are both rational, the Hitchin function $s$ does not exist. However, there exists  a differential of the third kind $\Omega$  (see \cite{DS2019}, formula (7))  which  naturally generalizes $(d\log s)/n$.
For $g>1$, in order to extend the considerations of Chebyshev dynamics from Section \ref{sect_Chebyshev}  to the cases of irrational $\c_1, \c_2$, we  employ  potential theory and harmonic analysis, see \cite{Si2011, Si2015a}.

\smallskip

Let us start with an arbitrary union of $d$  finite  real intervals
$$
E=[c_{2d},c_{2d-1}]\cup[c_{2d-2},c_{2d-3}]\cup\dots\cup[c_2,c_1] \,\quad\text{with}\quad  c_{2d}<c_{2d-1}<\dots<c_1 .
$$
A generic set $E$ does not support a Chebyshev polynomial, therefore the associated Akhiezer function is not defined. However, let us  consider the  compact  curve corresponding to the equation
\begin{equation}\label{eq:hyperisoharmoinic}
\mu^2=\prod_{j=1}^{2d}(z-c_j)\,.
\end{equation}
Let us again denote the two points at infinity of the curve by $\infty^+$ and $\infty^-$, where we have $\mu\sim z^{d}$ locally near $\infty^+$ and $\mu\sim -z^{d}$ at $\infty^-$. Now introduce
the differential of the third kind $\eta$  defined on this curve  having simple poles at the two points at infinity of the curve and subject to  the conditions:
\begin{equation}\label{eq:diff}
\int_{c_{2j+1}}^{c_{2j}}\eta =0, \qquad j=d-1,\; d-2, ..., 1,
\end{equation}
with some normalization, see \cite{Bogatyrev2012}: we can choose the normalization so that the residues of $\eta$ be $\pm 1$ at $\infty^{\pm}$.
\begin{remark}\label{rem:eta} We refer to the intervals $(c_{2j+1}, c_{2j})$ which belong to the complement of $E$ as {\it the gap intervals}.
The differential $\eta$ is of the form:
$$
\eta=\frac{k(z)}{\mu},
$$
where $k$ is a real polynomial of degree $d-1$. The polynomial $k$ has one zero in each of the $d-1$ gap intervals because of the condition
\eqref{eq:diff}. Thus it has exactly one zero in each of the gap intervals and no zeros outside the gap intervals. From there we also see that the polynomial
$k$  has a constant sign in each of the intervals $[c_{2k}, c_{2k-1}]$, $k=1, \dots, d$.
\end{remark}

\smallskip

{\it The equilibrium measure} ${\mathcal M}_E$ \cite{Si2011} is defined by:
$$
{\mathcal M}_E([c_{2k}, c_{2k-1}])=\frac{1}{\pi}\int_{c_{2k}}^{c_{2k-1}}|\eta|.
$$
We are interested in the behaviour of the equilibrium measure when some of the endpoints $\{c_{2d}, \dots, c_1\}$ of the intervals vary.
We are investigating variations which  keep one of the intervals unchanged and   keep one endpoint of all other intervals unchanged as well   . In total, we assume $d+1$ endpoints to remain unchanged  and their type as the right or the left endpoint to remain also unchanged. Thus, let us denote by $\tilde x$ the  subset of  $d+1$ elements  of the set  $\{c_j|j=1, \dots,2d\}$ which remain unchanged: both endpoints of one of the $d$ intervals and exactly one of the endpoints of the remaining intervals. Denote by $\hat u=(\hat u_1, \dots, \hat u_{d-1})$  the remaining $d-1$ endpoints, which are subject to variations.

\smallskip
Following \cite{DR2019b}, we define the map  $F_{\tilde x}: \mathbb R^{d-1} \to \mathbb R_{\geq 0}^{d-1}$ by
\begin{equation}\label{eq:freq}
F_{\tilde x}(\hat u):= (f_1, \dots, f_{d-1}),
\end{equation}
with
$$
f_j=\sum_{k=d+1-j}^{d}{\mathcal M}_E([c_{2k}, c_{2k-1}]),\quad j=1, \dots, d-1.
$$
We call the map $F_{\tilde x}(\hat u)$  the frequency map and its components $f_j$ the frequencies, following \cite{DR2019b}, because of their
interpretation in the theory of integrable billiards. We will say a bit more about this interpretation in Section \ref{sec_billiard}.

\smallskip
The frequency map $F_{\tilde x}(\hat u_1, \dots, \hat u_{d-1})$   with $\tilde x=\{\tilde x_1, \dots, \tilde x_{d+1}\}$ fixed
is  a local diffeomorphism.  This property was proved in
\cite{DR2019b}, Theorem 13, by using considerations similar  to the proof of the Bogataryev-Peherstorfer-Totik Theorem (Theorem 5.6.1 from \cite{Si2011}).

 The above property of the frequency map implies that if  ${\tilde x}$ is  fixed and  $(f_1, \dots, f_{d-1})$ are given, then ${\hat u}$  is uniquely determined via \eqref{eq:freq}. This property is a key ingredient in the definition of isoharmonic deformations; it generalizes the Peherstorfer-Schiefermayr results (Theorems 2.7 and 2.12 from \cite{PS1999} mentioned above) to the case when the set of intervals $E$ does not support  Pell's equation. The cases when $E$ does support  Pell's equation are characterized by the property that all frequencies are rational (see \cite{SY1992} and Section \ref{sec_billiard}). It turns out  that the differentials $\Omega_{\mathbb A}$ from \eqref{eq:Akhiezer} and $\eta$ from \eqref{eq:diff} essentially coincide in the case where all the frequencies  are rational. See Section \ref{sec:inher} for an additional
 comment on the frequency map and its injectivity.

\smallskip

We now establish a generalization of the fact that points at infinity on {\bf T}-curves are of a finite order. Thus, the following Lemma generalizes Lemma \ref{lemma_Abel_Chebyshev}:

\begin{lemma}
\label{lemma:z0fixed}
With an appropriate choice of a canonical homology basis for the curve \eqref{eq:hyperisoharmoinic}, the following relations connect the Abel map on the curve \eqref{eq:hyperisoharmoinic} of $\infty^+$  and the frequencies:
 $$\pm(\mathcal A_{\infty^-}(\infty^+))_j=\i f_j,\, j=1, \dots, d-1.$$
\end{lemma}

\smallskip

The proof follows from bilinear relations for differentials of the first and third kind, see for example \cite{Sp1957}, Theorem 10-6 and Corollary 10-4, taking into account the defining properties of $\eta$, see Remark \ref{rem:eta}.

\smallskip
We now introduce {\it one class of  isoharmonic deformations} as generalizations of the Chebyshev dynamics from Section \ref{sect_Chebyshev}.
Assume the frequencies $(f_1, \dots, f_{d-1})$ given together with $d+1$ endpoints $\tilde x$ as above. Two of these $d+1$ endpoints are the endpoints of the same interval (say $\tilde x_1, \tilde x_2$) and the remaining  $d-1$ endpoints denote by $\hat x$.  Given that the frequency map \eqref{eq:freq} is a local diffeomorphism,  the frequencies   $(f_1, \dots, f_{d-1})$ and the endpoints $\hat x$ uniquely determine the endpoints $\hat u_1, \dots, \hat u_{d-1}$, assuming that the type of the endpoint as a left or right endpoint of each of the points $\hat u$ is prescribed. Now, we start to deform smoothly the $d-1$ endpoints   $\hat x$ while keeping the frequencies $(f_1, \dots, f_{d-1})$ and the pair of endpoints $(\tilde x_1, \tilde x_2)$ unchanged. We define now the remaining $d-1$ endpoints $\hat u$ as functions of $\hat x$, i.e.
$$\hat u=\hat u(\hat x)= (\hat u_1(\hat x), \dots, \hat u_1(\hat x)),$$
such that
$$
F_{(\tilde x_1, \tilde x_2, \hat x)}(\hat u(\hat x)):= (f_1, \dots, f_{d-1}).
$$
Let $E_{\hat x}$ denote the union of $d$ intervals with the endpoints $\{\tilde x_1, \tilde x_2, \hat x, \hat u(\hat x)\}$ obtained from $E$
as just described. We will say that the deformation of the complement of $E$ into the complement of $E_{\hat x}$  is an {\it isoequilibrium deformation}. The complements are defined with respect to the extended complex plane.\\

Following \cite{Widom}, we introduce some further notions of  potential theory. We will consider a domain $V$ in the extended complex plane with the boundary consisting of  sufficiently smooth  Jordan curves $E_1, \dots, E_d$   (of class $C^{1+}$ in \cite{Widom}). Real Green's function with a pole at $w_0\in V$, denoted by $g(z, w_0)$, is defined by the following conditions:
\begin{itemize}
\item [(a)] $g(z, w_0)$ is harmonic in $V\setminus \{w_0\}$;
\item [(b1)] if $w_0\neq\infty$ then $g(z, w_0)+\ln|z-w_0|$ is harmonic around $w_0$;
\item [(b2)] if $w_0=\infty$ then $g(z, w_0)-\ln|z|$ is harmonic around $w_0=\infty$;
\item [(c)] $\lim_{z\rightarrow \zeta} g(z, w_0)=0$ for all $\zeta\in E_k, \, k=1,\dots, d$.
\end{itemize}

Given a boundary function $\psi$ on $E=E_1\cup \dots\cup E_d$, the Green function resolves the Dirichlet problem for $V$. Namely:
$$
h(z)=\frac{1}{2\pi}\int_E\psi(s)\frac{\partial} {\partial n_s} g(s, z)|ds|
$$
is harmonic in $V$ having $\psi$ as its boundary function on $E$. Here $n_s$ is the unit normal at $s\in E$ directed toward $V$. In the particular case of $\psi=\chi_{E_k}$, the characteristic function of $E_k$, one gets the so-called {\it harmonic measure} of $E_k$ corresponding to $w_0$:
\begin{equation}\label{eq:harmonicmeasure}
\omega_k(w_0)=\frac{1}{2\pi}\int_{E_k}\frac{\partial} {\partial n_s} g(s, w_0)|ds|.
\end{equation}
As harmonic functions, the Green functions and harmonic measures have their harmonic conjugates, and determine corresponding holomorphic functions. For the obtained holomorphic functions, the Green functions and harmonic measures are their real parts. For example, denote ${\hat g}(z, w_0)$ the harmonic conjugate
of $g(z, w_0)$ and $G(z, w_0)=g(z, w_0)+i {\hat g}(z, w_0)$ will also be referred to as the Green function.
The function $f=\exp(G)$ will be called {\it the complex Green function} \cite{Akh4}. If $w_0=\infty$ then the harmonic measure coincides with the equilibrium measure (see  \cite{Si2011}).

 One can find more about Green functions, harmonic measures and potential theory for example in \cite{Si2011, Si2015a, Walsh} and references therein. The Green functions and harmonic measures are conformal invariants.

Let us now  set $w_0=\infty$ and  apply the M\"obius transformation $\rho$ defined by $\rho(c_{2d})=0, \; \rho(c_{2d-1})= 1, \; \rho(c_1)=\infty$. We denote the images of the other endpoints by $x_j$ and $u_j$ as before so that each interval connects one of
the $x_j$'s to one of the $u_j$'s  and such that $\rho(c_2)=x_g$, where $g=d-1$ is the genus of the compactified curve \eqref{eq:hyperisoharmoinic}.
 Thus $\rho$ takes  the complement of $E=[c_{2d},c_{2d-1}]\cup[c_{2d-2},c_{2d-3}]\cup\dots\cup[c_2,c_1]$ to the complement of $E_{\bf x}=[0,1]\cup [x_1, u_1]\cup\dots \cup [x_{g-1}, u_{g-1}] \cup [x_{g}, \infty]$,  where and ${\bf x}=(x_1, \dots, x_g)$.  Here we assume that $u_j$ are the right endpoints for simplicity of notation.  Let us  denote the image of $w_0=\infty$ by $y_0:=\rho(w_0)$.

Note that the M\"obius transformation can be alternatively defined by sending any dependent right endpoint to infinity;
our definition $\rho(c_1)=\infty$ is a choice.

The curve \eqref{eq:hyperisoharmoinic} is thus transformed into the curve of the equation
\begin{equation}
\label{rho_curve}
v^2= {\Delta} (u) = u(u-1)\prod_{j=1}^{g-1}(u-x_j)(u-u_j)(u-x_g)
\end{equation}
with $u=\rho(z).$

The differential $\eta$ becomes  the differential $\hat \eta_{\bf x}$ (see for example \cite{Widom}, p. 227):
\begin{equation}\label{eq:etax}
\hat \eta_{\bf x}= \frac {\hat k(u) du}{\sqrt{{\Delta} (u)} (u-y_0)},
\end{equation}
where $\hat k$ is a polynomial of degree $g$ determined by the conditions
$$
\int_1^{x_1}\hat \eta_{\bf x}=\int_{u_j}^{x_{j+1}}\hat \eta_{\bf x}=0, \qquad j=1, \dots, g-1,
$$
and
$$
\hat k(y_0)=-\sqrt{{\Delta} (y_0)}.
$$

 The differential $\hat\eta_{\bf x}$ has simple poles at the points that are $\rho$-images of $\infty^+$ and $\infty^-$, let us denote them by $Q_0$ and $Q_0^*$, respectively. Thus $Q_0, Q_0^*$ are points on the hyperelliptic curve \eqref{rho_curve} above the point $y_0$ related to each other by the hyperelliptic involution.

  The condition of preservation of the equilibrium measure of the intervals  transforms to the condition of preservation of the harmonic measures with a pole at  $y_0$ of the intervals $\hat E_1=[0, 1]$, $\hat E_j=[x_j, u_j],$ $E_g=[x_g, \infty]$ with $j=1, \dots, g-1$, when $\bf x$ varies.
When $\bf x$ varies, we have a family of curves \eqref{rho_curve}; recall that $u_j$'s become functions of $\bf x$. Assuming, as in Lemma \ref{lemma_Abel_Chebyshev}, that the variation of $\bf x$ is small enough to allow for a consistent choice of a canonical homology bases in all the curves of the family, we can consider the corresponding family of Jacobians.

According to Lemma \ref{lemma:z0fixed}, the condition of preservation of the above harmonic measures of the intervals is equivalent
to the constancy of the coordinates of the point $z_0:=\mathcal A_{Q_0^*} (Q_0)$ over the family of the Jacobians of the compactified curves \eqref{rho_curve}.
Denoting now by $P_\infty$ the point at infinity of each curve of the family \eqref{rho_curve}, and noting that $A_{P_\infty}(Q_0) = -A_{P_\infty}(Q_0^*)$, we obtain the constancy of the vectors $\c_1, \c_2$ from \eqref{section_cond} for the family of curves.

\begin{definition} \label{def:isoharmonic} We call the deformations of the complement of $E_{\bf x}$ and the marked point  $y_0$, (which is the same as the deformations of $({\bf x}, y_0)$)  isoharmonic if the harmonic measures
$\omega_k(y_0)=\int_{\hat E_k}\hat \eta$ , $k=1, \dots, g$, are preserved.
\end{definition}

 Note that  the Green function for the complement of $E_{\bf x}$ with the pole at $y_0$ is given by
$$
G_{E_{\bf x}}(z, y_0)=\int_{0}^z\hat \eta_{\bf x}.
$$

\begin{corollary}\label{rem:injectivity}
For a fixed
value of the parameter $(u_1, \dots, u_g, y_0)=({\bf u}, y_0)$, the map of the argument $\bf x=(x_1, x_2, \dots, x_g)$
$$
\hat F_{({\bf u}, y_0)}(x_1, x_2, \dots, x_g) = (\omega_1(y_0), \dots, \omega_g(y_0))
$$
is invertible.
\end{corollary}

This statement can be seen as a corollary of Theorem 13 from \cite{DR2019b}. The proof goes along
the lines of the proof of the above mentioned  Bogataryev-Peherstorfer-Totik Theorem (Theorem 5.6.1 from \cite{Si2011}) once we observe the monotonicity of the M\"obius transformation $\rho$.  \\

 When $y_0$ remains fixed as the point at infinity ($y_0=\infty$) under an isoharmonic deformation, then the isoharmonic deformation is an isoequilibrium deformations, considered in the first part of this section.

\section{Constrained Schlesinger system}
\label{sect_constrainedSS}

The classical Schlesinger system introduced in  \cite{Schlesinger} is an integrable nonlinear system describing monodromy preserving deformations in the class of non-resonant matrix Fuchsian systems with $N+1$ logarithmic singularities. It is closely related to two problems, both called the Riemann-Hilbert inverse monodromy problem: one requiring to find a Fuchsian system with prescribed monodromy and another one requiring to find a matrix function with Fuchsian singularities at the $N+1$ points and prescribed monodromy at those points. In the case of matrix dimension two, the Schlesinger system reduces to the Garnier system \cite{Garnier1912, Garnier1926}, and is equivalent to a Painlev\'e VI equation if $N=3$.

In this paper, we set the matrix dimension to be two and consider the isomonodromic deformation of the non-resonant Fuchsian system for the function $\Phi(z)\in M_{2\times 2}(\mathbb C)$ with Fuchsian singularities at $z_1, \dots, z_N$ and $\infty: $
\begin{equation}
\label{Fuchsian}
\frac{d\Phi}{dz} = \sum_{k=1}^N\frac{A_k}{z-z_k}\Phi=A(z)\Phi\;, \qquad z\in\mathbb C\setminus\{z_1, \dots, z_N\}.
\end{equation}
Denoting by $A_\infty$ the residue of $A(z)dz$ at $z=\infty$, we have $A_\infty=-\sum_{k=1}^N A_k\,.$
Recall that monodromy, in matrix dimension two, refers to the monodromy representation $\rho:\pi_1(\mathbb C\setminus \B) \to GL(2, \mathbb C)$ for the discrete set $\B:=\{z_1, \dots, z_N\}\subset \mathbb C$ resulting from the analytical continuation of solutions to the linear system of ordinary differential equations \eqref{Fuchsian}
along closed paths in $\mathbb C\setminus \B$ based at some point in $\mathbb C$ away from $\B.$  Here $A_k$ are matrices independent of $z$.  Equivalently, the meromorphic $1$-form
\begin{equation}
\label{1form}
A(z)dz:= \sum_{k=1}^N \frac{A_k}{z-z_k}dz
\end{equation}
can be seen as the connection form of a flat meromorphic connection with simple poles at $\B\cup \infty$ in the trivial rank two vector bundle over $\mathbb CP^1$. In this case the monodromy representation is the holonomy representation of the connection.  The isomonodromic deformation problem is the problem of finding the dependence of the residue matrices $A_k$ on positions $z_k$ of singularities which results in the monodromy representation being constant under small variations of the set $\B$.

A system \eqref{Fuchsian} is non-resonant if eigenvalues of each residue matrix $A_k$ do not differ by an integer. In that case, a fundamental matrix of solutions behaves as follows close to a singularity:
\begin{equation}
\label{Phi}
\Phi(z)\underset{z\sim z_k}{=} \left( G_k+\mathcal O(z-z_k)\right)(z-z_k)^{T_k} C_k
\end{equation}
with some matrices $G_k, T_k, C_k$ independent of $z$. The residue matrices $A_k$ are then given by $A_k=G_kT_kG_k^{-1}$ and the monodromies of $\Phi$ are $M_k=C_k^{-1}{\rm e}^{2\pi\i T_k}C_k\,.$ The requirement of isomonodromy is thus equivalent to requiring the matrices $T_k$ and $C_k$ to be constant under small variations of the positions of singularities $z_1, \dots, z_n\,.$ In this case, differentiating \eqref{Phi} with respect to $z$ and $z_j$, we get the behaviour of the derivatives of $\Phi$ close to $z_k$
\begin{equation}
\label{Phi_sys}
\frac{d\Phi}{dz} \Phi^{-1} \underset{z\sim z_k}{=} \frac{A_k}{z-z_k} + \mathcal O(1)
\qquad\mbox{and}\qquad
\Phi_{z_k} \Phi^{-1} \underset{z\sim z_k}{=} - \frac{A_k}{z-z_k} + \mathcal O(1)\,,
\end{equation}
moreover, $\Phi_{z_k} \Phi^{-1}$ does not have singularity away from $z=z_k\,.$ Assuming a normalization of $\Phi$ which implies that $\Phi_{z_k} \Phi^{-1}$ vanishes at infinity, we see that a fundamental matrix solution $\Phi$ of the Fuchsian system \eqref{Fuchsian} satisfies also the following system
\begin{equation}
\label{iso}
\Phi_{z_k} \Phi^{-1} = - \frac{A_k}{z-z_k},  \qquad k=1,\dots, N\,.
\end{equation}
The Schlesinger system can then be derived as compatibility condition of \eqref{Fuchsian} and \eqref{iso}:
\begin{equation}
\label{Schlesinger}
\frac{\partial A_j}{\partial z_k} =
\frac{[A_k, A_j]}{z_k-z_j}\,, \qquad\qquad
\frac{\partial A_k}{\partial z_k} = -\sum_{j\neq k}\frac{[A_k,A_j]}{z_k-z_j},
\end{equation}
where  the second equation is equivalent to the condition $ A_{\infty}=-A_1 -\dots-A_N  = const.$
Solutions of \eqref{Schlesinger} provide residue matrices $A_k$ as functions of positions of Fuchsian singularities for which the linear system \eqref{Fuchsian} is isomonodromic.

Motivated by our families $\hat T$ and $T$ of hyperelliptic curves \eqref{T} and \eqref{TMobius}, we suggest the following generalization of the isomonodromic deformation problem. Let us fix two positions of Fuchsian singularities of the system \eqref{Fuchsian} at $0$ and $1$ and split the rest of the set $\B$ into two subsets, denoted by  $\{x_j\}_{j=1}^K$ and $\{u_j\}_{j=1}^{2g-K-1}$ so that $2g+1$ stands for the total number of elements in $\B$ and $K<2g$ is a positive  integer. Assume now that $x_j$ are allowed to vary independently whereas $u_j$ are functions of $\{x_j\}_{j=1}^K\,.$ The {\it constrained isomonodromic deformation problem} is the question of finding residue matrices $A_k$ and the functions $u_j(x_1, \dots, x_K)$ such that the monodromy representation of $\pi_1(\mathbb C\setminus \B)$ induced by the Fuchsian system \eqref{Fuchsian} stays constant under small variations of $\{x_j\}_{j=1}^K\,.$

\smallskip

Let us denote by $a_j$ an arbitrary element of the set $\B=\{0,1,x_1, \dots, x_{K}, u_1, \dots, u_{2g-1-K}\}$ and by $\Aaj$ the residue matrix corresponding to the singularity at $a_j.$ Then  the Fuchsian system \eqref{Fuchsian} takes the form for $u\in\mathbb C\setminus \B\,:$
\begin{equation}
\label{Fuchsian-2}
\frac{d\Phi}{du}=\Big( \frac{A_0}{u}+\frac{A_1}{u-1}+\sum_i \frac{A_{x_i}}{u-x_i}+\sum_j \frac{A_{u_j}}{u-u_j}\Big) \Phi\,.
\end{equation}
 A fundamental matrix of the Fuchsian system \eqref{Fuchsian-2} locally behaves as in \eqref{Phi}, and, assuming that matrices $T_k$ and $C_k$  are constant, the above derivation yields the following system for $\Phi:$
\begin{equation}
\label{iso-2}
\Phi_{x_i}\Phi^{-1}=-\frac{A_{x_i}}{u-{x_i}}-\sum_j \frac{A_{u_j}}{u-{u_j}}\frac{\partial u_j}{\partial x_i}.
\end{equation}
Computing now the compatibility condition of \eqref{Fuchsian-2} and \eqref{iso-2}, we obtained a system of equations for the residue matrices as functions of independent variables $x_1, \dots, x_K$. A solution to this system defines an isomonodromic Fuchsian system of the form \eqref{Fuchsian-2}. We call this compatibility condition the {\it constrained Schlesinger system} to reflect the fact that some positions of Fuchsian singularities are constrained to be functions of the independently varying ones.

\begin{theorem}\label{th:constrainedSS}
Denote by $a_j$ an arbitrary element of the set  $\B=\{0,1,x_1, \dots, x_{K}, u_1, \dots, u_{2g-1-K}\}$.  The constrained Schlesinger system has the form:

\begin{eqnarray}
\label{constrained}
&&\partial_{x_i} A_{a_j} = \frac{[A_{x_i}, A_{a_j}]}{x_i-a_j} + \sum_{k=1}^{2g-1-K}\frac{[A_{u_k}, A_{a_j}]}{u_k-a_j} \frac{\partial u_k}{\partial x_i}\,, \quad\mbox{for}\quad a_j\notin\{x_i, u_1, \dots, u_{2g-1-K}\};
\nonumber\\
&&\partial_{x_i} A_{u_m} = \frac{[A_{x_i}, A_{u_m}]}{x_i-u_m} +\!\!\!\!\sum_{\substack{k=1\\k\ne m}}^{2g-1-K} \frac{[ A_{u_k},A_{u_m}]}{u_k-u_m} \frac{\partial u_k}{\partial x_i}-\frac{\partial u_m}{\partial x_i}\!\!\sum_{\substack{a_j\in \B\\a_j\ne u_m}} \frac{[A_{a_j},A_{u_m}]}{a_j-u_m}\,,\;\;\mbox{for}\;\; 1\leq m\leq 2g-1-K;\qquad
\\
&&\partial_{x_i} A_{x_i} = -\sum_{\substack{a_j\in \B \\a_j\ne x_i}}\frac{[A_{x_i}, A_{a_j}]}{x_i-a_j} + \sum_{k=1}^{2g-1-K}\frac{[A_{u_k}, A_{x_i}]}{u_k-x_i} \frac{\partial u_k}{\partial x_i}\,.
\nonumber
\end{eqnarray}
\end{theorem}
Here again the last equation can be replaced by $\sum_{a_j\in \B} A_{a_j}=const.$

Note that we can interpret this system in two ways. First, we can see it as an underdetermined system if we say that
a solution to system \eqref{constrained} is a set of matrices $\Aaj$ together with derivatives $\frac{\partial u_m}{\partial x_i}$ for $a_j\in \B\,,$ and indices $i=1,\dots, K $ and $m=1, \dots, 2g-1-K\,.$ On the other hand, we can fix some functions $u_m(x_1, \dots, x_K)$ and look for the matrices $\Aaj$ such that \eqref{constrained} is satisfied. In the latter case, we obtained a determined system.

In this paper we construct a solution to the  constrained Schlesinger system in the case where the number of dependent variables $u_k$ is one less that the number of independent variables $x_j$, that is $K=g$,  and the matrices $\Aaj$ are traceless $2\times 2$ with eigenvalues $\pm \frac{1}{4}\,.$
 Our solution is constructed in terms of functions and differentials defined on the compact curves of the family \eqref{rho_curve} with their branch points playing the role of independent variables $x_k$ and of the functions $u_k(x_1, \dots, x_g)$.

\section{Surfaces associated with the constrained Schlesinger system}
\label{sect_surfaces}

\subsection{$\mathcal T$-family of hyperelliptic surfaces}
\label{sect_Tfamily}

Consider the family $\hat T:\hat {\mathcal H}\to \hat X^0$ of ${\bf T}$-curves from Section \ref{sect_Chebyshev}. This family admits a section $s_{\infty^+}$ consisting of points at infinity on the curves which are of a finite order, that is satisfy condition \eqref{Abel_Chebyshev}.
Let us apply M\"obius transformation  to move the section $s_{\infty^+}$ into the affine part of the {\bf T}-curves. To simplify the exposition, let us assume the following ordering of the endpoints of the right-most interval in the support of Chebyshev polynomials corresponding to our {\bf T}-curves:   $\hat x_{d-1}<\hat u_{d-1}\,.$ This assumption is not necessary as everything which follows can be adapted to any ordering of the endpoints.
Under this assumption, the  M\"obius  transformation in the $z$-sphere
\begin{equation}
\label{Mobius}
\rho(z) = \frac{z(1-\hat u_{d-1})}{z-\hat u_{d-1}}
\end{equation}
 sending $\{0,1,\hat u_{d-1}\}$ to $\{0,1,\infty\}$ is increasing on the set of the branch points of the {\bf T}-curves. Seeing the {\bf T}-curves \eqref{T} as two-fold ramified coverings of the $z$-sphere and applying  $\rho$ in each sheet of the fibers of the family $\hat T$,  we  obtain a new family of hyperelliptic curves which can be described as $ T: {\mathcal  H}\to  X^0 $ parameterized by a subset $X^0:=\rho(\hat X^0)$ of the set $ X:=\{ x_1, \dots,  x_{d-1}\,|\,  x_1<x_2<\dots< x_{d-1}\}\subset (\mathbb R\setminus [0, 1])^{d-1}$ such that the fiber over $( x_1, \dots,  x_{d-1})$ is the projective closure of the algebraic curve of the equation
\begin{equation}
\label{TMobius}
v^2=u(u-1)\prod_{j=1}^{d-1}(u- x_j)\prod_{j=1}^{d-2}(u- u_j)
\end{equation}
where $u=\rho(z)$ and  $\{ u_j=\rho(\hat u_j)\}\subset \mathbb R$ are functions of $\{ x_j=\rho(\hat x_j)\}\subset \mathbb R$ such that the set $\left(\cup_{j=1}^{d-2}[\rho^{-1}(x_j),  \rho^{-1}(u_j)]\right)\cup [\rho^{-1}(x_j), \rho^{-1}(\infty)]$ is the support of a Chebyshev polynomial. Here, for notational simplicity, let us assume that $[x,y]$ stands for the interval between ${\rm min}\{x,y\}$ and  ${\rm max}\{x,y\}\,.$

The canonical homology bases in the fibers of $\hat T:\hat {\mathcal H}\to \hat X^0$ transform by $\rho$ into canonical homology bases in the fibers of $ T: {\mathcal  H}\to  X^0 $, such that
 the projections of the basis cycles onto the $u$-sphere are independent of $\x\in X^0\,.$ Let us denote the obtained canonical basis in the homology of $\mathcal H_\x$ by $\{\a_1, \dots, \a_g; \b_1, \dots, \b_g\}$ without keeping track of the $\x\in X^0\,.$ From now on, we assume such a basis to be chosen.  Here as before, $g=d-1$ is the genus of the curves.

 Let us now fix some notation. The fibers $\mathcal H_\x$ of the family $T$ with $\x\in X^0$  seen as
 two-fold ramified coverings $u: \mathcal H_\x\to \mathbb CP^1$ are ramified over the set
 \begin{equation}
 \label{B}
 B:=\{0,1,x_1, \dots, x_g; u_1, \dots, u_{g-1}\}
 \end{equation}
  and the point $u=\infty$. Let us use notation $a_j$ for points of the set $B$, that is $ B=\{a_j\}_{j=1}^{2g+1}\,.$ We call the points of the set $B$ the branch points of the curves $\mathcal H_\x\,.$

Capital letters $P$ and $Q$ will be used to denote points on the curves $\mathcal H_\x$, for example $P=(u, v)$ and we then write $u$ for $u(P)$. For the ramification points of the covering $u$ we use the notation $P_{a_j} = (a_j, 0)$ and $P_\infty=(\infty, \infty)\,.$
Ramification points form sections of the family $T$, which we denote in the same way as ramification points themselves, for example, $P_\infty:X^0\to \mathcal H$ with $P_\infty(\x)=P_\infty\in\mathcal H_\x\,.$

Introduce the {\it standard local coordinates} on the surface $\mathcal H_\x$ as follows:
\begin{align}
\label{coordinates}
&\zeta_k(P)=\sqrt{u(P) - a_k} \quad\mbox{if}\quad P\sim P_{a_k},
\nonumber
\\
& \zeta_\infty(P)= \frac{1}{\sqrt{u(P)}}\quad\mbox{if}\quad P\sim P_\infty,
\\
& \zeta(P)=u(P)-u(Q) \quad\mbox{if}\quad P\sim Q \quad\mbox{and $Q$ is a regular point.}
\nonumber
\end{align}
Let us now denote
\begin{equation}
\label{y0}
y_0=\rho(\infty)
\end{equation}
and for each curve $\mathcal H_\x$ with $\x\in X^0$ denote by $Q_0(\x)$ the image of $s_{\infty^+}(\x)$ under $\rho;$ we have that $u(Q_0(\x))=y_0.$ Note that this is well defined as $\rho$ preserves the sheets of the covering. These points $Q_0(\x)$ define a section of the family $T$ which we denote also by $Q_0$.
By applying the hyperelliptic involution $(u,v)\mapsto (u, -v)$ on $\mathcal H_\x$ to the points $Q_0(\x)$, we obtain the section $Q_0^*$ of $T$, consisting of the images of $s_{\infty^-}(\x)$ under the M\"obius transformation.

Let ${\bf \omega} =(\omega_1, \dots, \omega_g)^t$ be the vector of holomorphic differentials on  $\mathcal H_\x$ normalized with respect to the above canonical homology basis by the condition
\begin{equation}
\label{hol_norm}
\oint_{\a_j}\omega_k= \delta_{jk}\,.
\end{equation}
 These differentials form sections of the vector bundle over $X^0$ whose fiber over $\x\in X^0$ is a space of all holomorphic differentials on $\mathcal H_\x\,.$ Note that the differentials $\omega_j$ on $\mathcal H_\x$ are linear combinations of $\frac{du}{v},\frac{udu}{v}, \dots, \frac{u^{g-1}du}{v}$  and thus we have $\omega(P^*)=-\omega(P)$ for all points $P\in\mathcal H_\x\,.$ We use $\mathbb B_\x$ to denote the matrix of $\b$-periods of $\{\omega_k\}$, the Riemann matrix of $\mathcal H_\x.$

By the change of variables given by $\rho$, the condition \eqref{Abel_Chebyshev} for $s_{\infty^+}(\x)$ to be the point of order $n$ becomes
\begin{equation}
\label{Mobius_finite}
n\mathcal A_{s_{\infty^-}(\x)}(s_{\infty^+}(\x))=n\int_{Q_0^*(\x)}^{Q_0(\x)}\omega =  2n\int_{P_\infty(\x)}^{Q_0(\x)}\omega  \equiv 0
\end{equation}
where the second equality is obtained using $\omega(P^*)=-\omega(P)\,.$ Basing the Abel map on $\mathcal H_\x$ at $P_\infty(\x)$, we obtain a point $Q_0(\x)$ of order $2n$ on each fiber of the family $T$. Thus we obtain that the section $Q_0:X^0\to\mathcal H$ is of a finite order, that is formed by points of a finite order. We can rewrite \eqref{Mobius_finite} by introducing rational vectors $\c_1, \c_2\in \mathbb Q^g$ as follows:
\begin{equation}
\label{q0}
\int_{P_\infty(\x)}^{Q_0(\x)}\omega =\c_1+\mathbb B_\x \c_2\,.
\end{equation}
It is important to note that, in this relation, while $\omega, Q_0$ and $\mathbb B$ depend on $\x\in X^0$, the vectors $\c_1$ and $\c_2$ are constant. This is because the point  $Q_0(\x)$ is of finite order on $\mathcal H_\x$ for any $\x$ and therefore $\c_1$ and $\c_2$ are rational for any $\x$ thus cannot vary continuously with $\x$.   Let us also note that $\c_1$ and $\c_2$ cannot be simultaneously half-integer vectors as, by our construction, the points $Q_0$ do not coincide with a ramification point. In other words, we have $Q_0\neq Q_0^*$ for the curves of the family $T$.

In the case $d=2$ the family $T$ reduces to the genus one family $\pi:\mathcal E \to \mathcal B$ from \cite{Hitchin1} discussed in the introduction with one  independently  varying branch point $x_1=x$ and no dependent branch points. In this case, $u(Q_0)$ gives rise to the Picard solution of PVI$(0,0, 0, \frac{1}{2}).$ However, a section $Q_0$ satisfying \eqref{q0} exists on any family of elliptic curves. In higher genera, it is non-trivial to describe a family of curves admitting a section satisfying \eqref{q0}. This,  along with the consideration in Section \ref{sect_isoharmonic}, motivate us to introduce the following terminology.

\begin{definition}
\label{def:Tfamily}
A triple $(T, s_\infty, s)$ is called $\mathcal T$-family, if $ T: {\mathcal  H}\to  X$ is a smooth fibration with fibers given by compactified hyperelliptic curves  ${\mathcal  H}_\x$ for $\x\in X$ with a consistent choice of canonical homology bases in all fibers, $s_\infty:X\to\mathcal H$ is a section such that $s_\infty(\x)$ is a point at infinity of $\mathcal H_\x$ for all $\x\in X$ and
$s:X\to \mathcal H$ is a section of $T$ such that
\begin{equation}
\label{c1c2_general}
\mathcal A_{s_\infty(\x)}(s(\x)) =\c_1+\c_2\mathbb B_\x,
\end{equation}
where $\mathcal A_{s_\infty(\x)}$ is the Abel map of $\mathcal H_\x$ based at $s_\infty(\x)$ for the given choice of the homology bases,  $\mathbb B_\x$ is the corresponding Riemann matrix of $\mathcal H_\x$,  and $\c_1, \c_2\in\mathbb C^g$ are constant vectors independent of $\x\in X\,.$
\end{definition}

From Section \ref{sect_Chebyshev}, we see that $\bf T$-curves form a natural $\mathcal T$-family $(\hat T, s_\infty, s)$ with the sections $s=s_{\infty^+}\,,\; s_\infty=s_{\infty^-}$ and $\c_1, \c_2$ being real vectors with rational components. Section \ref{sect_isoharmonic} shows how to construct a $\mathcal T$-family with  $\c_1, \c_2$ having   irrational components  as well.

\subsection{Abelian differentials on the hyperelliptic curves}
\label{sect_differentials}

Consider one hyperelliptic curve $\mathcal H_\x$ from the $\mathcal T$-family $(T, P_\infty, Q_0)$ of Section \ref{sect_Tfamily}. This curve is defined by equation \eqref{TMobius}, where the set of branch points is  $B=\{a_j\}_{j=1}^{2g+1}=\{0,1,x_1, \dots, x_{g}, u_1, \dots, u_{g-1}\}$. Recall that there is a chosen canonical homology basis for $\mathcal H_\x$ denoted by $\{\a_1, \dots, \a_g; \b_1, \dots, \b_g\}$.  Here we list the Abelian differentials on $\mathcal H_\x$ which will be useful for us. From now on, we drop the dependence on $\x$ in our notation, writing, for example, $P_\infty$ and $Q_0$ for points on $\mathcal H_\x\,.$

\begin{paragraph}{Holomorphic differentials}
We have already introduced a basis of normalized $1$-forms \eqref{hol_norm} in the space of holomorphic differentials. We also need
the following holomorphic non-normalized differential on $\mathcal H_\x$:
\begin{equation}
\label{phi}
\varphi(P) = \frac{du}{\sqrt{\prod_{a_j\in B} (u-a_j)}} \,.
\end{equation}
Recall that we write $u$ for $u(P)$ with $P\in\mathcal H_\x$. We also need to introduce the concept of evaluation of Abelian differentials at a point of the Riemann surface. We define the evaluation of an Abelian differential $\Upsilon$ at a point $Q\in\mathcal H_\x$ as the constant term of the Taylor series expansion of the differential with respect to the standard local parameter $\xi$ from the list \eqref{coordinates} at $Q$, that is
\begin{equation}
\label{evaluation}
\Upsilon(Q) = \frac{\Upsilon(P)}{d\xi(P)}\Big{|}_{P=Q}.
\end{equation}
For the differential $\varphi$ this gives (recall that $y_0=u(Q_0)$)
\begin{equation}
\label{phi_Q0}
\varphi(Q_0) = \frac{1}{\sqrt{\prod_{a_j\in B} (y_0-a_j)}}
\end{equation}
and
\begin{equation}
\label{phi_Pk}
\varphi(P_{a_k}) = \frac{2}{\sqrt{\prod_{a_j\in B\setminus\{a_k\}} (a_k-a_j)}}\,,
\end{equation}
where the evaluation at a regular point $Q_0$ is done with respect to the local parameter $\xi=u-y_0$ and the evaluation at a ramification point $P_{a_k}$ is done using $\xi=\zeta_k$ from \eqref{coordinates}. One can easily see that $\varphi(P_\infty)=0\,;$ in fact $P_\infty$ is the only zero of $\varphi$, which is thus of order $2g-2\,.$

Together with the basis of normalized differentials \eqref{hol_norm} in the space of holomorphic $1$-forms on $\mathcal H_\x$,  we introduce another basis normalized by the values at $g$ points on the surface: the $g-1$ ramification points $P_{u_j}$ corresponding to the dependant branch points and the point $Q_0$. More precisely, we define holomorphic differentials $v_1, \dots, v_g$ by the conditions
\begin{equation}
\label{vcond}
v_i(P_{u_j}) =\delta _{ij}, \
\qquad
v_i(Q_0)=\delta _{ig},  \quad\text{with} \;\;  1\leq i \leq g,\;\;1\leq j \leq g-1.
\end{equation}
Here the evaluation of the differentials is done as introduced in \eqref{evaluation}. These differentials admit an explicit description in terms of the variable $u$ as follows:
\begin{align}
\label{v}
v_i(P)&=\frac{\varphi (P)\prod_{\alpha=1, \alpha\neq i}^{g-1} (u-u_\alpha)(u-y_0)}{\varphi (P_{u_i})\prod_{\alpha, \alpha\neq i} (u_i-u_\alpha)(u_i-y_0)}, \qquad   i=1, \cdots , g-1,
\\
\label{vg}
v_g(P)&=\frac{\varphi (P) \prod_{\alpha=1}^{g-1} (u-u_\alpha)}{\varphi (Q_0)\prod_{\alpha=1}^{g-1} (y_0-u_\alpha)}\,.
\end{align}
Note that the zeros of $v_j$ at ramification points are of second order and  $v_1, \dots, v_{g-1}$ vanish also at $Q_0^*.$
\end{paragraph}
\begin{paragraph}{Meromorphic differentials} The fundamental tool for our work is the Riemann bidifferential $W(P,Q)$ with $P,Q\in\mathcal H_\x$, which can be defined as the unique bidifferential on $\mathcal H_\x$ having the following three properties:
\begin{itemize}
\item Symmetry: $W(P,Q) = W(Q,P);$
\item No singularity except for a second order pole along the diagonal $P=Q$ with biresidue $1$: for $\xi$ being a local parameter near $P=Q$, the bidifferential has the following local expansion:
\begin{equation*}
W(P,Q) \underset{P\sim Q}{=} \left( \frac{1}{(\xi(P) - \xi(Q))^2}  + {\cal O}(1) \right)d\xi(P) d\xi(Q);
\end{equation*}
\item Normalization by vanishing of all $\a$-periods: $\oint_{\a_k} W(P,Q) = 0.$
\end{itemize}
Due to the symmetry,  the above integral can be computed with respect to either $P$ or $Q$.
Clearly, $W$ depends on the choice of a canonical homology basis. As a consequence of this definition we have: $\oint_{\b_k}W(P,Q) = 2\pi{\rm i}\,\omega_k(P).$

The Riemann bidifferential admits a rather explicit representation in terms of theta-functions. We use a different approach working in terms of the coordinates $u$ and $v$ of the algebraic curve. In these coordinates, it is difficult to write a satisfactory expression for $W(P,Q)$ with $P$ and $Q$ being arbitrary varying points on the surface. However, evaluating $W$ with respect to $Q$ at specific points, we obtain a $1$-form on the surface with a second order pole at the point in question normalized by vanishing of all $\a$-periods. For the $1$-forms obtained in this way, their singularity structure allows us to write expressions in terms of the coordinates $u$, $v$ and some normalizing constants. More precisely, here are such expressions for the differentials of the second kind $W(P, \Pinfty)$ and $W(P,\Paj)$.

For $P=(u,v)$ being a point on the curve and denoting by $I_k$  the normalizing constant given by $I_k = \oint_{\a_k}u(P)^g\varphi(P)\,,$ we have
\begin{equation}
\label{Wxinfty}
W(P, \Pinfty) =  -\frac{u(P)^g\varphi(P)}{2} + \frac{1}{2} \sum_{k=1}^g I_k\omega_k(P)\,.
\end{equation}
Similarly, for the constants $\beta^{(j)}_k$ defined by vanishing of all $\a$-periods of the right hand side, we have
\begin{equation}
\label{Waj}
W(P,\Paj)=\frac{1}{u-a_j} \frac{\varphi(P)}{\varphi(\Paj)} - \sum_{k=1}^g\beta^{(j)}_k\omega_k(P)\,.
\end{equation}
We will use \eqref{Waj} for the ramification points with $a_j\in\{0,1,x_1,\dots, x_g\}\,.$ For the points $\Pualpha$ corresponding to dependent branch points $u_1, \dots, u_{g-1}$, we will need a similar expression in terms of the second basis of holomorphic differentials $v_1,
\dots, v_g$ given by \eqref{v} and \eqref{vg}:
\begin{equation}
\label{Wualpha}
W(P,\Pualpha)=\frac{1}{u-u_\alpha} \frac{\varphi(P)}{\varphi(\Pualpha)} - \sum_{k=1}^g\gamma_k^{(\alpha)} v_k(P)\,.
\end{equation}
Here again $\gamma_k^{(\alpha)}$ are normalizing constants such that the $\a$-periods of \eqref{Wualpha} vanish. Note that the ``constants'' $I_k\,,\;\beta^{(j)}_k$ and $\gamma_k^{(\alpha)}$ depend on the branch points $x_1, \dots, x_g\,.$

Let us now introduce the following differential of the third kind which will allow us to write a solution to the constrained Schlesinger system. Let us write $\Omega_{Q_0, Q^*_0}$ for the  differential of the third kind with simple poles at $Q_0$ and $Q_0^*$ with residues $1$ and $-1$, respectively, normalized by the vanishing of all of its $\a$-periods. Such a differential exist for any regular point $Q_0$  of the surface, for us $Q_0$ is the point $Q_0(\x)\in\mathcal H_\x$ defined by \eqref{q0}. This differential can be constructed in terms of the bidifferential $W$ as the integral of $W(P, Q)$ from $Q_0^*$ to $Q_0\,,$ as can be seen by comparing \eqref{Omega} and \eqref{Omega_W} below. Using the differential $\Omega_{Q_0, Q^*_0}$ and the column vector $\omega$ of the holomorphic normalized differentials as well as the constant vector $\c_2$ fixed in \eqref{q0}, we define:
\begin{equation}
\label{Omega}
\Omega(P):=\Omega_{Q_0, Q^*_0}(P) - 4 \pi\i \c_2^t\omega(P)\,.
\end{equation}
Thus  $\Omega$ is the differential of the third kind normalized by the condition $\oint_{\a_j}\Omega=-4\pi\i \c_{2j}$. One can see that its $\b$-periods are $\oint_{\b_j}\Omega=4\pi\i \c_{1j}$ where $\c_1$ is the constant vector from \eqref{q0}. Here $\c_{1j}$ and $\c_{2j}$ are the $j$th components of $\c_1$ and $\c_2$. We can also write $\Omega$ using the Riemann bidifferential:
\begin{equation}
\label{Omega_W}
\Omega(P)=\int_{Q^*_0}^{Q_0} W(P, Q) - 4 \pi\i \c_2^t\omega(P)\,.
\end{equation}
Similarly to \eqref{Wualpha}, we can write an expression for $\Omega$ in the coordinates $u, v$ of the curve as follows
\begin{equation}
\label{Omega-delta}
\Omega(P)=\frac{\varphi(P)}{\varphi(Q_0)(u-y_0)} + \sum_{j=1}^g \delta_jv_j(P)\,,
\end{equation}
where $\delta_j\in\mathbb C$ are normalizing constants and $v_j$ are the holomorphic differentials \eqref{v}.
\end{paragraph}

As is easy to see, all differentials in this section defined on the curves of the $\mathcal T$-family $T:\mathcal H \to X^0$ form well-defined objects over the whole family. Namely, the holomorphic differentials are sections of a vector bundle over $X^0$ whose fiber at $\x\in X^0$ is the space of all holomorphic differentials on $\mathcal H_\x$. This applies to the differentials $v_j$ since $Q_0$ is a section of the $\mathcal T$-family. Analogously, this $\mathcal T$-family induces a vector bundle $V_{Q_0Q_0^*} \to X^0$ whose fiber is the vector space of all the differentials having poles along the sections $Q_0$ and $Q_0^*$ of $T.$ The differentials $\Omega$ form a section of this bundle. We keep the same notation for the differentials and for the corresponding sections and do not specify the $\x$-dependence of the differentials in what follows, assuming, for example,  that $\Omega=\Omega(\x)\,.$

\section{Theorem 2: Solution to constrained Schelsinger system.}
\label{sect_solution}

Here we give a solution to system \eqref{constrained} using the above $\mathcal T$-family $(T, P_\infty, Q_0)$ of hyperelliptic curves \eqref{TMobius}. The solution is given for the case of $g$ being the genus of the curves and $K=g\,.$ The independent variables are the elements of the set $X^0$ of the independently varying branch points and the functions $u_m$ are given by the dependent branch points of the curves of the family.  All functions defined in this section depend naturally on $\x\in X^0$ although we do not reflect this dependence in our notation.

Consider the following $2\times 2$ traceless matrix
\begin{equation}
\label{A}
A(u)=\left(\begin{array}{cc}A^{11} & A^{12} \\A^{21} & -A^{11}\end{array}\right) =\sum_{a_j\in B} \frac{A_{a_j}}{u-a_j}
\end{equation}
where, using the differentials from Section \ref{sect_differentials},  we define
\begin{equation}
\label{A12_g}
A^{12}(u) := \frac{t}{2}\frac{\Omega(P)}{\varphi(P)} \frac{(u-y_0)^g}{v^2}
\end{equation}
with some arbitrary complex parameter $t$. Note that $A^{12}(u)$ is a well-defined function of $u\in\mathbb CP^1$ with simple poles in the set $B=\{0,1, x_1, \dots, x_g, u_1, \dots, u_{g-1}\}\subset\mathbb C$ and thus for any $a_j\in B$ we have
\begin{equation}
\label{A12_aj_g}
A^{12}_{a_j} =\underset{u=a_j}{\rm res} A^{12}(u) du =  \frac{t}{4}\Omega(P_{a_j})\varphi(P_{a_j})(a_j-y_0)^g\,,
\end{equation}
where the evaluation of $\Omega$ and $\phi$ at a ramification point $\Paj$ is done according to \eqref{evaluation}.
The sum of $A^{12}_{a_j}$ vanishes as a sum of residues of a differential on a compact surface:
\begin{equation}
\label{sum_A12}
\sum_{a_j\in B}A^{12}_{a_j} =0\,.
\end{equation}

Let us also introduce
\begin{equation}
\label{betas_g}
\beta_{a_j} = A_{a_j}^{12}\left(\frac{1}{(g-1)!}\frac{\partial^{g-1}}{\partial y_0^{g-1}}\left\{\frac{1}{\varphi(Q_0)(a_j-y_0)}\right\} - \frac{g}{2}\Omega(P_\infty)\right)\,.
\end{equation}

\begin{theorem}
\label{thm_main_g}
Let $(T, P_\infty, Q_0)$ be the $\mathcal T$-family of curves $\mathcal H_\x$ defined by equation \eqref{TMobius} and
 $\c_1, \c_2\in\mathbb Q^g$ the associated constant vectors from \eqref{q0}.
 As before, the family is defined over a parameter set $X^0=\{x_1, \dots, x_g\}$ and a canonical homology bases $\{\a_1, \dots, \a_g; \b_1, \dots, \b_g\}$  are chosen consistently for all curves such that the projections of the basic cycles on the $u$-sphere are independent of $\mathbf x \in X^0$. Let $B$ stand for the set of branch points, $ B=\{0,1,x_1, \dots, x_g; u_1, \dots, u_{g-1}\}$.
Let  $y_0=y_0(\mathbf x)=u(Q_0(\mathbf x))$ and  $\varphi$ and $\Omega$ be as defined in Section \ref{sect_differentials}, the sections of the appropriate vector bundles over $X^0$ induced by  the $\mathcal T$-family. Then the dependent branch points $u_1, \dots, u_{g-1}$ of the curves $\mathcal H_\x$ and the traceless matrices
\begin{equation*}
A_{a_j}=\left(\begin{array}{cc}\Aaj^{11} & \Aaj^{12} \\ \Aaj^{21} & -\Aaj^{11}\end{array}\right) \qquad \mbox{with} \quad a_j\in B
\end{equation*}
where
\begin{align*}
&A^{12}_{a_j} =  \frac{t}{4}\Omega(P_{a_j})\varphi(P_{a_j})(a_j-y_0)^g\;,
\\
& A^{11}_{a_j} = -\frac{1}{4} - \frac{g}{2t}\beta_{a_j}\;,
\\
& A^{21}_{a_j} = -\frac{g \betaaj}{4t^2}\,\frac{t + g\beta_{a_j}}{A^{12}_{a_j}}= \frac{\frac{1}{16} - (\Aaj^{11})^2}{\Aaj^{12}}\;,
\end{align*}
with the quantities $\betaaj$ being defined by \eqref{betas_g}, and $t\in\mathbb C$ being an arbitrary constant,  satisfy the constrained Schlesinger system \eqref{constrained} with respect to the variables $\x=(x_1, \dots, x_g)\in X^0$.

Moreover, derivatives of $u_k$ and $y_0$ with respect to the variables $\x=(x_1, \dots, x_g)\in X^0$ are given by
\begin{equation}
\label{umder_A}
 \frac{\partial u_m}{\partial x_i}=-\frac{\Axi^{12}}{\Aum^{12}} \frac{(u_m-y_0)^{g-1}\prod_{\alpha\neq m}(x_i-u_\alpha)}{(x_i-y_0)^{g-1}\prod_{\alpha\neq m}(u_m-u_\alpha)}\,,
\end{equation}
and
\begin{equation}
\label{dy0}
\frac{\partial y_0}{\partial x_i}=-\frac {A_{x_i}^{12}}t \frac {\prod_{\alpha=1}^{g-1}(x_i-u_{\alpha})}{\varphi (Q_0)(x_i-y_0)^g \prod _{\alpha=1}^{g-1}(y_0-u_{\alpha}) }\,.
\end{equation}

\end{theorem}

This theorem  is
proved in Section \ref{sect_proof}  except for formulas \eqref{umder_A} and \eqref{dy0} which are proved in Section \ref{sect_variation_pts}. Note that the presence of the arbitrary parameter $t$ in Theorem \ref{thm_main_g} reflects the invariance of the constrained Schlesinger system \eqref{constrained} under the following simultaneous, for all $a_j\in B,$ rescaling of the anti-diagonal elements of the matrices: $\Aaj^{12}\mapsto t\Aaj^{12} $ and $\Aaj^{21}\mapsto \frac{1}{t}\Aaj^{21}$.

\begin{remark}\label{rmk_c1c2_general}
{\rm
Theorem \ref{thm_main_g} and its proof remain valid in a more general case, namely if we assume that the family of curves
\eqref{TMobius} is a $\mathcal T$-family induced by the isoharmonic deformations from Section \ref{sect_isoharmonic}. In this case we have  $\c_1, \c_2\in\mathbb R^g\,.$ Moreover, Theorem \ref{thm_main_g} is valid if \eqref{TMobius} is any $\mathcal T$-family of hyperelliptic curves in the sense of Definition \ref{def:Tfamily}.  In this case, the constant vectors $\c_1, \c_2$ may be complex as well as the variables $x_j$ and the functions $u_k\,.$
}
\end{remark}

\begin{remark}\label{rem:onegaeta}
\rm{
In the case of real hyperelliptic curves with   appropriately chosen  bases of cycles (see \cite{Bogatyrev2012}),
the differential $\Omega=\Omega_{Q_0^*, Q_0}$ (see \eqref{Omega})  is equal  to the differential $\hat \eta_{\bf x}$ from \eqref{eq:etax} and therefore
the condition $\c_2=0$ is satisfied. In this case,
 the constrained Schlesinger equations solved in  Theorem \ref{thm_main_g} govern the isoharmonic deformations. If in addition, $\hat c_1\in \mathbb Q^g$ then the obtained constrained isomonodromic deformations provide the dynamics of the hyperelliptic $\bf T$-curves and associated real Chebyshev polynomials.}
\end{remark}

\section{Combinatorics of Bell polynomials. Useful identities}
\label{sect_identities}

In this section, we collect some of the identities which will help us to prove Theorem \ref{thm_main_g}. For the purpose of this section, all the involved quantities may be regarded as defined for one fixed curve $\mathcal H_\x$ of the $\mathcal T$-family $(T, P_\infty, Q_0),$ the compactified hyperelliptic curve of equation \eqref{TMobius}.

\subsection{Bell polynomials. Coefficients $\beta_{a_j}$}
\label{sect_Bell}
Let us introduce polynomials $L$, which will be useful for working with the coefficients $\betaaj$ defined by \eqref{betas_g}:
\begin{equation}
\label{L}
L_l(z_1, \dots, z_l) = \sum_{\substack{p_1+2p_2+\dots +lp_l=l \\ p_1,\dots, p_l\geq 0}} \frac{(-1)^{\sum_{k=1}^l p_k}\,l!\, z_1^{p_1} z_2^{p_2}\dots z_l^{p_l}}{2^{\sum_{k=1}^l p_k} \prod_{k=1}^l p_k! \prod_{k=1}^l k^{p_k}}\,.
\end{equation}
For example, we have
\begin{equation*}
L_0=1, \qquad L_1(z_1) = -\frac{z_1}{2}, \qquad L_2(z_1, z_2) = \frac{z_1^2}{4} - \frac{z_2}{2},
\qquad L_3(z_1, z_2, z_3) = -\frac{z_1^3}{8} + \frac{3}{4}z_1z_2 -z_3\,.
\end{equation*}
The polynomials $L_l$  are related to the {\it complete exponential Bell polynomials} by a simple change of variables. Indeed, the $l$-th {\it complete exponential Bell polynomial} is given by
\begin{equation*}
B_l(y_1, \dots, y_l)=\sum_{k=1}^l \sum_{\substack{p_1+p_2+\dots+p_{l-k+1}=k\\ p_1+2p_2+\dots (l-k+1)p_{l-k+1}=n}} \frac{n!}{p_1!p_2!\cdots p_{l-k+1}!}
\left( \frac{y_1}{1!}\right)^{p_1}\left( \frac{y_2}{2!}\right)^{p_2}\cdots \left( \frac{y_{l-k+1}}{(l-k+1)!}\right)^{p_{l-k+1}}
\end{equation*}
and we have
\begin{equation*}
B_l(y_1, \dots, y_l)= L_l\left(-\frac{z_1}{2\cdot 0!}, -\frac{z_2}{2\cdot1!}, \dots, -\frac{z_l}{2\cdot(l-1)!}\right)\,.
\end{equation*}

As a corollary of a similar relation for the Bell polynomials, our polynomials satisfy the binomial relation, which can also be derived from Proposition \ref{prop_Ll} and the Leibniz rule for differentiation:
\begin{equation}
\label{L-Leibniz}
\sum_{k=0}^n{n\choose k}L_k(z_1, \dots, z_k) L_{n-k}(y_1, \dots, y_{n-k})=L_n(z_1+y_1, \dots, z_n+y_n)\,.
\end{equation}
We will use the following sums over all the branch points $a_j$
\begin{equation*}
\Sigma_k:=\sum_{a_j\in B} \frac{1}{(a_j-y_0)^k}
\end{equation*}
for integer values of $k\,.$ To shorten the expressions, we will write
\begin{equation}
\label{Ll-notation}
L_l:=L_l(\Sigma_1, \dots, \Sigma_l)\,.
\end{equation}

To understand better the quantities $\beta_{a_j}$, we need the following result.
\begin{proposition}
\label{prop_Ll}
Let the polynomials $L_l$ be as above and the value $\varphi(Q_0)$ be defined by \eqref{phi_Q0}. Then for any integer $n$ and any integer $l\geq 0$
\begin{equation}
\label{Ll}
\frac{\partial^{l}}{\partial y_0^l}\left\{\frac{1}{\varphi(Q_0)}\right\} = \frac{L_l}{\varphi(Q_0)}
\qquad\mbox{and}\qquad
\frac{\partial^{l}\varphi^n(Q_0)}{\partial y_0^l} = \varphi^n(Q_0) L_l(-n\Sigma_1, \dots, -n\Sigma_l)\,.
\end{equation}
\end{proposition}
{\it Proof.} We prove the first equality, the proof for the second one being entirely similar.
As is easy to see, the derivatives in question can be written in the form of the right hand side of \eqref{Ll} with some polynomial in $\Sigma_1, \dots, \Sigma_l$ of the following form
\begin{equation*}
\frac{\partial^{l}}{\partial y_0^l}\left\{\frac{1}{\varphi(Q_0)}\right\} = \frac{1}{\varphi(Q_0)} \sum_{\substack{p_1+2p_2+\dots +lp_l=l \\ p_1,\dots, p_l\geq 0}} C_{p_1 p_2\dots p_l} \Sigma_1^{p_1}\Sigma_2^{p_2}\dots \Sigma_l^{p_l}\,.
\end{equation*}
It remains thus to find the coefficients $C_{p_1 p_2\dots p_l}\,.$  By examining the derivatives, we notice that these coefficients are solutions of the recursion relation
\begin{equation*}
C_{p_1 p_2\dots p_l} = -\frac{1}{2}C_{p_1-1, p_2\dots p_l} + \sum_{n=1}^{l-1} n(p_n+1) C_{p_1\dots p_n+1,p_{n+1}-1, p_{n+2}\dots p_l}\,.
\end{equation*}
It is straightforward to verify that the coefficients of the polynomial \eqref{L} satisfy this recursion with the initial value $C_0 = 1\,.$ Note that in this notation, one can omit the indices of a coefficient $C$ which are equal to zero and are placed at the end of the string of indices, that is, for example,  $C_1=C_{1,0,0,0}\,.$
$\Box$

Let us thus define
\begin{equation}
\label{Coefficients}
C_{p_1 p_2\dots p_l} = \frac{(-1)^{\sum_{k=1}^l p_k}\,l!\, }{2^{\sum_{k=1}^l p_k} \prod_{k=1}^l p_k! \prod_{k=1}^l k^{p_k}}\,.
\end{equation}
\begin{corollary}
\label{corollary_betas}
For any $a_j\in B$, the quantities $\beta_{a_j}$ defined by \eqref{betas_g} can be rewritten as
\begin{equation}
\label{betas}
\beta_{a_j} = A_{a_j}^{12}\left(\frac{1}{\varphi(Q_0)} \sum_{l=0}^{g-1} \frac{L_l}{l!(a_j-y_0)^{g-l}} - \frac{g}{2}\Omega(P_\infty)\right)\,.
\end{equation}
\end{corollary}

{\it Proof.} This follows from Proposition \ref{prop_Ll} by applying the Leibniz rule for derivative of a product.
$\Box$

In what follows, we will often deal with quantities of the form
\begin{equation}
\label{difference_betas_1}
\frac{\betaai}{\Aai^{12}}-\frac{\betaaj}{\Aaj^{12}} = \frac{1}{\varphi(Q_0)} \sum_{l=0}^{g-1} \left( \frac{L_l}{l!(a_i-y_0)^{g-l}} - \frac{L_l}{l!(a_j-y_0)^{g-l}} \right)\,,
\end{equation}
which can be conveniently written as
\begin{equation}
\label{difference_betas_2}
\frac{\betaai}{\Aai^{12}}-\frac{\betaaj}{\Aaj^{12}} = \frac{a_j-a_i}{\varphi(Q_0)} \sum_{l=0}^{g-1}  \frac{L_l}{l!}\sum_{k=0}^{g-l-1}\frac{1}{(a_i-y_0)^{k+1}(a_j-y_0)^{g-l-k}} \,.
\end{equation}

\subsection{Identities for sums over branch points}
\label{sect_sumidentities}
Note that
\begin{equation}
\label{sumA12_g}
\sum_{a_j\in B}A^{12}_{a_j} =0
\end{equation}
as a sum of residues of the differential $A^{12}(u)du$, see \eqref{A12_g}. Moreover, we have the following two lemmas.
\begin{lemma}
\label{lemma_residues} Let $B$ \eqref{B} be the set of finite branch points of the compactified hyperelliptic curve $\mathcal H_\x$ of equation \eqref{TMobius} and the functions $\Aaj^{12}$ be defined by \eqref{A12_aj_g}. Let the point $Q_0$ on $\mathcal H_\x$ be defined by \eqref{q0} and $\Omega$ and $\phi$ be differentials \eqref{Omega} and \eqref{phi} on $\mathcal H_\x$, respectively.
Then the following identities hold on
\begin{eqnarray}
\label{res0}
&&\sum_{a_j\in B} \frac{\Aaj^{12}}{(a_j-y_0)^s}=0\qquad\mbox{for}\quad 0\leq s\leq g-1\,;
\\
\label{res1}
&&\sum_{a_j\in B} \frac{\Aaj^{12}}{(a_j-y_0)^g}=-t\varphi(Q_0)\,;
\\
\label{res2}
&&\sum_{a_j\in B} \frac{\Aaj^{12}}{(a_j-y_0)^s}=-t\,\underset{P=Q_0}{\rm res}\frac{\Omega(P)\varphi(P)}{(u-y_0)^{s-g}du}\qquad\mbox{for}\quad s\geq g+1\,.
\end{eqnarray}
\end{lemma}

{\it Proof.} These identities follow from the vanishing of the sum of residues of the differential $\frac{t\Omega(P)\varphi(P)}{2(u-y_0)^{s-g}du}$ on the compact surface $\mathcal H_\x$.
$\Box$

\begin{lemma}
\label{lemma_residues_g}
In the situation of Lemma \ref{lemma_residues}, the following relations hold
\begin{eqnarray}
\label{res3}
&& \frac{1}{2}\sum_{a_j\neq u_m} \frac{\Omega(\Paj)\varphi(\Paj)}{a_j-u_m} +\underset{P=\Pum}{\rm res} \frac{\Omega(P) \varphi(P)}{(u-u_m)du}+ 2\frac{\varphi(Q_0)}{y_0-u_m}=0\,;
\\
\label{res4}
&& \frac{1}{2}\sum_{a_j\neq u_m} \Omega(\Paj)W(\Paj, \Pum) +\underset{P=\Pum}{\rm res} \frac{\Omega(P) W(P, \Pum)}{du}+ 2W(Q_0, \Pum)=0\,.
\end{eqnarray}
\end{lemma}
{\it Proof.}
 Equality  \eqref{res3} is obtained  as the sum of residues of the differential $\frac{\Omega(P) \varphi(P)}{(u-u_m)du}$ on $\mathcal H_\x$.  The vanishing of the sum of residues of the differential $\frac{\Omega(P) W(P, \Pum)}{du}$ gives \eqref{res4}. $\Box$

Note also that, due to \eqref{Wualpha}, and defining relations \eqref{vcond} for the differentials $v_j$,  we have
\begin{equation*}
\frac{1}{\varphi(P_{u_m})}\underset{P=P_{u_m}}{\rm res} \frac{\Omega(P) \varphi(P)}{(u-u_m)du} -
\underset{P=P_{u_m}}{\rm res} \frac{\Omega(P) W(P, P_{u_m})}{du} =
\underset{P=P_{u_m}}{\rm res} \frac{\Omega(P)}{du} \sum_{k=1}^g \gamma_k^{(m)} v_k(P)
=\frac{1}{2}\Omega(\Pum)\gamma_m^{(m)}\,.
\end{equation*}
This together with \eqref{res3} allows us to obtain for the  residue in \eqref{res4}:
\begin{equation}
\label{res-diff_g}
\underset{P=P_{u_m}}{\rm res} \frac{\Omega(P) W(P, P_{u_m})}{du}
=
-\frac{1}{\varphi(P_{u_m})}\left(\frac{1}{2}\sum_{a_j\neq u_m} \frac{\Omega(\Paj)\varphi(\Paj)}{a_j-u_m} + 2\frac{\varphi(Q_0)}{y_0-u_m}\right) -
\frac{1}{2}\Omega(\Pum)\gamma_m^{(m)}\,.
\end{equation}

\begin{lemma} Let the notation be as in Lemma \ref{lemma_residues} and  let $k$ be an integer $0\leq k< g\,.$ The following identity holds
\label{lemma_simplify}
\begin{equation*}
 \sum_{\substack{a_j\in B\\a_j \neq u_m}} \frac {A_{a_j}^{12}}{(a_j-y_0)^{g-k}(a_j-u_m) } = \frac 1{(u_m-y_0)^{g-k}}  \sum_{\substack{a_j\in B\\a_j \neq u_m}}\frac {A_{a_j}^{12}}{a_j-u_m} +\frac {(g-k)A_{u_m}^{12}}{(u_m-y_0)^{g-k+1}} + \delta_{k,0}\frac{t \varphi(Q_0)}{u_m-y_0}\,.
\,
\end{equation*}

\end{lemma}
{\it Proof.} Taking a factor of $(a_j-y_0)$ in the denominator and splitting $\frac {1}{(a_j-y_0)(a_j-u_m) }$ into a sum of partial fractions, we can then use the vanishing of the sum of residues \eqref{res0} if $k>0$ or \eqref{res1} if $k=0$ to obtain
\begin{equation*}
 \sum_{a_j \neq u_m} \frac {A_{a_j}^{12}}{(a_j-y_0)^{g-k}(a_j-u_m) } = \frac 1{u_m-y_0} \sum_{a_j \neq u_m} \frac {A_{a_j}^{12}}{(a_j-y_0)^{g-k-1}(a_j-u_m) } +\frac {A_{u_m}^{12}}{(u_m-y_0)^{g-k+1}}+ \delta_{k,0}\frac{t \varphi(Q_0)}{u_m-y_0}.
\end{equation*}
Applying this procedure successively $g-k$ times, we prove the lemma.
$\Box$

\subsection{Some rational identities}
\label{sect_rational}

In this section we list some simple rational relations that will be useful in our calculation. In all the identities we assume $N\geq 1\,.$

\begin{lemma}
\label{lemma_rational1}
For any $x$ and $y$ not in the set $\{u_\alpha\}_{\alpha=1}^{N}\subset \mathbb C$ and $0\leq s \leq N-1$, we have
\begin{align}
\label{rat1}
&\sum_{k=1}^{N} \frac{1}{(x-u_k)\prod_{ \alpha \neq k}(u_k-u_\alpha)}=\frac{1}{\prod_{\alpha =1}^{N}(x-u_{\alpha})};
\\
\label{rat2}
&\sum_{k=1}^{N} \frac {(u_k-y)^s}{(x-u_k)\prod_{\alpha \neq k}(u_k-u_{\alpha})}=\frac {(x-y)^s}{\prod_{\alpha=1}^{N}(x-u_{\alpha})};
\\
\label{rat3}
&\sum_{k=1}^{N} \frac {(u_k-y)^{N}}{(x-u_k)\prod_{\alpha \neq k}(u_k-u_{\alpha})}=\frac {(x-y)^{N}}{\prod_{\alpha=1}^{N}(x-u_{\alpha})}-1.
\end{align}
\end{lemma}

{\it Proof.}
In each equality, the expressions in the left and the right hand side, as functions of $x$, have simple poles at $x=u_k$ with equal residues, have no other singularities, and vanish at infinity.
$\Box$

\begin{corollary}
\label{corollary_rational}
For $y$ not in the set $\{u_\alpha\}_{\alpha=1}^{N}\subset \mathbb C$ and $0\leq s\leq N-2$, we have
\begin{align}
\label{rat4}
\sum_{k=1}^{N} \frac {(u_k-y)^{N-1}}{\prod_{\alpha \neq k}(u_k-u_{\alpha})}=1;
\\
\label{rat5}
\sum_{k=1}^{N} \frac {(u_k-y)^s}{\prod_{\alpha \neq k}(u_k-u_{\alpha})}=0.
\end{align}
\end{corollary}
{\it Proof.} The first identity follows by setting $x=y$ in \eqref{rat3} and \eqref{rat5} is obtained by differentiating \eqref{rat4} with respect to $y$. $\Box$

\begin{corollary}
\label{corollary_rational2}
For $x$ and $y$ not in the set $\{u_\alpha\}_{\alpha=1}^{N}\subset \mathbb C$ and $s\in\mathbb N$, $s>0,$ we have
\begin{equation*}
\sum_{k=1}^{N} \frac {1}{(u_k-y)^s(x-u_k)\prod_{\alpha \neq k}(u_k-u_{\alpha})}
=\frac{1}{(x-y)^s\prod_{\alpha =1}^{N}(x-u_{\alpha})} - \frac{1}{(s-1)!} \frac{\partial^{s-1}}{\partial y^{s-1}} \left\{ \frac{1}{(x-y)\prod_{\alpha =1}^{N}(y-u_{\alpha})}\right\}\,.
\end{equation*}
\end{corollary}
{\it Proof.} We apply the first identity of Lemma \ref{lemma_rational1} to the left hand side represented in the form
\begin{equation*}
\frac{1}{(n-1)!} \frac{\partial^{n-1}}{\partial y^{n-1}} \sum_{k=1}^{N} \frac{1}{x-y}\left( \frac{1}{u_k-y} + \frac{1}{x-u_k}\right) \frac {1}{\prod_{\alpha \neq k}(u_k-u_{\alpha})}\,.
\end{equation*}
$\Box$

\begin{lemma}
\label{lemma_rational2}
For $y$ not in the set $\{u_\alpha\}_{\alpha=1}^{N}\subset \mathbb C$ and $r\geq 0$, we have
\begin{multline}
\label{rat-alpha2}
\sum_{\substack{k=1\\ k\neq m}}^{N}\frac{1}{(u_k-y)^{r+1}(u_k-u_m)\prod_{\beta\neq k}(u_k-u_\beta)}
=
\frac{\sum_{\beta\neq m} \frac{1}{u_m-u_\beta}}{(u_m-y)^{r+1}\prod_{k\neq m}(u_m-u_k)}
\\
+ \frac{r+1}{(u_m-y)^{r+2}\prod_{k\neq m}(u_m-u_k)}
+\frac{1}{r!} \frac{\partial^r}{\partial y^r}\left\{ \frac{1}{(u_m-y)\prod_{k=1}^{N}(y-u_k)}  \right\}.
\end{multline}
\end{lemma}

{\it Proof.} Let us first prove the equality for $r=0.$ In this case, the left hand side can be written as
\begin{equation*}
\frac{1}{u_m-y}\left( \sum_{\substack{k=1\\ k\neq m}}^{N} \frac{1}{(u_k-u_m)\prod_{\beta\neq k}(u_k-u_\beta)}
+ \sum_{\substack{k=1\\ k\neq m}}^{N} \frac{1}{(y-u_k)\prod_{\beta\neq k}(u_k-u_\beta)}\right).
\end{equation*}
The result then follows by using \eqref{rat1} from Lemma \ref{lemma_rational1} in the second sum. The identity for the other values of $r$ can be obtained from the case $r=0$ by the $r$-fold differentiation with respect to $y$. $\Box$

\begin{lemma}
\label{lemma_rational3}
For $y$ not in the set $\{u_\alpha\}_{\alpha=1}^{N}\subset \mathbb C$ and $0\leq s\leq N-1$, we have
\begin{align}
\label{rat-alpha3}
&\sum_{\substack{k=1\\ k\neq m}}^{N}\frac{(u_k-y)^s}{(u_k-u_m)\prod_{\beta\neq k}(u_k-u_\beta)} =-\frac{s(u_m-y)^{s-1}}{\prod_{k\neq m}(u_m-u_k)} + \frac{(u_m-y)^s}{\prod_{k\neq m}(u_m-u_k)} \sum_{\substack{\beta=1\\\beta\neq m}}^N \frac{1}{u_m-u_\beta};
\\
\label{rat-alpha4}
&\sum_{\substack{k=1\\ k\neq m}}^{N}\frac{(u_k-y)^N}{(u_k-u_m)\prod_{\beta\neq k}(u_k-u_\beta)} =1-\frac{N(u_m-y)^{N-1}}{\prod_{k\neq m}(u_m-u_k)} + \frac{(u_m-y)^N}{\prod_{k\neq m}(u_m-u_k)} \sum_{\substack{\beta=1\\\beta\neq m}}^N \frac{1}{u_m-u_\beta}.
\end{align}
\end{lemma}

{\it Proof.} We prove \eqref{rat-alpha3} by representing its left hand side as
\begin{equation*}
- \frac{\partial}{\partial u_m}   \sum_{\substack{k=1\\ k\neq m}}^{N} \frac{(u_k-y)^s}{(u_m-u_k)\prod_{\substack{\beta\neq k\\\beta\neq m}}(u_k-u_\beta)}
=
- \frac{\partial}{\partial u_m} \left\{   \frac{(u_m-y)^s}{\prod_{\beta\neq m}(u_m-u_\beta)}
\right\}
\end{equation*}
where the equality is obtained by applying \eqref{rat2} in the case $s\leq N-2$ and \eqref{rat3} in the case $s=N-1$ with $x=u_m$ and the set of $N-1$ distinct points $\{u_\alpha\}$ from which $u_m$ is removed. To prove \eqref{rat-alpha4}, we first rewrite its left hand side in the form
\begin{multline*}
\sum_{\substack{k=1\\ k\neq m}}^{N}\frac{(u_k-y)^{N-1}}{(u_k-u_m)\prod_{\substack{\beta\neq k\\ \beta \neq m}}(u_k-u_\beta)}
+
(u_m-y)\sum_{\substack{k=1\\ k\neq m}}^{N}\frac{(u_k-y)^{N-1}}{(u_k-u_m)^2\prod_{\substack{\beta\neq k\\ \beta \neq m}}(u_k-u_\beta)}
\\
= 1 - \frac{(u_m-y)^{N-1}}{\prod_{\beta \neq m}(u_m-u_\beta)}
- (u_m-y) \frac{\partial}{\partial u_m} \sum_{\substack{k=1\\ k\neq m}}^{N}\frac{(u_k-y)^{N-1}}{(u_m-u_k)\prod_{\substack{\beta\neq k\\ \beta \neq m}}(u_k-u_\beta)}
\end{multline*}
where the equality is obtained by using \eqref{rat3} as above in the first sum. It remains to use \eqref{rat3} in the second sum as well and then compute the derivative.
$\Box$

\section{Variational formulas}
\label{sect_variational}

In this section, we study the dependence of the quantities related to the $\mathcal T$-family $(T, P_\infty, Q_0)$ of hyperelliptic curves on the point $\x=(x_1, \dots, x_g)$ in the  parameter set $X^0\,.$ To this end, we use the Rauch variational formulas from \cite{Fay92} in the form written in \cite{KokoKoro}. These formulas allow us to find derivatives of the quantities involved in the statement of Theorem \ref{thm_main_g} with respect to the independetly varying branch points $x_1, \dots, x_g$ of $\mathcal H_\x\,.$ From now on, we adopt a slightly different terminology for the families of hyperelliptic curves,  the terminology used to describe the Rauch variation. Namely, varying a branch point of a hyperelliptic curve $\mathcal H_\x$ results in the variation of the complex structure given by the local charts \eqref{coordinates} on the associated topological surface. Thus our $\mathcal T$-family of hyperelliptic curves is regarded as a family of complex structures on a compact orientable surface of genus $g$. The complex structures are parameterized by branch points, thus varying the position of a branch point results in a variation of the complex structure. All the quantities defined on our curves depend on the complex structure and thus become functions of the branch points. The Rauch variational formulas allow us to describe this variation.

\subsection{Rauch variational formulas}
\label{sect_Rauch}

In this subsection, we assume that all finite branch points of the curves \eqref{TMobius}, the points of the set $B=\{a_j\}_{j=1}^{2g+1}$,  can vary independently of each other. Therefore we cannot use the condition \eqref{q0} for the point $Q_0$; this point will be considered simply as some regular point of the curve. We also assume that such a variation leaves the vector $(a_1, \dots, a_{2g+1})$ in a small open ball inside $\mathbb C^{2g+1}\setminus\{\Delta_{ij}\}_{i\neq j}$, where $\Delta_{ij}$ are diagonals defined by $a_i=a_j$ for  $i\neq j$  ranging through $1, \dots, 2g+1\,.$ The complex structure of the Riemann surface associated with the curve $\mathcal H_\x$ \eqref{TMobius} is defined by the local coordinates  \eqref{coordinates}. This complex structure varies under the variation of the branch points $a_j\in B$ and therefore all Abelian differentials defined on the surface vary accordingly. We can thus consider our differentials as depending on $a_1, \dots, a_{2g+1}$ and a point $P$ of the surface.

To measure the dependence of the differentials on the branch points, we use the following Rauch derivative, see \cite{Fay92, KokoKoro}, defined as derivative with respect to one of the branch points while fixing the point $P$ on a varying surface by the requirement that its $u$-coordinate stays fixed under the variation: $u(P)={\rm const}.$ We denote this derivative as follows:
\begin{equation}
\label{Rauch}
\frac{\partial^{{\rm Rauch}}}{\partial a_k} \Upsilon(P) := \frac{\partial}{\partial a_k}\Big{|}_{u(P)=const} \Upsilon(P)
\end{equation}
for an Abelian differential $\Upsilon(P)$ defined on the Riemann surface of the curve \eqref{TMobius}. In the case of the Riemann bidifferential we need to require that both $u(P)$ and $u(Q)$ stay fixed. We have the following Rauch variational formula for the  $W$, see \cite{KokoKoro}:
\begin{equation*}
\frac{\partial^{{\rm Rauch}}}{\partial a_k}W(P,Q) := \frac{\partial}{\partial a_k}\Big{|}_{\substack{u(P)=const\\u(Q)=const}}W(P,Q) = \frac{1}{2}W(P, \Paj)W(\Paj, Q).
\end{equation*}
This variation of the Riemann bidifferential is a master-formula which implies the following Rauch formulas, via $\omega_j(P)=\oint_{\b_j}W(P,Q)/(2\pi\i)$ and $\mathbb B_{jk} = \oint_{\b_k}\omega_j\,:$
\begin{equation}
\label{RauchB}
\frac{\partial^{{\rm Rauch}} \omega_j(P)}{\partial a_k} = \frac{1}{2}\omega_j(\Pak) W(P, \Pak),
\qquad\qquad
\frac{\partial^{{\rm Rauch}} \mathbb B_{ij}}{\partial a_k} = \frac{\pi \i}{2}\omega_j(\Pak)\omega_i(\Pak)\,.
\end{equation}
Assuming now that $Q_0$ is a point on the Riemann surface $\mathcal H_\x$ with a fixed projection $y_0$ on the $u$-sphere, $u(Q_0)=y_0$, independent of the branch points, and using definition \eqref{Omega_W} of the differential $\Omega$, we derive the following Rauch variational formula
\begin{equation}
\label{RauchOmega}
\frac{\partial^{{\rm Rauch}} \Omega(P)}{\partial a_k} = \frac{1}{2}\Omega(\Pak) W(P, \Pak)\,.
\end{equation}
Evaluating \eqref{RauchOmega} at $\Paj$ for $k\neq j$ we have
\begin{equation}
\label{RauchOmegaaj}
\frac{\partial^{{\rm Rauch}} \Omega(\Paj)}{\partial a_k} = \frac{1}{2}\Omega(\Pak) W(\Paj, \Pak)\,.
\end{equation}

In what follows, we will also need a formula that allows us to differentiate $\Omega(\Pak)$ with respect to the branch point coinciding with the argument, that is with respect to $a_k.$ As the right hand side of \eqref{RauchOmegaaj} is not defined for $a_k=a_j$,  we need the following lemma in this case.
\begin{lemma}\label{lemma_epsilon_g}
Let $a_k\in B$ be an arbitrary branch point of the curve $\mathcal H_\x$ \eqref{TMobius} and $Q_0, Q_0^*\in\mathcal H_\x$ be two regular points $u$-coordinate of which, denoted by $y_0=u(Q_0)=u(Q_0^*)$, is fixed and independent of the branch points. Let $\Omega$ be the differential of the third kind  having simple poles at the points $Q_0$ and $Q_0^*$ defined by \eqref{Omega} in which $\c_2\in\mathbb C^g$ is a constant  vector and $\omega$ is a vector of holomorphic  normalized differentials. The following variational formula holds:
\begin{equation*}
\frac{\partial \Omega(P_{a_k})}{\partial a_k} = - \frac{1}{2}\sum_{\substack{a_j\in B \\ a_j\neq a_k}}\Omega(\Paj)W(\Pak,\Paj) - 2W(Q_0, \Pak)
\,.
\end{equation*}
\end{lemma}
{\it Proof.} Let $\varepsilon$ be a complex number with small absolute value and let  the curve $\mathcal H_\x^\varepsilon$ be obtained from the curve $\mathcal H_\x$ by performing a shift by $\varepsilon$ in every sheet of the covering $u:\mathcal H_\x \to \mathbb C.$ More precisely, we define
\begin{equation*}
\mathcal H_\x^\varepsilon:=\{P^\varepsilon=(u+\varepsilon, v)\,|\, P=(u,v)\in \mathcal H_\x\},
\end{equation*}
an algebraic curve which we see as a two-fold covering of the $u$-plane ramified at the points $\Paj^\varepsilon:=(a_j+\varepsilon, 0)$ for ${a_j\in B}$. Denote also by $\mathcal H_\x^\varepsilon$ and $\mathcal H_\x$ the compact Riemann surfaces corresponding to the two algebraic curves. The complex structure \eqref{coordinates} on $\mathcal H_\x$ is not affected by the shift $u\mapsto u+\varepsilon$ and thus the two Riemann surfaces coincide. In other words, we have the biholomorphic map $\varepsilon: \mathcal H_\x\to \mathcal H_\x^\varepsilon$ sending a point $P$ to $P^\varepsilon$ which extends to the compact Riemann surfaces.

Let us denote by $W^\varepsilon$ the Riemann bidifferential on the surface $\mathcal H_\x^\varepsilon$. We may assume that the chosen canonical homology basis $\{\a_j, \b_j\}$ on $\mathcal H_\x$ transforms into a canonical basis $\{\a_j^\varepsilon, \b_j^\varepsilon\}$ on $\mathcal H_\x^\varepsilon$. Since $\varepsilon$ is a biholomorphic map, the pull-back of
$W^\varepsilon$ by this map to $\mathcal H_\x$ coincides with $W$ on $\mathcal H_\x$ due to the unicity of the Riemann bidifferential for a fixed canonical basis in the homology. In other words, we have
\begin{equation}
\label{Wpullback}
W(P,Q)=W^\varepsilon(P^\varepsilon, Q^\varepsilon)\,.
\end{equation}
In a similar way, we denote $\omega_j^\varepsilon = \oint_{\b_j^\varepsilon}W^\varepsilon(\cdot, Q)/(2\pi\i)$ the holomorphic differentials on $\mathcal H_\x^\varepsilon$ normalized by $\oint_{\a_k^\varepsilon}\omega_j^\varepsilon=\delta_{ij}\,.$ Note that, similarly to \eqref{Wpullback}, we have
\begin{equation}
\label{omegapullback}
\omega_j^\varepsilon(P^\varepsilon)=\omega_j(P)\,, \qquad j=1, \dots, g\,.
\end{equation}
 Let us also define $\Omega^\varepsilon$ on $\mathcal H_\x^\varepsilon$ by \eqref{Omega_W}, that is
\begin{equation*}
\Omega^\varepsilon(P)=\int_{{Q^*_0}^\varepsilon}^{{Q_0}^\varepsilon} W^\varepsilon(P, Q) - 4 \pi\i \c_2^t\omega^\varepsilon(P)
\end{equation*}
where $P\in\mathcal H_\x^\varepsilon$ and ${Q_0}^\varepsilon, {Q_0^*}^\varepsilon$ are points on $\mathcal H_\x^\varepsilon$ whose $u$-coordinate is $y_0+\varepsilon\,.$ Evaluating $\Omega^\varepsilon$ at a ramification point $\Pak^\varepsilon$ as in \eqref{evaluation}, we obtain $\Omega^\varepsilon(\Pak^\varepsilon)$ as a function of $\{a_j+\varepsilon\}$ and $y_0+\varepsilon\,:$
\begin{equation*}
\Omega^\varepsilon(\Pak^\varepsilon)=\int_{{Q^*_0}^\varepsilon}^{{Q_0}^\varepsilon} W^\varepsilon(\Pak^\varepsilon, Q) - 4 \pi\i \c_2^t\omega^\varepsilon(\Pak^\varepsilon)
\end{equation*}
which we can rewrite as
\begin{equation*}
\Omega^\varepsilon(\Pak^\varepsilon)=\int_{{Q^*_0}}^{{Q_0}} W^\varepsilon(\Pak^\varepsilon, Q^\varepsilon) - 4 \pi\i \c_2^t\omega^\varepsilon(\Pak^\varepsilon) =\Omega(\Pak)\,,
\end{equation*}
where in the last equality we used \eqref{Wpullback} and \eqref{omegapullback}. From the equality $\Omega^\varepsilon(\Pak^\varepsilon)=\Omega(\Pak)$ we deduce
\begin{equation*}
\frac{d}{d\varepsilon}\Omega^\varepsilon(\Pak^\varepsilon)=0\,.
\end{equation*}
On the other hand, since $\Omega^\varepsilon(\Pak^\varepsilon)$ is a function of $\{a_j+\varepsilon\}, \; y_0+\varepsilon\,$, we have
\begin{equation*}
0=\frac{d}{d\varepsilon}\Big{|}_{\varepsilon=0}\Omega^\varepsilon(\Pak^\varepsilon)=\left(\sum_{a_j\in B} \frac{\partial\Omega^\varepsilon(\Pak^\varepsilon)}{\partial (a_j+\varepsilon)} + \frac{\partial\Omega^\varepsilon(\Pak^\varepsilon)}{\partial (y_0+\varepsilon)}\right)\Big{|}_{\varepsilon=0}=\sum_{a_j\in B} \frac{\partial\Omega(\Pak)}{\partial a_j} + \frac{\partial\Omega(\Pak)}{\partial y_0}\,.
\end{equation*}
From here we deduce the expression for the derivative of $\Omega(\Pak)$ with respect to $a_k$. Noting that the partial derivative $\frac{\partial\Omega(\Paj)}{\partial a_k}$ coincides with the Rauch derivative \eqref{Rauch}, \eqref{RauchOmegaaj}, we obtain
\begin{equation*}
\frac{\partial\Omega(\Pak)}{\partial a_k}=-\sum_{\substack{a_j\in B \\ a_j\neq a_k}} \frac{\partial^{{\rm Rauch}}\Omega(\Pak)}{\partial a_j} - \frac{\partial\Omega(\Pak)}{\partial y_0}
\,.
\end{equation*}
From definition \eqref{Omega_W} of $\Omega$, we see that
\begin{equation*}
 \frac{\partial\Omega(\Pak)}{\partial y_0}= W(\Pak, Q_0)-W(\Pak, Q_0^*)=2W(\Pak, Q_0)\,,
\end{equation*}
where in the last equality we used the anti-invariance of $W(\Pak, \cdot)$ with respect to the hyperelliptic involution, which follows, for example, from \eqref{Waj}. With this and \eqref{RauchOmegaaj} for the Rauch derivatives, we prove the lemma.
$\Box$

\subsection{Variational formulas on the ${\mathcal T}$-families of curves}
\label{sect_Tvariation}

Here we go back to regarding $\mathcal H_\x$ \eqref{TMobius} as a curve  of the ${\mathcal T}$-family $(T, P_\infty, Q_0)$ where the section $Q_0$ satisfies \eqref{q0} and thus
the branch points $x_1, \dots, x_g$ vary independently while the branch points $u_1, \dots, u_{g-1}$ are functions of $x_1, \dots, x_g$. The three branch points at $0,\,1$ and $\infty$ remain fixed.

\subsubsection{Variation of dependent branch points on the ${\mathcal T}$-families of curves}
\label{sect_variation_pts}

To obtain derivatives of the dependent branch points $u_1, \dots, u_{g-1}$ and of the projection $y_0=u(Q_0)$ of the point $Q_0$ with respect to the independent branch points $x_1, \dots, x_g$, let us differentiate relation \eqref{q0} defining the section $Q_0.$
Differentiating this relation with respect to an independent branch point $x_i$ with the help of the Rauch variational formulas \eqref{Rauch}, we have
\begin{multline*}
\frac{1}{4}{\bf \omega}(p_{x_i})\int_{Q_0^*}^{Q_0}W(p_{x_i}, p)
+\frac{1}{4}\sum_{k=1}^{g-1}{\bf \omega}(p_{u_k})\int_{Q_0^*}^{Q_0}W(p_{u_k}, p) \frac{\partial u_k}{\partial x_i}
 + {\bf \omega}(Q_0) \frac{\partial y_0}{\partial x_i}
 \\
 = \pi \i{\bf \omega}(p_{x_i}){\bf \omega}^T(p_{x_i}) \c_2 + \pi \i\sum_{k=1}^{g-1}{\bf \omega}(p_{u_k}){\bf \omega}^T(p_{u_k}) \c_2\frac{\partial u_k}{\partial x_i}\,.
\end{multline*}
Here we used the anti-invariance of differential $W(\Pxi, P)$ with respect to the hyperelliptic involution, $W(\Pxi, P)=-W(\Pxi, P^*),$ which follows from \eqref{Waj}. This anti-invariance implies $\int_\Pinfty^{Q_0}W(\Pxi, P) =- \int_\Pinfty^{Q_0^*}W(\Pxi, P)\,.$ The above relation can be rewritten using the differential $\Omega$ \eqref{Omega_W} as follows:
\begin{equation}
\label{tempdiff_g}
\frac{1}{4}{\bf \omega}(P_{x_i})\Omega(P_{x_i})
+\frac{1}{4}\sum_{m=1}^{g-1}{\bf \omega}(P_{u_m})\Omega(P_{u_m}) \frac{\partial u_m}{\partial x_i}
 + {\bf \omega}(Q_0) \frac{\partial y_0}{\partial x_i} = 0\,.
\end{equation}

Recall that $\omega$ stands for the $g$-component column vector of holomorphic normalized differentials. Thus for $1\leq m \leq g-1$ and a fixed $i$, this is a linear system of equations for the $g$ unknown functions $\frac{\partial u_m}{\partial x_i}$ and $\frac{\partial y_0}{\partial x_i}$, solving which by Cramer's rule, we obtain
\begin{equation}
\label{umder}
 \frac{\partial u_m}{\partial x_i}=-\frac{1}{4}\frac{ \Omega(\Pxi)}{\Omega(\Pum)} \frac{{\rm det}\,M_m}{{\rm det}\,M}
\end{equation}
and
\begin{equation}
\label{y0der}
 \frac{\partial y_0}{\partial x_i}=-\frac{1}{4}\Omega(\Pxi) \frac{{\rm det}\,M_g}{{\rm det}\,M}\,,
\end{equation}
where $M$ stands for the $g\times g$ matrix $[\omega(P_{u_1}), \dots, \omega(P_{u_{g-1}}), \omega(Q_0)]$  and $M_m$ is the matrix $M$ with the $m$th column replaced by $\omega(\Pxi).$ Note that expressions \eqref{umder} are invariant under replacing the column vector $\omega$ of holomorphic normalized differentials $\omega_j$ \eqref{hol_norm} by a column vector whose components form any other basis in the space of holomorphic 1-forms on our surface $\surf$. Replacing the vector $\omega$ by the column vector $v=(v_1, \dots, v_g)^t$ of the differentials \eqref{v}, \eqref{vg}, the matrix $M$ gets replaced by the identity matrix and we obtained a simpler form of the variational formulas for $u_m$ and $y_0$:
\begin{equation}
\label{umder_g}
 \frac{\partial u_m}{\partial x_i}=-\frac{ \Omega(\Pxi)}{\Omega(\Pum)} v_m(\Pxi)\,,
\end{equation}
\begin{equation}
\label{y0der_g}
\frac{\partial y_0}{\partial x_i}=-\frac{1}{4}\Omega(\Pxi)v_g(\Pxi)\,.
\end{equation}

Rewriting these derivatives in terms of quantities $\Axi^{12}$ introduced in Section \ref{sect_solution}, we prove \eqref{umder_A} and \eqref{dy0} from Theorem \ref{thm_main_g}.
Quite often it will be convenient to use \eqref{tempdiff_g} for the last derivative to separate it in two parts, making explicit the appearance of derivatives of $u_\alpha$, namely
\begin{equation}
\label{dy0-2}
\frac{\partial y_0}{\partial x_i}=-\frac {A_{x_i}^{12}}{t \varphi (Q_0)(x_i-y_0)^g}-\sum_{\alpha = 1}^{g-1} \frac {A_{u_{\alpha}}^{12}}{t \varphi (Q_0)(u_{\alpha}-y_0)^g} \frac {\partial u_{\alpha}}{\partial x_i}\,.
\end{equation}

\begin{remark}
\label{rmk_dynamics}
{\rm
Let us invert the M\"obius transformation \eqref{Mobius} to find $\hat x_1, \dots, \hat x_g, \hat u_1, \dots, \hat u_g$ with $g=d-1\,:$
\begin{equation*}
\hat u_g=1-y_0, \qquad \hat u_j=\frac{u_j(1-y_0)}{u_j-y_0}, \qquad
\hat x_k=\frac{x_k(1-y_0)}{x_k-y_0}, \quad j=1, \dots, g-1, \quad k=1, \dots, g\,.
\end{equation*}
The question asked at the beginning of Section \ref{sect_Chebyshev} concerned the dependence of $\hat u_j$ on $\hat x_k$. Given that equations \eqref{umder_g}  and \eqref{y0der_g} together with  the above system provide us with expressions for $\frac{\partial \hat u_j}{\partial  x_i}$ and $\frac{\partial \hat x_k}{\partial  x_i}$,
the derivatives $\frac{\partial \hat u_j}{\partial \hat x_k}$, $k=1, \dots, g,$ for each $j=1, \dots, g$ can be found as solutions of the following two linear systems
\begin{eqnarray*}
&&\frac{\partial \hat u_j}{\partial  x_i} = \sum_{k=1}^g \frac{\partial \hat u_j}{\partial \hat x_k} \frac{\partial \hat x_k}{\partial x_i} , \qquad i=1, \dots, g; %
\\
&&\frac{\partial \hat u_g}{\partial  x_i} = - \sum_{k=1}^g \frac{\partial y_0}{\partial \hat x_k} \frac{\partial \hat x_k}{\partial x_i} , \qquad i=1, \dots, g.
\end{eqnarray*}
}\end{remark}

Assume now that  $\c_1$ from Definition \ref{def:Tfamily} is a rational vector and $\c_2=0$.
In this case,  $\hat x_1, \dots, \hat x_g, \hat u_1, \dots, \hat u_g$ are endpoints of the support of a Chebyshev polynomial, and we obtain an answer to the question asked at the beginning of Section \ref{sect_Chebyshev} about Chebyshev dynamics. For the general value of the parameter $\c_1 \in \mathbb R^{d-1}$ with $\c_2=0$ it   solves the isoharmonic deformations and thus provides the resolution of a more general problem of constrained variation of Jacobi inversion for hyperelliptic curves of any genus.\\

\begin{theorem}\label{th:jacobiconst}
Equations \eqref{umder_g}  and \eqref{y0der_g}  describe a solution to the constrained variation of Jacobi inversion for hyperelliptic curves.
\end{theorem}

\subsubsection{Variation of further quantities on the ${\mathcal T}$-families of curves}
\label{sect_furtherTvar}

In order to prove Theorem \ref{thm_main_g} we need to find derivatives of $\Omega(\Paj)$ with respect to the independently varying branch points $x_i$, $i=1, \dots, g-1,$ and different from $a_j\,.$ Obtaining these derivatives for $a_j\neq u_\alpha$ and for $\Omega(\Pinfty)$ is a straightforward application of the Rauch variational formulas from Section \ref{sect_Rauch}. However, differentiating $\Omega(\Pualpha)$ with respect to $x_i$ is essentially different due to the presence of a dependent variable in the argument of $\Omega\,.$ As can easily be seen, Rauch formulas are not adapted for such a case. In this subsection, we obtain these derivatives, starting with the simpler ones.
\begin{proposition}
\label{proposition_Omega_infty}
Let $x_i$ be an independent branch point of the curve $\mathcal H_\x$ \eqref{TMobius}, $i=1, \dots, g,$ and $\Omega$ be the differential of the third kind given by \eqref{Omega_W}. The evaluation of $\Omega$ at $P_\infty$ is understood as defined by \eqref{evaluation}. Let $\Axi^{12}$ be the function given by \eqref{A12_aj_g} with $a_j=x_i$.
The following variational formula holds:
\begin{equation*}
\frac{\partial \Omega(P_\infty)}{\partial x_i } = -\frac{\Axi^{12}}{t}\frac{\prod_{\alpha=1}^{g-1}(x_i-u_\alpha)}{(x_i-y_0)^{g-1}}\,.
\end{equation*}
\end{proposition}
{\it Proof.}
Using the Rauch formulas and the fact that the branch points $u_1, \dots, u_{g-1}$ as well as $y_0$ are functions of $x_1, \dots, x_g$, from definition \eqref{Omega_W} of $\Omega$, we have
\begin{equation*}
\frac{\partial \Omega(P_\infty)}{\partial x_i }= \frac{1}{2} W(\Pxi, P_\infty) \Omega(\Pxi) + \frac{1}{2} \sum_{\alpha=1}^{g-1} W(\Pualpha, P_\infty) \Omega(\Pualpha) \frac{\partial u_\alpha}{\partial x_i} + 2 W(Q_0, P_\infty) \frac{\partial y_0}{\partial x_i}\,.
\end{equation*}
After plugging in the expression \eqref{Wxinfty} for $W(P, \Pinfty)$ evaluated for $P=\Pxi, \,\Pualpha, \,Q_0$ we see that the terms containing the normalization constants $I_k$ cancel out due to relation \eqref{tempdiff_g}. Thus, using expressions \eqref{umder_A} and \eqref{dy0} for the derivatives of $u_\alpha$ and $y_0$,  in terms of $\Axi^{12}$ we get
\begin{equation*}
\frac{\partial \Omega(P_\infty)}{\partial x_i }= -\frac{\Axi^{12}}{t(x_i-y_0)^{g-1}}
\left(
-\frac{x_i^g}{x_i-y_0} +
\sum_{\alpha=1}^{g-1}\frac{u_\alpha^g \prod_{\beta\neq\alpha}(x_i-u_\beta)}{(u_\alpha-y_0) \prod_{\beta\neq \alpha}(u_\alpha-u_\beta)}
+ \frac{y_0^g \prod_{\alpha=1}^{g-1}(x_i-u_\alpha)}{(x_i-y_0) \prod_{\alpha=1}^{g-1}(y_0-u_\alpha)}
\right).
\end{equation*}
Examining the rational function in the parenthesis as a function of $x_i$, we find that it has no poles in the $x_i$-sphere except at $x_i=\infty$ and is thus a polynomial of degree $g-1.$  Since it vanishes at $x_i=u_\alpha$, we conclude that this function is $\prod_{\alpha}(x_i-u_\alpha)\,.$
$\Box$

\begin{proposition}
\label{proposition_Omega_aj}
Let $x_i$ be an independent branch point of the curve $\mathcal H_\x$ \eqref{TMobius}, $i=1, \dots, g\,,$ and let $a_j$ be another branch point different from the dependent ones and from $x_i$, that is $a_j\in B\setminus\{u_1, \dots, u_{g-1}, x_i\}.$ Let $\Axi^{12}$ and $\Aaj^{12}$ be  given by \eqref{A12_aj_g}.
The following variational formula holds for the differential $\Omega$ defined by \eqref{Omega} and evaluated at the ramification point $\Paj$ according to \eqref{evaluation}:
\begin{equation*}
\frac{\partial \Omega(\Paj)}{\partial x_i } = \frac{\Omega(\Paj)}{2(x_i-a_j)}\frac{\Axi^{12}}{\Aaj^{12}} \frac{(a_j-y_0)^{g-1}}{(x_i-y_0)^{g-1}}\frac{\prod_{\alpha=1}^{g-1}(x_i-u_\alpha)}{\prod_{\alpha=1}^{g-1}(a_j-u_\alpha)}\,.
\end{equation*}
\end{proposition}
{\it Proof.}
Similarly to the proof of Proposition \ref{proposition_Omega_infty}, we have
\begin{equation*}
\frac{\partial \Omega(\Paj)}{\partial x_i }= \frac{1}{2} W(\Pxi, \Paj) \Omega(\Pxi) + \frac{1}{2} \sum_{\alpha=1}^{g-1} W(\Pualpha, \Paj) \Omega(\Pualpha) \frac{\partial u_\alpha}{\partial x_i} + 2 W(Q_0, \Paj) \frac{\partial y_0}{\partial x_i}.
\end{equation*}
Plugging in \eqref{Waj} for $W(P, \Paj)$ with $P=\Pxi, \Pualpha$ and $Q_0,$
we see that the terms containing the normalization constants $\beta_k^{(j)}$  cancel out again due to  \eqref{tempdiff_g}.
After plugging \eqref{umder_A} and \eqref{dy0} for the derivatives of $u_\alpha$ and $y_0$  in terms of $\Axi^{12}$, and rewriting the whole expression in terms of quantities $A^{12}$, we get
\begin{multline*}
\frac{\partial \Omega(\Paj)}{\partial x_i }= \frac{2\Axi^{12} }{t\varphi(P_{a_j})(x_i-y_0)^g} \left(\frac{1}{x_i-a_j}
\right.
\\
\left.
-(x_i-y_0)\prod_{\beta=1}^{g-1}(x_i-u_\beta)\sum_{\alpha=1}^{g-1} \frac{1}{(u_\alpha-a_j)(u_\alpha-y_0)(x_i-u_\alpha)\prod_{\beta\neq \alpha}(u_\alpha-u_\beta)}
+
   \frac {\prod_{\alpha=1}^{g-1}(x_i-u_{\alpha})}{(a_j-y_0) \prod _{\alpha=1}^{g-1}(y_0-u_{\alpha}) }
 \right)\,.
\end{multline*}
It remains to use rational identity \eqref{rat1} from Lemma \ref{lemma_rational1} for the sum over $\alpha$ to finish the proof.
$\Box$

Now we turn to differentiating $\Omega(\Pum)$ with respect to an independent branch point $x_i$. Applying the chain rule, we will need to differentiate  $\Omega(\Pum)$ with respect to all dependent branch points and, in particular, with respect to $u_m$. The Rauch formulas \eqref{RauchOmegaaj} from Section \ref{sect_Rauch} do not allow for such differentiation, we need to use Lemma \ref{lemma_epsilon_g} instead.

\begin{proposition} Let $m\in\{1, \dots, g-1\}$ be fixed, $u_m$ be one of the dependent branch points of the curve $\mathcal H_\x$ \eqref{TMobius}, and $x_i$ be an independent branch point of the same curve. Let $\Aaj^{12}$ and $\Aum^{12}$ be defined by \eqref{A12_aj_g}.  The following variational formula holds for the differential $\Omega$ defined by \eqref{Omega} and evaluated at the ramification point $\Pum$ according to \eqref{evaluation}:
\label{prop_Omega_um_diff_g}
\begin{equation*}
\frac{\partial \Omega(P_{u_m})}{\partial x_i} =-\frac{\Omega(P_{u_m})}{2}
 \left( \frac{1}{x_i-u_m} - \sum_{\substack{ \alpha=1\\ \alpha\neq m}}^{g-1}\frac{1}{u_m-u_\alpha}  +\frac{1}{A_{u_m}^{12}}\sum_{\substack{a_j\in B\\a_j\neq u_m}} \frac{A_{a_j}^{12}}{a_j-u_m} + \frac{g-1}{u_m-y_0}      \right)\frac{\partial u_m}{\partial x_i}\,.
\end{equation*}
\end{proposition}

{\it Proof.}
Differentiating $\Omega(P_{u_m})$ with respect to $x_i$ with the help of Rauch variational formulas \eqref{RauchOmegaaj} and Lemma \ref{lemma_epsilon_g}, we get
\begin{multline*}
\frac{\partial \Omega(P_{u_m})}{\partial x_i} = \frac{\partial^{{\rm Rauch}} \Omega(P_{u_m})}{\partial x_i} +\frac{\partial \Omega(P_{u_m})}{\partial y_0}\frac{\partial y_0}{\partial x_i}
 + \sum_{\substack{k=1 \\ k\neq m}}^{g-1}\frac{\partial^{{\rm Rauch}} \Omega(P_{u_m})}{\partial u_k} \frac{\partial u_k}{\partial x_i} + \frac{\partial \Omega(P_{u_m})}{\partial u_m} \frac{\partial u_m}{\partial x_i}
 \\
 =\frac{1}{2}\Omega(\Pxi)W(\Pxi, \Pum) + 2W(Q_0, \Pum) \frac{\partial y_0}{ \partial x_i} + \frac{1}{2}\sum_{\substack{k=1 \\ k\neq m}}^{g-1} \Omega(\Puk) W(\Puk, \Pum) \frac{\partial u_k}{\partial x_i}
 \\
 -\left(\frac{1}{2}\sum_{\substack{a_j\in B \\ a_j\neq u_m}} \Omega(\Paj) W(\Paj, \Pum)  + 2W(Q_0, \Pum) \right)\frac{\partial u_m}{\partial x_i}\,.
\end{multline*}
Using now \eqref{res4} in the last line and \eqref{Wualpha} together with \eqref{vcond} to write $W(\cdot, \Pum)$, we have
\begin{multline*}
\frac{\partial \Omega(P_{u_m})}{\partial x_i}
 =\frac{1}{2}\Omega(\Pxi)\left( \frac{1}{x_i-u_m} \frac{\varphi(\Pxi)}{\varphi(\Pum)} - \sum_{\alpha=1}^g\gamma_\alpha^{(m)} v_\alpha(\Pxi)\right)
  + 2\frac{\partial y_0}{ \partial x_i}\left( \frac{1}{y_0-u_m} \frac{\varphi(Q_0)}{\varphi(\Pum)} - \gamma_g^{(m)} \right)
  \\
  + \frac{1}{2}\sum_{\substack{k=1\\k\neq m}}^{g-1} \Omega(\Puk) \left( \frac{1}{u_k-u_m} \frac{\varphi(\Puk)}{\varphi(\Pum)} - \gamma_k^{(m)} \right) \frac{\partial u_k}{\partial x_i}
 +\underset{P=\Pum}{\rm res} \frac{\Omega(P) W(P, \Pum)}{du} \frac{\partial u_m}{\partial x_i}\,.
\end{multline*}
Plugging in \eqref{umder_g} and \eqref{y0der_g} for derivatives of $u_m$ and $y_0$ in terms of differentials $\Omega$ and $v$, we see that all terms containing normalization constants $\gamma_k^{(m)}$ disappear, except the term with $\gamma_m^{(m)}$. This remaining term cancels the corresponding term coming out of the residue after we express this residue as in \eqref{res-diff_g}:
\begin{multline}
\label{derOmegaum_temp}
\frac{\partial \Omega(P_{u_m})}{\partial x_i}
 =\frac{1}{2}\frac{\Omega(\Pxi)}{\varphi(\Pum)}\left( \frac{\varphi(\Pxi)}{x_i-u_m}
  - \frac{v_g(\Pxi)\varphi(Q_0)}{y_0-u_m}   -\sum_{\substack{k=1\\k\neq m}}^{g-1}   \frac{v_k(\Pxi)\varphi(\Puk)}{u_k-u_m}  \right)
 \\
 +\left(\frac{1}{2}\sum_{\substack{a_j\in B \\ a_j\neq u_m}}  \frac{\Omega(\Paj)\varphi(\Paj)}{a_j-u_m} + 2\frac{\varphi(Q_0)}{y_0-u_m}\right) \frac{ \Omega(\Pxi)}{\Omega(\Pum)} \frac{v_m(\Pxi)}{\varphi(P_{u_m})}\,.
\end{multline}
Let us now plug in explicit expressions \eqref{v} and \eqref{vg} for the differentials $v_k$ and us the sum of residues \eqref{res1} for $\varphi(q_0)$ in the last term.  We then evaluate the sum over $k$ using the rational identity \eqref{rat-alpha2} from Section \ref{sect_rational}:
\begin{multline*}
\sum_{\substack{k=1\\k\neq m}}^{g-1}   \frac{v_k(\Pxi)\varphi(\Puk)}{u_k-u_m}
=
\varphi(\Pxi)\prod_{\alpha}(x_i-u_\alpha) \sum_{\substack{k=1\\k\neq m}}^{g-1} \left( \frac{1}{x_i-u_k}+\frac{1}{u_k-y_0}\right)\frac{1}{(u_k-u_m)\prod_{\alpha\neq k}(u_k-u_\alpha)}
\\
=\frac{\varphi(\Pxi)\prod_{\alpha\neq m}(x_i-u_\alpha)}{\prod_{\alpha\neq m}(u_m-u_\alpha)} \frac{x_i-y_0}{u_m-y_0}\left( \sum_{\substack{\alpha=1\\ \alpha\neq m}}^{g-1} \frac{1}{u_m-u_\alpha} - \frac{1}{x_i-u_m} +\frac{1}{u_m-y_0} \right)
\\
+\frac{\varphi(\Pxi)}{x_i-u_m} - \frac{\varphi(\Pxi)\prod_{\alpha=1}^{g-1}(x_i-u_\alpha)}{(y_0-u_m)\prod_{\alpha=1}^{g-1}(y_0-u_\alpha)}\,.
\end{multline*}
Note that the terms in the last line cancel the first two terms in our expression for the derivative of $\Omega(\Pum).$ For the last line in \eqref{derOmegaum_temp}, due to \eqref{res1} we have
\begin{multline*}
\frac{1}{2}\sum_{\substack{a_j\in B\\a_j\neq u_m}} \frac{\Omega(\Paj)\varphi(\Paj)}{a_j-u_m} + 2\frac{\varphi(Q_0)}{y_0-u_m}
=
\frac{2}{t} \sum_{\substack{a_j\in B\\a_j\neq u_m}} \frac{\Aaj^{12}}{(a_j-y_0)^g(a_j-u_m)}
-\frac{2}{t} \sum_{a_j\in B} \frac{\Aaj^{12}}{(a_j-y_0)^g(y_0-u_m)}
\\
=-\frac{2}{t(y_0-u_m)} \sum_{\substack{a_j\in B\\a_j\neq u_m}} \frac{\Aaj^{12}}{(a_j-y_0)^{g-1}(a_j-u_m)}
+\frac{2}{t}  \frac{\Aum^{12}}{(u_m-y_0)^{g+1}}\,.
\end{multline*}
Now we can simplify the sum over the branch points using Lemma \ref{lemma_simplify} with $k=1$. Combining these results, re-expressing everything in terms of $A_{a_j}^{12}$ and using \eqref{umder_A} to write the overall factor as derivative of $u_m$, we prove the proposition.
$\Box$

Assuming again that $u_\alpha$ stands for a dependent branch poing and $a_j$ denotes any of the remaining branch points, we have in a straightforward way for the holomorphic non-normalized differential $\varphi$ defined by \eqref{phi}, \eqref{phi_Q0}, \eqref{phi_Pk}:
\begin{equation}
\label{der_phi_aj_g}
\frac{\partial \varphi(P_{a_j})}{\partial x_i} = \frac{\varphi(P_{a_j})}{2}\left( \frac{1}{a_j-x_i} + \sum_{\alpha=1}^{g-1}\frac{1}{a_j-u_\alpha} \frac{\partial u_\alpha}{\partial x_i} \right)\qquad \text{if}\quad a_j\neq u_k\,.
\end{equation}
\begin{equation}
\label{der_phi_um}
\frac{\partial \varphi(P_{u_m})}{\partial x_i} = \frac{\varphi(P_{u_m})}{2}\left( \frac{1}{u_m-x_i} -\frac{\partial u_m}{\partial x_i}\sum_{\substack{a_j\in B\\a_j\neq u_m}} \frac{1}{u_m-a_j}  +\sum_{\substack{\alpha=1\\ \alpha \neq m}}^{g-1} \frac{1}{u_m-u_\alpha} \frac{\partial u_\alpha}{\partial x_i}  \right)
\end{equation}
\begin{equation}
\label{phi_q0_g}
\frac{\partial \varphi(Q_0)}{\partial x_i} = \frac{\varphi(Q_0)}{2}\left(\sum_{a_j\in B} \frac{1}{a_j-y_0} \frac{\partial y_0}{\partial x_i} -\frac{1}{x_i-y_0}
- \sum_{\alpha=1}^{g-1}\frac{1}{u_\alpha-y_0} \frac{\partial u_\alpha}{\partial x_i}
  \right)
\,.
\end{equation}

\section{Proof of Theorem \ref{thm_main_g}}
\label{sect_proof}

It will be convenient for what follows to write the constrained Schlesinger system \eqref{constrained} in matrix components:

\begin{eqnarray}
\label{Aaj12}
&&\frac{\partial A_{a_j}^{12}}{\partial x_i} = 2\frac{A^{11}_{x_i} A_{a_j}^{12} - A^{12}_{x_i} A_{a_j}^{11} }{x_i-a_j} + 2 \sum_{k=1}^{g-1} \frac{A^{11}_{u_k} A_{a_j}^{12} - A^{12}_{u_k} A_{a_j}^{11} }{u_k-a_j} \frac{\partial u_k}{\partial x_i}
\quad\mbox{for} \quad a_j\notin\{u_1, \dots, u_{g-1}\};\qquad\qquad\qquad
\\
&&\frac{\partial A_{u_m}^{12}}{\partial x_i} = 2\frac{A^{11}_{x_i} A_{u_m}^{12} - A^{12}_{x_i} A_{u_m}^{11} }{x_i-u_m}
+ 2 \sum_{\substack{k=1\\k\neq m}}^{g-1} \frac{A^{11}_{u_k} A_{u_m}^{12} - A^{12}_{u_k} A_{u_m}^{11} }{u_k-u_m} \frac{\partial u_k}{\partial x_i}
-2\frac{\partial u_m}{\partial x_i} \!\!\sum_{\substack{a_j\in B \\a_j\neq u_m}}\frac{  \Aaj^{11}\Aum^{12}-\Aaj^{12}\Aum^{11}}{a_j-u_m} ; \qquad\quad
\label{Aum12}
\end{eqnarray}
\begin{eqnarray}
\label{Aaj11}
&&\frac{\partial \Aaj^{11}}{\partial x_i} = \frac{\Axi^{12}\Aaj^{21} - \Axi^{21}\Aaj^{12}}{x_i-a_j} + \sum_{k=1}^{g-1} \frac{\Auk^{12}\Aaj^{21} - \Auk^{21}\Aaj^{12}}{u_k-a_j} \frac{\partial u_k}{\partial x_i}
\quad\mbox{for} \quad a_j\notin\{u_1, \dots, u_{g-1}\}; \qquad\qquad\qquad
\\
\label{Aum11}
&&\frac{\partial A_{u_m}^{11}}{\partial x_i} = \frac{\Axi^{12} \Aum^{21} - \Axi^{21}\Aum^{12}}{x_i-u_m}
+ \sum_{\substack{k=1\\k\neq m}}^{g-1} \frac{  \Auk^{12}\Aum^{21}-\Auk^{21}\Aum^{12}}{u_k-u_m} \frac{\partial u_k}{\partial x_i}
 - \frac{\partial u_m}{\partial x_i}\!\!\sum_{\substack{a_j\in B \\a_j\neq u_m}} \frac{  \Aaj^{12}\Aum^{21}-\Aaj^{21}\Aum^{12}}{a_j-u_m};
 \qquad\qquad\quad
\end{eqnarray}

\begin{eqnarray}
\label{Aaj21}
&&\frac{\partial \Aaj^{21}}{\partial x_i} = 2 \frac{\Axi^{21} \Aaj^{11} - \Axi^{11}\Aaj^{21}}{x_i-a_j} + 2\sum_{k=1}^{g-1}
\frac{\Auk^{21}\Aaj^{11} - \Auk^{11}\Aaj^{21}}{u_k-a_j} \frac{\partial u_k}{\partial x_i}
\quad\mbox{for} \quad a_j\notin\{u_1, \dots, u_{g-1}\};\qquad\qquad\qquad
\\
\label{Aum21}
&&\frac{\partial A_{u_m}^{21}}{\partial x_i} = 2 \frac{\Axi^{21} \Aum^{11} - \Axi^{11}\Aum^{21}}{x_i-u_m} + 2\sum_{\substack{k=1\\k\neq m}}^{g-1}
\frac{\Auk^{21}\Aum^{11} - \Auk^{11}\Aum^{21}}{u_k-u_m} \frac{\partial u_k}{\partial x_i}
-2\frac{\partial u_m}{\partial x_i}\!\!\!\sum_{\substack{a_j\in B \\a_j\neq u_m}}
\frac{  \Aaj^{21}\Aum^{11}-\Aaj^{11}\Aum^{21}}{a_j-u_m} \,. \qquad\quad
\end{eqnarray}

\subsection{Proof for (12)-components}

\subsubsection{ Proof for $\Aaj^{12}$ with $a_j\notin \{u_1, \dots, u_{g-1}\}$}
Let us prove that the functions defined in Theorem \ref{thm_main_g} satisfy equation \eqref{Aaj12} of the constrained Schlesinger system, which in terms of  $\betaaj$ it takes the following form:
\begin{equation}
\label{Aaj_system}
\frac{\partial \Aaj^{12}}{\partial x_i} = \frac{\Axi^{12}-\Aaj^{12}}{2(x_i-a_j)} - \frac{g}{t}\frac{\Aaj^{12}\Axi^{12}}{x_i-a_j} \left( \frac{\betaxi}{\Axi^{12}} - \frac{\betaaj}{\Aaj^{12}} \right) + \sum_{\alpha=1}^{g-1} \left[
\frac{\Aualpha^{12}-\Aaj^{12}}{2(u_\alpha-a_j)} - \frac{g}{t}\frac{\Aaj^{12}\Aualpha^{12}}{u_\alpha-a_j} \left( \frac{\betaualpha}{\Aualpha^{12}} - \frac{\betaaj}{\Aaj^{12}} \right)
\right]\frac{\partial u_\alpha}{\partial x_i}\,.
\end{equation}
We want to prove that the right hand side coincides with the derivative of the function $\Aaj^{12}$ from  Theorem \ref{thm_main_g}. We can compute this derivative  using Proposition \ref{proposition_Omega_aj} for $\partial_{x_i}\Omega(\Paj)$ and \eqref{der_phi_aj_g} for $\partial_{x_i}\varphi(\Paj)$ as well as \eqref{dy0} for the derivative of $y_0:$

\begin{multline}
\label{Aaj_derivative}
\frac{\partial \Aaj^{12}}{\partial x_i} = \frac{\partial }{\partial x_i}\left\{ \frac{t}{4}\Omega(\Paj) \varphi(\Paj)(a_j-y_0)^g\right\}
=
\frac{\Aaj^{12}}{2}\left\{  \frac{\partial_{x_i} \Omega(\Paj)}{\Omega(\Paj)} + \frac{\partial_{x_i} \varphi(\Paj)}{\varphi(\Paj)} -\frac{2g}{(a_j-y_0)}\frac{\partial y_0}{\partial x_i}
\right\}
\\
=
\frac{\Aaj^{12}}{2}\left\{  \frac{1}{x_i-a_j}\frac{\Axi^{12}}{\Aaj^{12}} \frac{(a_j-y_0)^{g-1}}{(x_i-y_0)^{g-1}}\frac{\prod_{\alpha=1}^{g-1}(x_i-u_\alpha)}{\prod_{\alpha=1}^{g-1}(a_j-u_\alpha)}
 + \left( \frac{1}{a_j-x_i} + \sum_{\alpha=1}^{g-1}\frac{1}{a_j-u_\alpha} \frac{\partial u_\alpha}{\partial x_i} \right) \right.
 \\
 \left.+\frac{2gA_{x_i}^{12}\prod_{\alpha=1}^{g-1}(x_i-u_{\alpha})}{t(a_j-y_0)\varphi (Q_0)(x_i-y_0)^g \prod _{\alpha=1}^{g-1}(y_0-u_{\alpha})}
\right\}.
\end{multline}
Thus we need to prove that the difference of right hand sides of \eqref{Aaj_system} and \eqref{Aaj_derivative} vanishes. As is easy to see the terms in \eqref{Aaj_derivative} coming out of the derivative $\partial_{x_i}\varphi(\Paj)$ cancel the corresponding terms in \eqref{Aaj_system}.

Let us work on the first term in the sum over $\alpha$ from \eqref{Aaj_system}. Plugging in \eqref{umder_A} for the derivative $\partial_{x_i} u_\alpha$, we obtain an expression which we can reduce using the rational identity \eqref{rat3} as follows:
\begin{multline*}
\sum_{\alpha=1}^{g-1}
\frac{\Aualpha^{12}}{2(u_\alpha-a_j)} \frac{\partial u_\alpha}{\partial x_i}
=-\frac{\Axi^{12}\prod_{\beta}(x_i-u_\beta)}{2(x_i-y_0)^{g-1}}\sum_{\alpha=1}^{g-1}\frac{(u_\alpha-y_0)^{g-1}}{(u_\alpha-a_j)(x_i-u_\alpha)\prod_{\beta\neq \alpha}(u_\alpha-u_\beta)}
\\
=-\frac{\Axi^{12}\prod_{\beta}(x_i-u_\beta)}{2(x_i-y_0)^{g-1}(x_i-a_j)}\sum_{\alpha=1}^{g-1}\left( \frac{1}{x_i-u_\alpha} - \frac{1}{a_j-u_\alpha}\right)\frac{(u_\alpha-y_0)^{g-1}}{\prod_{\beta\neq \alpha}(u_\alpha-u_\beta)}
\\
=-\frac{\Axi^{12}}{2(x_i-a_j)}
+ \frac{\Axi^{12}}{2(x_i-a_j)}\frac{(a_j-y_0)^{g-1}\prod_\alpha(x_i-u_\alpha)}{(x_i-y_0)^{g-1}\prod_{\alpha}(a_j-u_\alpha)}.
\end{multline*}
Note that the first term of this expression will cancel the first term of \eqref{Aaj_system} and the second one will cancel the first term of the right hand side of \eqref{Aaj_derivative}. Thus for the difference between \eqref{Aaj_derivative} and \eqref{Aaj_system}, which we denote $X_{a_j}$, we have
\begin{multline}
\label{system-derivative}
\frac{X_{a_j}}{\Aaj^{12}}:=  \frac{g}{t}\frac{\Axi^{12}}{x_i-a_j} \left( \frac{\betaxi}{\Axi^{12}} - \frac{\betaaj}{\Aaj^{12}} \right)
+ \sum_{\alpha=1}^{g-1}
  \frac{g}{t}\frac{\Aualpha^{12}}{u_\alpha-a_j} \left( \frac{\betaualpha}{\Aualpha^{12}} - \frac{\betaaj}{\Aaj^{12}} \right)
\frac{\partial u_\alpha}{\partial x_i}
\\
 + \frac{gA_{x_i}^{12}\prod_{\alpha=1}^{g-1}(x_i-u_{\alpha})}{t(a_j-y_0)\varphi (Q_0)(x_i-y_0)^g \prod _{\alpha=1}^{g-1}(y_0-u_{\alpha})} \,.
\end{multline}
Plugging in \eqref{difference_betas_2} for differences $ \frac{\betaai}{\Aai^{12}} - \frac{\betaaj}{\Aaj^{12}}$ and \eqref{umder_A} for the derivative of $u_\alpha$, we have
\begin{multline*}
\frac{tX_{a_j}\varphi(Q_0)}{g\Aaj^{12}\Axi^{12}}=  - \sum_{l=0}^{g-1}  \frac{L_l}{l!}\sum_{k=0}^{g-l-1}\frac{1}{(x_i-y_0)^{k+1}(a_j-y_0)^{g-l-k}}
+\frac{\prod_{\alpha=1}^{g-1}(x_i-u_{\alpha})}{(a_j-y_0)(x_i-y_0)^g \prod _{\alpha=1}^{g-1}(y_0-u_{\alpha})}
\\
+\frac{\prod_{\alpha=1}^{g-1}(x_i-u_\alpha)}{(x_i-y_0)^{g-1}}
 \sum_{l=0}^{g-1}  \frac{L_l}{l!}\sum_{k=0}^{g-l-1}\frac{1}{(a_j-y_0)^{g-l-k}}
\sum_{\alpha=1}^{g-1} \frac{(u_\alpha-y_0)^{g-k-2}}{(x_i-u_\alpha)\prod_{\beta\neq \alpha}(u_\alpha-u_\beta)} \,.
\end{multline*}
Now it remains to use the rational identity \eqref{rat2} in the last line in all sums over $\alpha$ noting that the sum with $l=0$ and $k=g-1$ needs to be singled out and split into two sums of partial fractions. This yields cancellation of all the terms and shows that the difference $X_{a_j}$ between \eqref{Aaj_derivative} and \eqref{Aaj_system} is zero.
$\Box$

\subsubsection{  Proof for $\Aaj^{12}$ with $a_j=u_m$}
Here we want to prove that the functions defined in Theorem \ref{thm_main_g} satisfy equation \eqref{Aum12} of the constrained Schlesinger system. Let us first rewrite \eqref{Aum12}  in terms of  $\betaaj$; it takes the form:
\begin{multline}
\label{Aum_system}
\hspace{-0.3cm}
\frac{\partial \Aum^{12}}{\partial x_i} = \frac{\Axi^{12}-\Aum^{12}}{2(x_i-u_m)} - \frac{g}{t}\frac{\Aum^{12}\Axi^{12}}{x_i-u_m} \left( \frac{\betaxi}{\Axi^{12}} - \frac{\betaum}{\Aum^{12}} \right)
+ \sum_{\substack{\alpha=1\\\alpha\neq m}}^{g-1} \left[
\frac{\Aualpha^{12}-\Aum^{12}}{2(u_\alpha-u_m)} - \frac{g}{t}\frac{\Aum^{12}\Aualpha^{12}}{u_\alpha-u_m} \left( \frac{\betaualpha}{\Aualpha^{12}} - \frac{\betaum}{\Aum^{12}} \right)
\right]\frac{\partial u_\alpha}{\partial x_i}
\\
-\frac{\partial u_m}{\partial x_i} \sum_{\substack{a_j\in B\\a_j\neq u_m}} \left[
\frac{\Aaj^{12}-\Aum^{12}}{2(a_j-u_m)} - \frac{g}{t}\frac{\Aum^{12}\Aaj^{12}}{a_j-u_m} \left( \frac{\betaaj}{\Aaj^{12}} - \frac{\betaum}{\Aum^{12}} \right)\right]\,.
\end{multline}
Let us now compute the derivative of $\Aum^{12}$ defined by \eqref{A12_aj_g} with $a_j=u_m$ with respect to an independently varying branch point $x_i$. Putting the derivatives from Proposition \ref{prop_Omega_um_diff_g} and  \eqref{der_phi_um} together, we have
\begin{multline}
\label{Aum_derivative}
\frac{\partial \Aum^{12}}{\partial x_i} = \frac{\partial }{\partial x_i}\left\{ \frac{t}{4}\Omega(\Pum) \varphi(\Pum)(u_m-y_0)^g\right\}
=
\frac{\Aum^{12}}{2}\left\{  \frac{2g}{u_m-y_0}\left( \frac{\partial u_m}{\partial x_i}-\frac{\partial y_0}{\partial x_i}\right)
\right.
\\
-\frac{\partial u_m}{\partial x_i}
 \left( \frac{1}{x_i-u_m} - \sum_{\substack{ \alpha=1\\ \alpha\neq m}}^{g-1}\frac{1}{u_m-u_\alpha}  +\frac{1}{A_{u_m}^{12}}\sum_{\substack{a_j\in B\\a_j\neq u_m}} \frac{A_{a_j}^{12}}{a_j-u_m} + \frac{g-1}{u_m-y_0}    \right)
\\
\left.
+\frac{1}{u_m-x_i} -\frac{\partial u_m}{\partial x_i}\sum_{\substack{a_j\in B\\a_j\neq u_m}} \frac{1}{u_m-a_j}  +\sum_{\substack{\alpha=1\\ \alpha \neq m}}^{g-1} \frac{1}{u_m-u_\alpha} \frac{\partial u_\alpha}{\partial x_i}
\right\}
\end{multline}
We want to prove that the difference between this derivative and the right hand side of \eqref{Aum_system} vanishes. Note first that, in this difference, the terms in the last line of \eqref{Aum_derivative} as well as the sum over the branch points in the second line cancel the corresponding terms in \eqref{Aum_system}. Let us now compute, using \eqref{difference_betas_1} for the difference of ${\betaaj}/{\Aaj^{12}}$, the following sum over $a_j$ in \eqref{Aum_system}:
\begin{multline*}
\sum_{\substack{a_j\in B\\a_j\neq u_m}}\frac{\Aaj^{12}}{a_j-u_m} \left( \frac{\betaaj}{\Aaj^{12}} - \frac{\betaum}{\Aum^{12}} \right)
=
\frac{1}{\varphi(Q_0)}\sum_{l=0}^{g-1}\frac{L_l}{l!}\left(\sum_{\substack{a_j\in B\\a_j\neq u_m}}\frac{\Aaj^{12}}{(a_j-u_m)(a_j-y_0)^{g-l}}
-
\frac{\sum_{a_j\neq u_m}\frac{\Aaj^{12}}{a_j-u_m}}{(u_m-y_0)^{g-l}}      \right)\,.
\end{multline*}
Due to Lemma \ref{lemma_simplify} this simplifies to
\begin{equation*}
\sum_{\substack{a_j\in B\\a_j\neq u_m}}\frac{\Aaj^{12}}{a_j-u_m} \left( \frac{\betaaj}{\Aaj^{12}} - \frac{\betaum}{\Aum^{12}} \right)
=
\frac{A_{u_m}^{12}}{\varphi(Q_0)}\sum_{l=0}^{g-1}\frac {(g-l)L_l}{l!(u_m-y_0)^{g-l+1}} + \frac{t }{u_m-y_0}
\,.
\end{equation*}
 Combining all these results together, we compute the difference between the right hand sides of \eqref{Aum_system} and \eqref{Aum_derivative}. Denoting the result by $X_{u_m}$, after some simplification, we have
\begin{multline}
\label{Xum}
X_{u_m} = \frac{\Axi^{12}}{2(x_i-u_m)}
+ \frac{g}{t}\frac{(A_{u_m}^{12})^2}{\varphi(Q_0)}\frac{\partial u_m}{\partial x_i}\sum_{l=0}^{g-1}\frac {(g-l)L_l}{l!(u_m-y_0)^{g-l+1}}
+ \frac{g\Aum^{12}}{u_m-y_0} \frac{\partial y_0}{\partial x_i}
\\
-\frac{g}{t}\frac{\Aum^{12}\Axi^{12}}{x_i-u_m} \left( \frac{\betaxi}{\Axi^{12}} - \frac{\betaum}{\Aum^{12}} \right)
+ \sum_{\substack{\alpha=1\\ \alpha\neq m}}^{g-1} \left[
\frac{\Aualpha^{12}}{2(u_\alpha-u_m)} - \frac{g}{t}\frac{\Aum^{12}\Aualpha^{12}}{u_\alpha-u_m} \left( \frac{\betaualpha}{\Aualpha^{12}} - \frac{\betaum}{\Aum^{12}} \right)
\right]\frac{\partial u_\alpha}{\partial x_i}
\\
+\frac{\Aum^{12}}{2}\frac{\partial u_m}{\partial x_i}
 \left( \frac{1}{x_i-u_m} - \sum_{\substack{\alpha=1\\ \alpha\neq m}}^{g-1}\frac{1}{u_m-u_\alpha}  + \frac{g-1}{u_m-y_0}    \right)\,.
\end{multline}
Now we plug in expressions from \eqref{umder_A} and \eqref{dy0} for derivatives $\partial_{x_i}u_\alpha$ and $\partial_{x_i} y_0$, respectively, and also use \eqref{difference_betas_1} for the differences $\frac{\betaaj}{\Aaj^{12}} - \frac{\betaum}{\Aum^{12}}$. This, after changing the order of summation and splitting some fractions into partial fractions, becomes
\begin{multline*}
\frac{X_{u_m}}{\Axi^{12}} = \frac{1}{2(x_i-u_m)}
- \frac{g}{t} \frac{A_{u_m}^{12}\prod_{\alpha\neq m}(x_i-u_\alpha)}{\varphi(Q_0)(x_i-y_0)^{g-1}\prod_{\alpha\neq m}(u_m-u_\alpha)}\sum_{l=0}^{g-1}\frac {(g-l)L_l}{l!(u_m-y_0)^{2-l}}
\\
- \frac{g}{t}\frac{\Aum^{12}\prod_{\alpha=1}^{g-1}(x_i-u_{\alpha})}{\varphi (Q_0)(u_m-y_0)(x_i-y_0)^g \prod _{\alpha=1}^{g-1}(y_0-u_{\alpha}) }
-\frac{g}{t}\frac{\Aum^{12}}{\varphi(Q_0)(x_i-u_m)}  \sum_{l=0}^{g-1} \left( \frac{L_l}{l!(x_i-y_0)^{g-l}} - \frac{L_l}{l!(u_m-y_0)^{g-l}} \right)
\\
- \frac{\prod_{\beta=1}^{g-1}(x_i-u_\beta)}{2(x_i-u_m)(x_i-y_0)^{g-1}}\sum_{\alpha\neq m} \left( \frac{1}{u_\alpha-u_m} + \frac{1}{x_i-u_\alpha}\right)
 \frac{(u_\alpha-y_0)^{g-1}}{\prod_{\beta\neq \alpha}(u_\alpha-u_\beta)}
\\
+ \frac{g}{t}\frac{\Aum^{12}\prod_{\beta=1}^{g-1}(x_i-u_\beta)}{\varphi(Q_0)(x_i-y_0)^{g-1}(x_i-u_m)}\sum_{l=0}^{g-1}\frac{L_l}{l!}\sum_{\alpha\neq m}  \left( \frac{1}{u_\alpha-u_m} + \frac{1}{x_i-u_\alpha}\right)
\frac{(u_\alpha-y_0)^{l-1}}{\prod_{\beta\neq \alpha}(u_\alpha-u_\beta)}
\\
- \frac{g}{t}\frac{\Aum^{12}\prod_{\beta=1}^{g-1}(x_i-u_\beta)}{\varphi(Q_0)(x_i-y_0)^{g-1}(x_i-u_m)}\sum_{l=0}^{g-1}  \frac{L_l}{l!(u_m-y_0)^{g-l}}\sum_{\alpha\neq m}  \left( \frac{1}{u_\alpha-u_m} + \frac{1}{x_i-u_\alpha}\right)
 \frac{(u_\alpha-y_0)^{g-1}}{\prod_{\beta\neq \alpha}(u_\alpha-u_\beta)}
\\
- \frac{(u_m-y_0)^{g-1}\prod_{\alpha\neq m}(x_i-u_\alpha)}{2(x_i-y_0)^{g-1}\prod_{\alpha\neq m}(u_m-u_\alpha)}
 \left( \frac{1}{x_i-u_m} - \sum_{\alpha\neq m}\frac{1}{u_m-u_\alpha}  + \frac{g-1}{u_m-y_0}    \right)\,.
\end{multline*}
It remains to compute the the sums over alpha using Lemmas \ref{lemma_rational1} and \ref{lemma_rational3}. This is a quite lengthy but absolutely straightforward calculation, which shows that $X_{u_m}=0.$ To simplify this calculation, one can notice that the terms containing the factor of $\Aum^{12}$ cancel each other.
$\Box$

\subsection{Proof for (11)-components}

\subsubsection{Technical lemmas}

We first prove the technical lemmas that we use in the proof of Theorem \ref{thm_main_g} for the (11)-components of the matrices. In all these lemmas,
$C_{p_1\dots p_l} $ are the coefficients \eqref{Coefficients} of the polynomial $L_l$ \eqref{L}; $a_j$ is an arbitrary element of the set of branch points $B:=\{0,1,x_1, \dots, x_{g}, u_1, \dots, u_{g-1}\}$ of our hyperelliptic curve; and $\betaaj$ are given by \eqref{betas}.

\begin{lemma} For any two distinct branch points $a_i$ and $a_j$, we have
\label{lemma_A}
\begin{multline*}
\sum_{l=0}^{g-1} \sum_{\substack{p_1+2p_2+\dots +lp_l=l \\ p_1,\dots, p_l\geq 0}}\frac{C_{p_1\dots p_l} \Sigma_1^{p_1}\dots\Sigma_l^{p_l}}{l!\varphi(Q_0)(a_j-y_0)^{g-l}}\left\{\frac{1}{2(a_i-y_0)}-\sum_{k=1}^l \frac{kp_k}{\Sigma_k}\frac{1}{(a_i-y_0)^{k+1}} \right\}
\\
=-\frac{1}{2(a_i-a_j)} \left( \frac{\betaai}{\Aai^{12}}-\frac{\betaaj}{\Aaj^{12}}\right)\,.
\end{multline*}
\end{lemma}

{\it Proof.} Using \eqref{difference_betas_1} in the right hand side and the identity
\begin{equation}
\label{eq:identityLemmaA1}
\frac{1}{(a_i-y_0)^{m}}-\frac{1}{(a_j-y_0)^{m}}=
(a_j-a_i)\sum_{k=0}^{m-1}\frac{1}{(a_j-y_0)^{m-k}(a_i-y_0)^{k+1}}\,,
\end{equation}
we rewrite this relation in the form
\begin{multline*}
-\sum_{l=0}^{g-1} \sum_{\substack{p_1+2p_2+\dots +lp_l=l \\ p_1,\dots, p_l\geq 0}}\frac{C_{p_1\dots p_l} \Sigma_1^{p_1}\dots\Sigma_l^{p_l}}{l!(a_j-y_0)^{g-l}}\sum_{k=1}^l \frac{kp_k}{\Sigma_k}\frac{1}{(a_i-y_0)^{k+1}}
= \frac{1}{2} \sum_{l=0}^{g-2}\frac{L_l}{l!}\sum_{k=1}^{g- l-1}\frac{1}{(a_j-y_0)^{g- l-k}(a_i-y_0)^{k+1}}\,;
\end{multline*}
note that one of the terms in the left hand side cancelled the term with $k=0$ on the right and the term with $l=g-1$ on the right is not included as it has no contribution. Now we work on the left hand side, writing explicitly the coefficients $C_{p_1\dots p_l}$ \eqref{Coefficients} of the polynomials $L_l$. Since each term of the polynomial $L_l$ is divided by $\Sigma_k$, the exponent $p_k$ is reduced by one. Therefore the condition $p_1+2p_2+\dots +lp_l=l$ becomes $\hat p_1+2\hat p_2+\dots +l\hat p_l=l-k$ where $\hat p_j=p_j$ for $j\neq k$ and $\hat p_k = p_k-1\,.$ Noting that the term with $l=0$ in the left has no contribution,
 and introducing $s=l-k$, we rewrite the left hand side in the form
\begin{equation*}
\frac{1}{2}\sum_{l=1}^{g-1}\sum_{k=1}^l\frac{L_{l-k}}{(l-k)!(a_j-y_0)^{g-l}(a_i-y_0)^{k+1}}
=\frac{1}{2}\sum_{s=0}^{g-2}\sum_{k=1}^{g-1-l}\frac{L_s}{s!(a_j-y_0)^{g-s-k}(a_i-y_0)^{k+1}}
\end{equation*}
which proves the lemma.
$\Box$

\begin{lemma}
\label{lemma_B}
Let $m$ be a fixed integer, $1\leq m\leq g-1\,.$ The following identities holds
\begin{multline*}
\sum_{\substack{a_j\in B\\a_j\ne u_m}} \frac{\frac{\betaaj}{\Aaj^{12}} - \frac{\betaum}{\Aum^{12}}}{2(a_j-u_m)}
+
\sum_{l=0}^{g-1} \sum_{\substack{p_1+2p_2+\dots +lp_l=l \\ p_1,\dots, p_l\geq 0}}\frac{C_{p_1\dots p_l} \Sigma_1^{p_1}\dots\Sigma_l^{p_l}}{l!\varphi(Q_0)(u_m-y_0)^{g-l}}\left\{\sum_{k=1}^l \frac{kp_k}{\Sigma_k}\frac{1}{(u_m-y_0)^{k+1}} + \frac{g-l-1/2}{u_m-y_0}\right\}
\\
=g\sum_{l=0}^g\frac{L_l}{l!\varphi(Q_0)(u_m-y_0)^{g-l+1}} = \frac{1}{(g-1)!} \frac{\partial^g}{\partial y_0^g} \left\{ \frac{1}{\varphi(Q_0)(u_m-y_0)}\right\}
\,.
\end{multline*}
\end{lemma}
{\it Proof.} First, let us note that the last equality in the lemma is obtained by the Leibniz rule using \eqref{Ll} from Proposition \ref{prop_Ll} for derivatives of $\varphi(Q_0)$. Now, using Lemma \ref{lemma_A} for each term of the sum over the branch points in the left hand side and changing the order of summation over $k$ and over $a_j$, we have in the left hand side
\begin{multline*}
\sum_{l=0}^{g-1}\sum_{\substack{p_1+2p_2+\dots +lp_l=l \\ p_1,\dots, p_l\geq 0}}\frac{C_{p_1\dots p_l} \Sigma_1^{p_1}\dots\Sigma_l^{p_l}}{l!\varphi(Q_0)(u_m-y_0)^{g-l}}
\left\{\sum_{k=1}^l\frac{kp_k}{\Sigma_k}\left(\Sigma_{k+1}-\frac{1}{(u_m-y_0)^{k+1}}\right)
\right.
\\
\left.
-\frac{1}{2}\left( \Sigma_1-\frac{1}{(u_m-y_0)}\right)+ \sum_{k=1}^l\frac{kp_k}{\Sigma_k}\frac{1}{(u_m-y_0)^{k+1}}+
\frac{g-l-\frac{1}{2}}{(u_m-y_0)} \right\}
\\
=\sum_{l=0}^{g-1}\sum_{\substack{p_1+2p_2+\dots +lp_l=l \\ p_1,\dots, p_l\geq 0}}\frac{C_{p_1\dots p_l} \Sigma_1^{p_1}\dots\Sigma_l^{p_l}}{l!\varphi(Q_0)(u_m-y_0)^{g-l}}
\left\{\sum_{k=1}^l\frac{kp_k\Sigma_{k+1}}{\Sigma_k}-\frac{\Sigma_1}{2}+\frac{g-l}{u_m-y_0}\right\}
=\frac{\partial}{\partial y_0} \sum_{l=0}^{g-1}\frac{L_l}{l!\varphi(Q_0)(u_m-y_0)^{g-l}},
\end{multline*}
where the last equality is obtained by straightforward differentiation with respect to $y_0$, using \eqref{phi_Q0} for $\varphi(Q_0)$. The last obtained sum can be converted by the Leibniz rule into the $y_0$-derivative of order $g-1$ of the product $(\varphi(Q_0)(u-y_0))^{-1}$,
which proves the lemma.
$\Box$

\begin{corollary} Let $a_j\in B$ be a branch point of our curve. The following identity holds
\label{lemma_C}
\begin{equation*}
\sum_{l=0}^{g-1} \sum_{\substack{p_1+2p_2+\dots +lp_l=l \\ p_1,\dots, p_l\geq 0}}\frac{C_{p_1\dots p_l} \Sigma_1^{p_1}\dots\Sigma_l^{p_l}}{l!(a_j-y_0)^{g-l-1}}\left\{\sum_{k=1}^l \frac{kp_k\Sigma_{k+1}}{\Sigma_k} -\frac{l}{a_j-y_0} - \frac{\Sigma_1}{2}\right\}=\frac{L_g}{(g-1)!} \,.
\end{equation*}
\end{corollary}

{\it Proof.} This follows by noting that the last line of the last formula in the proof of Lemma \ref{lemma_B} holds after replacing $u_m$  with an arbitrary branch point $a_j$. Then we replace the right hand side in that line by
\begin{equation*}
g\sum_{l=0}^g\frac{L_l}{l!\varphi(Q_0)(a_j-y_0)^{g-l+1}}\,.
\end{equation*}
$\Box$

\begin{lemma} Let $L$ be the polynomials defined by \eqref{L}. Then the following identity holds
\label{lemma_D}
\begin{multline*}
\sum_{k=1}^{r}\frac{L_{r-k}(2\Sigma_1, \dots, 2\Sigma_{r-k})}{k!(r-k)!}
 \sum_{l=0}^{k-1}\binom{k}{l}L_l(-2\Sigma_1,\dots, -2\Sigma_{l})\frac{\partial^{k-l}}{\partial y_0^{k-l}} \prod_{\alpha=1}^{g-1}(y_0-u_{\alpha})
=\frac{1}{r!}\frac{ \partial^{r}}{\partial y_0^{r}}\prod_{\alpha=1}^{g-1}(y_0-u_{\alpha})\,.
\end{multline*}
\end{lemma}

{\it Proof.} This is a corollary of the Leibniz rule \eqref{L-Leibniz}. To see this, note that we can start the sum over $k$ from the value $k=0$. Adding and subtracting the term with $l=k$ and changing the order of summation, we have for the left hand side:
\begin{multline*}
\sum_{l=0}^{r}\sum_{k=l}^r\frac{L_{r-k}(2\Sigma_1, \dots, 2\Sigma_{r-k}) L_l(-2\Sigma_1,\dots, -2\Sigma_{l})\frac{\partial^{k-l}}{\partial y_0^{k-l}} \prod_{\alpha}(y_0-u_{\alpha})}
{l!(r-k)!(k-l)!}
\\
-\sum_{k=0}^r\frac{L_{r-k}(2\Sigma_1, \dots, 2\Sigma_{r-k}) L_k(-2\Sigma_1,\dots, -2\Sigma_{k})}{k!(r-k)!}\,.
\end{multline*}
Note that last sum vanishes due to the Leibniz rule \eqref{L-Leibniz}. Changing the summation index in the remaining sum over $k$ to $s=k-l$, and applying \eqref{L-Leibniz} to the sum over $s$, and then again to the sum over $l$, we prove the lemma.
$\Box$

\begin{lemma} For any natural $N$, some numbers $u_1, \dots, u_N\in\mathbb C$ and $1\leq m\leq N$, the following identity holds
\label{lemma_E}
\begin{equation*}
\sum_{r=1}^{N}\frac{(u_m-y_0)^{r-1}}{r!} \frac{ \partial^{r}}{\partial y_0^{r}}\prod_{\alpha=1}^{N}(y_0-u_{\alpha})
=\frac{\prod_{\alpha=1}^{N}(y_0-u_{\alpha})}{y_0-u_m}\;.
\end{equation*}
\end{lemma}

{\it Proof.}
Let us write the derivative in the following way, separating the terms in two groups depending on whether they contain $u_m$ or not. For $1\leq r\leq N$ we have:
\begin{equation*}
\frac{1}{r!}\frac{\frac{ \partial^{r}}{\partial y_0^{r}}\prod_{\alpha}(y_0-u_{\alpha})}{\prod_{\alpha}(y_0-u_{\alpha})}  =
\sum_{\substack {i_1<i_2<\dots<i_r\\ i_j\ne m}}\frac{1}{\prod_{j=1}^r(y_0-u_{i_j})}+\frac{1}{y_0-u_m}\sum_{\substack{i_1<i_2<\dots<i_{r-1}\\ i_j\ne m}}\frac{1}{\prod_{j=1}^{r-1}(y_0-u_{i_j})}\,.
\end{equation*}
Plugging this into the left hand side of the identity of the lemma and introducing the summation index $r'=r-1$ we complete the proof:
\begin{multline*}
\sum_{r=1}^{N}\frac{(u_m-y_0)^{r-1}}{r!\prod_{\alpha}(y_0-u_{\alpha})} \frac{ \partial^{r}}{\partial y_0^{r}}\prod_{\alpha}(y_0-u_{\alpha})
\\
=\sum_{r=1}^{N-1}\sum_{\substack{i_1<i_2<\dots<i_r\\ i_j\ne m}}\frac{(u_m-y_0)^{r-1}}{\prod_{j=1}^r(y_0-u_{i_j})}- \sum_{r'=0}^{N-1}\sum_{\substack{i_1<i_2<\dots<i_{r'}\\ i_j\ne m}}\frac{(u_m-y_0)^{r'-1}}{\prod_{j=1}^{r'}(y_0-u_{i_j})}=
\frac{1}{y_0-u_m}\,.
\end{multline*}
$\Box$

\begin{lemma} Let $N$ be a natural number. The following identity, where $u_j\in\mathbb C$ are distinct, $\alpha$ ranges in the set $\{1, \dots, N\}$ and $m$ is a fixed element of the same set, holds:
\label{lemma_F}
\begin{equation*}
\sum_{r=0}^{N-1}\frac{(u_m-y_0)^{r}}{r!} \frac{ \partial^{r}}{\partial y_0^{r}}\prod_{\substack{\alpha=1\\\alpha\ne m}}^N(y_0-u_{\alpha}) =\prod_{\substack{\alpha=1\\\alpha\ne m}}^N(u_m-u_{\alpha})\,.
\end{equation*}
\end{lemma}

{\it Proof.} Both sides are  monic polynomials of degree $N-1$ in $u_m$ with the same zeros: $u_m=u_{\alpha}$, $\alpha \ne m$. The fact that the polynomial in the left hand side vanishes when $u_m=u_\alpha$ can be proven by induction on $N$ using the Leibniz rule for the $r$-fold derivative.
Namely, denote the polynomial in the left hand side by $P_N(u_m).$  We have
\begin{equation*}
P_N(u_m)= \sum_{r=0}^{N-1}\frac{(u_m-y_0)^{r}}{r!} \sum_{k=0}^r {r\choose k} \frac{ \partial^{k}}{\partial y_0^{k}}\prod_{\substack{\alpha=1\\\alpha\ne m}}^{N-1}(y_0-u_{\alpha}) \frac{\partial^{r-k}}{\partial y_0^{r-k}}(y_0-u_N)\,.
\end{equation*}
Since the only non zero terms in the sum over $k$ are with $k=r$ and $k=r-1$, we have only two sums over $r$ in the right hand side.
One can then note that in one of the sums, $r$ ranges from $0$ to $N-2$ and in the other, from $1$ to $N-1$ and therefore
\begin{equation*}
P_N(u_m)= (y_0-u_N)P_{N-1}(u_m) +(u_m-y_0)P_{N-1}(u_m)=(u_m-u_N)P_{N-1}(u_m)\,.
\end{equation*}
Thus if $P_{N-1}(u_m)$ vanishes for $u_m=u_\alpha$ with $\alpha=1, \dots, N-1$ then $P_N(u_m)$ vanishes for $u_m=u_\alpha$ with $\alpha=1, \dots, N$. It is easy to check that $P_2(u_m) = \pm(u_1-u_2)$.
$\Box$

\subsubsection{  Proof for $\Aaj^{11}$ with $a_j\notin \{u_1, \dots, u_{g-1}\}$}
Now we want to prove that the functions defined in Theorem \ref{thm_main_g} satisfy equations \eqref{Aaj11}. Rewriting \eqref{Aaj11} in terms of coefficients $\beta_{a_j}$, we have
\begin{multline}
\label{betaaj_system}
\frac{\partial \beta_{a_j}}{\partial x_i} = \frac{\Axi^{12} \frac{\betaaj}{\Aaj^{12}}  - \Aaj^{12} \frac{\betaxi}{\Axi^{12}}  }{2(x_i-a_j)}
 + \frac{g}{2t} \frac{\Axi^{12}\Aaj^{12}}{x_i-a_j}\left( \left( \frac{\betaaj}{\Aaj^{12}}\right)^2 - \left( \frac{\betaxi}{\Axi^{12}}\right)^2\right)
 \\
 + \sum_{\alpha=1}^{g-1} \left[
 \frac{\Aualpha^{12} \frac{\betaaj}{\Aaj^{12}}  - \Aaj^{12} \frac{\betaualpha}{\Aualpha^{12}}  }{2(u_\alpha-a_j)}
 + \frac{g}{2t} \frac{\Aualpha^{12}\Aaj^{12}}{u_\alpha-a_j}\left( \left( \frac{\betaaj}{\Aaj^{12}}\right)^2 - \left( \frac{\betaualpha}{\Aualpha^{12}}\right)^2\right)
 \right]\frac{\partial u_\alpha}{\partial x_i}.
\end{multline}
Thus we need to prove that the quantities from Theorem \ref{thm_main_g} satisfy \eqref{betaaj_system}. Let us now differentiate $\betaaj$ defined by \eqref{betas} with $a_j$ being a branch point different from the dependent branch points $u_1, \dots, u_{g-1}$. We have
\begin{equation}
\label{dbetaaj}
\frac{\partial \beta_{a_j}}{\partial x_i} = \frac{\betaaj}{\Aaj^{12}} \frac{\partial \Aaj^{12}}{\partial x_i} + \Aaj^{12} \frac{\partial}{\partial x_i}\left\{  \frac{1}{\varphi(Q_0)} \sum_{l=0}^{g-1} \frac{L_l}{l!(a_j-y_0)^{g-l}}   \right\} -  \Aaj^{12}\frac{g }{2}\frac{\partial\Omega(P_\infty)}{\partial x_i}\,.
\end{equation}
We have already proved \eqref{Aaj_derivative} which we use now to differentiate $\Aaj^{12}$. The derivative of $\Omega(\Pinfty)$ is given in Proposition \ref{proposition_Omega_infty}. It remains to calculate the derivative of the middle term, which we compute in a straightforward way:
\begin{multline}
\label{middle_term}
\frac{\partial}{\partial x_i}\left\{  \frac{1}{\varphi(Q_0)} \sum_{l=0}^{g-1} \frac{L_l}{l!(a_j-y_0)^{g-l}}   \right\}=
\sum_{l=0}^{g-1}\frac{1}{l!} \sum_{\substack{p_1+2p_2+\dots +lp_l=l \\ p_1,\dots, p_l\geq 0}}\frac{C_{p_1\dots p_l} \Sigma_1^{p_1}\dots\Sigma_l^{p_l}}{\varphi(Q_0)(a_j-y_0)^{g-l}}\times
\\
\times\left[\sum_{k=1}^l \frac{kp_k}{\Sigma_k} \left(  \Sigma_{k+1}\frac{\partial y_0}{\partial x_i}  - \frac{1}{(x_i-y_0)^{k+1}}  - \sum_{\alpha=1}^{g-1} \frac{1}{(u_\alpha-y_0)^{k+1}} \frac{\partial u_\alpha}{\partial x_i} \right)
 +\frac{g-l}{a_j-y_0}\frac{\partial y_0}{\partial x_i}  \right)
 \\
 \left.+ \frac{1}{2} \left( -\Sigma_1 \frac{\partial y_0}{\partial x_i} + \frac{1}{x_i-y_0} + \sum_{\alpha=1}^{g-1}
\frac{1}{u_\alpha-y_0}\frac{\partial u_\alpha}{\partial x_i}\right)
\right]\,.
\end{multline}

We rewrite this using \eqref{dy0-2} for  $\partial y_0/\partial x_i$. We want to prove that the derivative \eqref{dbetaaj} of $\betaaj$ satisfies \eqref{betaaj_system}. Let us thus subtract the right hand side of \eqref{dbetaaj} from the right hand side of \eqref{betaaj_system}. Denoting the result by $T_{a_j}$
and after some simplification, we obtain
\begin{multline}
\label{Taj}
\frac{T_{a_j}}{\Aaj^{12}}=
\frac{ \frac{\betaaj}{\Aaj^{12}}-  \frac{\betaxi}{\Axi^{12}}  }{2(x_i-a_j)}
 - \frac{g}{2t} \frac{\Axi^{12}}{x_i-a_j}\left(  \frac{\betaaj}{\Aaj^{12}} -  \frac{\betaxi}{\Axi^{12}}\right)^2
 + \sum_{\alpha=1}^{g-1} \left[
 \frac{  \frac{\betaaj}{\Aaj^{12}}-  \frac{\betaualpha}{\Aualpha^{12}}  }{2(u_\alpha-a_j)}
 - \frac{g}{2t} \frac{\Aualpha^{12}}{u_\alpha-a_j}\left(  \frac{\betaaj}{\Aaj^{12}} - \frac{\betaualpha}{\Aualpha^{12}}\right)^2
 \right]\frac{\partial u_\alpha}{\partial x_i}
\\
+ \sum_{l=0}^{g-1}\frac{1}{l!} \sum_{\substack{p_1+2p_2+\dots +lp_l=l \\ p_1,\dots, p_l\geq 0}}\frac{C_{p_1\dots p_l} \Sigma_1^{p_1}\dots\Sigma_l^{p_l}}{\varphi(Q_0)(a_j-y_0)^{g-l}}\left( \sum_{k=1}^l \frac{kp_k}{\Sigma_k}\left[\frac{\Sigma_{k+1}\Axi^{12}}{t\varphi(Q_0)(x_i-y_0)^g} + \frac{1}{(x_i-y_0)^{k+1}}\right]\right.
\\
\left. +\frac{g-l}{a_j-y_0} \frac{\Axi^{12}}{t\varphi(Q_0)(x_i-y_0)^g} - \frac{\Axi^{12}\Sigma_1}{2t\varphi(Q_0)(x_i-y_0)^g} - \frac{1}{2(x_i-y_0)}\right)
\\
 + \sum_{\substack{\alpha=1}}^{g-1}\frac{\partial u_\alpha}{\partial x_i}\sum_{l=0}^{g-1}\frac{1}{l!} \sum_{\substack{p_1+2p_2+\dots +lp_l=l \\ p_1,\dots, p_l\geq 0}}\frac{C_{p_1\dots p_l} \Sigma_1^{p_1}\dots\Sigma_l^{p_l}}{\varphi(Q_0)(a_j-y_0)^{g-l}}\left( \sum_{k=1}^l \frac{kp_k}{\Sigma_k} \left[  \frac{\Sigma_{k+1}\Aualpha^{12}}{t\varphi(Q_0)(u_\alpha-y_0)^g} + \frac{1}{(u_\alpha-y_0)^{k+1}}\right] \right.
\\
\left. +\frac{g-l}{a_j-y_0} \frac{\Aualpha^{12}}{t\varphi(Q_0)(u_\alpha-y_0)^g} - \frac{\Aualpha^{12}\Sigma_1}{2t\varphi(Q_0)(u_\alpha-y_0)^g} - \frac{1}{2(u_\alpha-y_0)}\right)
- \frac{g }{2t} \frac{\Axi^{12}\prod_{\alpha=1}^{g-1}(x_i-u_\alpha)}{(x_i-y_0)^{g-1}}\,.
\end{multline}
Applying Lemma \ref{lemma_A} twice,  with $a_i=x_i$ and with $a_i=u_\alpha$, and then Lemma \ref{lemma_C} twice as well, we have
\begin{multline*}
\frac{t}{g}\frac{T_{a_j}}{\Aaj^{12}}=
 -  \frac{\Axi^{12}}{2(x_i-a_j)}\left(  \frac{\betaaj}{\Aaj^{12}} -  \frac{\betaxi}{\Axi^{12}}\right)^2
 -  \sum_{\alpha=1}^{g-1}
  \frac{\Aualpha^{12}}{2(u_\alpha-a_j)}\left(  \frac{\betaaj}{\Aaj^{12}} - \frac{\betaualpha}{\Aualpha^{12}}\right)^2
\frac{\partial u_\alpha}{\partial x_i}
 -  \frac{\Axi^{12}\prod_{\alpha=1}^{g-1}(x_i-u_\alpha)}{2(x_i-y_0)^{g-1}}
\\
+ \frac{\Axi^{12}}{(x_i-y_0)^g}\sum_{l=0}^g\frac{L_l}{l!\varphi^2(Q_0)(a_j-y_0)^{g-l+1}}
 + \sum_{\substack{\alpha=1}}^{g-1}\frac{\partial u_\alpha}{\partial x_i}\frac{\Aualpha^{12}}{(u_\alpha-y_0)^g}\sum_{l=0}^g\frac{L_l}{l!\varphi^2(Q_0)(a_j-y_0)^{g-l+1}}
\,.
\end{multline*}
Plugging in \eqref{umder_A} for the derivatives of $u_\alpha$ and using \eqref{rat1} from Lemma \ref{lemma_rational1}, we arrive at
\begin{multline*}
\frac{t}{g}\frac{T_{a_j}}{\Aaj^{12}\Axi^{12}}=
 -  \frac{1}{2(x_i-a_j)}\left(  \frac{\betaaj}{\Aaj^{12}} -  \frac{\betaxi}{\Axi^{12}}\right)^2
  -  \frac{\prod_{\alpha=1}^{g-1}(x_i-u_\alpha)}{2(x_i-y_0)^{g-1}}
 \\
  +  \frac{\prod_{\beta=1}^{g-1}(x_i-u_\beta)}{2(x_i-y_0)^{g-1}}\sum_{\alpha=1}^{g-1}
\frac{(u_\alpha-y_0)^{g-1}}{(u_\alpha-a_j)(x_i-u_\alpha)\prod_{\beta\neq \alpha}(u_\alpha-u_\beta)}\left(  \frac{\betaaj}{\Aaj^{12}} - \frac{\betaualpha}{\Aualpha^{12}}\right)^2
\\
+ \frac{\prod_{\beta=1}^{g-1}(x_i-u_\beta)}{(x_i-y_0)^{g}\prod_{\alpha=1}^{g-1}(y_0-u_\alpha)}\sum_{l=0}^g\frac{L_l}{l!\varphi^2(Q_0)(a_j-y_0)^{g-l+1}}
\,.
\end{multline*}
Recall that we want to prove that this expression vanishes identically. As a next step, we need to compute the squares in the first and second lines.
Let us first rewrite such a square using \eqref{difference_betas_1} and \eqref{difference_betas_2} as follows:
\begin{multline}
\label{square}
\left( \frac{\betaai}{\Aai^{12}} - \frac{\betaaj}{\Aaj^{12}}\right)^2
= \frac{a_j-a_i}{\varphi^2(Q_0)} \sum_{l_1=0}^{g-1} \left( \frac{L_{l_1}}{l_1!(a_i-y_0)^{g-l_1}} - \frac{L_{l_1}}{l_1!(a_j-y_0)^{g-l_1}} \right)
\\
\times
\sum_{l_2=0}^{g-1}  \frac{L_{l_2}}{l_2!}\sum_{k=0}^{g-l_2-1}\frac{1}{(a_i-y_0)^{g-l_2-k}(a_j-y_0)^{k+1}}\,.
\end{multline}
This yields
\begin{multline*}
\frac{t}{g}\frac{\varphi^2(Q_0)T_{a_j}}{\Aaj^{12}\Axi^{12}}=
   \frac{1}{2} \sum_{l_1, l_2=0}^{g-1}\frac{L_{l_1}L_{l_2}}{l_1!l_2!} \sum_{k=0}^{g-l_2-1} \frac{1}{(x_i-y_0)^{2g-l_1-l_2-k}(a_j-y_0)^{k+1}}
   \\
   -\frac{1}{2} \sum_{l_1, l_2=0}^{g-1}\frac{L_{l_1}L_{l_2}}{l_1!l_2!} \sum_{k=0}^{g-l_2-1} \frac{1}{(a_j-y_0)^{g-l_1+k+1}(x_i-y_0)^{g-l_2-k}}
 \\
  +  \frac{\prod_{\beta=1}^{g-1}(x_i-u_\beta)}{2(x_i-y_0)^{g-1}}\sum_{l_1,l_2=0}^{g-1}\frac{L_{l_1}L_{l_2}}{l_1!l_2!}\sum_{k=0}^{g-l_2-1}\frac{1}{(a_j-y_0)^{g-l_1+k+1}}\sum_{\alpha=1}^{g-1}
\frac{(u_\alpha-y_0)^{l_2+k-1}}{(x_i-u_\alpha)\prod_{\beta\neq \alpha}(u_\alpha-u_\beta)}
\\
  -  \frac{\prod_{\beta=1}^{g-1}(x_i-u_\beta)}{2(x_i-y_0)^{g-1}}\sum_{l_1,l_2=0}^{g-1}\frac{L_{l_1}L_{l_2}}{l_1!l_2!}\sum_{k=0}^{g-l_2-1}\frac{1}{(a_j-y_0)^{k+1}}\sum_{\alpha=1}^{g-1}
\frac{(u_\alpha-y_0)^{l_1+l_2+k-g-1}}{(x_i-u_\alpha)\prod_{\beta\neq \alpha}(u_\alpha-u_\beta)}
\\
+ \frac{\prod_{\beta=1}^{g-1}(x_i-u_\beta)}{(x_i-y_0)^{g}\prod_{\alpha=1}^{g-1}(y_0-u_\alpha)}\sum_{l=0}^g\frac{L_l}{l!(a_j-y_0)^{g-l+1}}
-  \varphi^2(Q_0)\frac{\prod_{\alpha=1}^{g-1}(x_i-u_\alpha)}{2(x_i-y_0)^{g-1}}\,.
\end{multline*}
It remains to compute the sums over $\alpha$ in the third and fourth lines by applying Lemma \ref{lemma_rational1} and Corollary \ref{corollary_rational2}.
This is a tedious but straightforward calculation. To apply the identities from Lemma \ref{lemma_rational1} and Corollary \ref{corollary_rational2}, we need to split the sums according to the value of exponents:
We split the sum in the next to the last line by introducing $n=l_2+k$ and $N=l_1+\underbrace{l_2+k}_\text{$n$}$ as follows
\begin{multline}
\label{split}
\sum_{l_1,l_2=0}^{g-1}\sum_{k=0}^{g-l_2-1} F(l_1,l_2, k) =\sum_{n=0}^{g-1}\sum_{l_1=0}^{g-1}\sum_{l_2=0}^{n} F(l_1,l_2, n-l_2)
\\
 = \underbrace{\sum_{N=0}^{g-1}\sum_{n=0}^{N} \sum_{l_2=0}^{n} F(N-n,l_2, n-l_2)}_\text{$-g-1\leq l_1+l_2+k-g-1\leq -2$} + \underbrace{\sum_{n=1}^{g-1} \sum_{l_2=0}^{n} F(g-n,l_2, n-l_2)}_\text{$N=g,\; l_1+l_2+k-g-1=-1$}
 +\underbrace{\sum_{N=g+1}^{2g-2} \sum_{n=N-g+1}^{g-1} \sum_{l_2=0}^{n} F(N-n,l_2, n-l_2)}_\text{$0\leq l_1+l_2+k-g-1\leq g-3$}\,.
\end{multline}
It is convenient to also perform the same splitting in the very first sum of our expression.
This yields
\begin{multline*}
\frac{t}{g}\frac{\varphi^2(Q_0)T_{a_j}}{\Aaj^{12}\Axi^{12}}=
   {\black \frac{\prod_{\beta=1}^{g-1}(x_i-u_\beta)}{2(x_i-y_0)^{g}\prod_{\alpha=1}^{g-1}(y_0-u_\alpha)}\sum_{l_1=0}^{g}\frac{L_{l_1}}{l_1!(a_j-y_0)^{g-l_1+1}} }
   +{\black \frac{\prod_{\beta=1}^{g-1}(x_i-u_\beta)}{2(x_i-y_0)^{g}\prod_{\alpha=1}^{g-1}(y_0-u_\alpha)}\frac{L_g}{g!(a_j-y_0)} }
\\
  + \frac{\prod_{\beta=1}^{g-1}(x_i-u_\beta)}{2(x_i-y_0)^{g-1}}\sum_{N=0}^{g-1}\sum_{n=0}^{N}\sum_{l_2=0}^{n}\frac{L_{N-n}L_{l_2}}{(N-n)!l_2!(a_j-y_0)^{n-l_2+1}}\frac{1}{(g-N)!}\frac{\partial^{g-N}}{\partial y_0^{g-N}}
\left\{\frac{1}{(x_i-y_0)\prod_{ \alpha=1}^{g-1}(y_0-u_\alpha)} \right\}
\\
  +  \frac{\prod_{\beta=1}^{g-1}(x_i-u_\beta)}{2(x_i-y_0)^{g}\prod_{\alpha=1}^{g-1}(y_0-u_\alpha)}\sum_{n=1}^{g-1}\sum_{l_2=0}^{n}\frac{L_{g-n}L_{l_2}}{(g-n)!l_2!(a_j-y_0)^{n-l_2+1}}
-  \varphi^2(Q_0)\frac{\prod_{\alpha=1}^{g-1}(x_i-u_\alpha)}{2(x_i-y_0)^{g-1}}\,.
\end{multline*}
In the line with derivatives, let us represent $\frac{1}{(a_j-y_0)^{n-l_2+1}}$ as an $n-l_2$-fold derivative with respect to $y_0$ and use the Leibniz rule, together with \eqref{Ll} to do the sum over $l_2$ and then the sum over $n$ and the sum over $N$. In the sum over $N$ we need to add and subtract the term with $N=g$ before being able to apply the Leibniz rule. At the same time, let us use the same technique in the first and the last remaining sums of the expression. This brings the second line to the form containing the following derivative:
\begin{equation*}
\frac{\partial^{g}}{\partial y_0^{g}}
\left\{\frac{1}{\varphi^2(Q_0)(a_j-y_0)(x_i-y_0)\prod_{ \alpha=1}^{g-1}(y_0-u_\alpha)} \right\}
\end{equation*}
Note that due to \eqref{phi_Q0} this is the $g$-fold derivative of a polynomial of degree $g$ and therefore is equal to $g!$. The polynomial in question is the product of all the factors of the form $(y_0-a_k)$ where $a_k$ runs in the set of all branch points except the points $a_j, x_i, u_1, \dots, u_{g-1}\,.$
Adding and subtracting a top term in some of the sums if necessary and  then using the Leibniz rule again to transform the sums, we see that our expression vanishes.
This proves that $T_{a_j}$ from \eqref{Taj} is zero and finishes the proof of the fact that  quantities $\Aaj^{11}$ satisfy equations \eqref{Aaj11}.

\subsubsection{  Proof for $\Aaj^{11}$ with $a_j=u_m$}

Here we prove that the functions given in Theorem \ref{thm_main_g} satisfy equations \eqref{Aum11}. It is straightforward to see that \eqref{Aum11} with $\Aaj^{21}, \Aaj^{12}$ and $\Aaj^{11}$ given by Theorem \ref{thm_main_g} is equivalent to
\begin{multline}
\label{dbetaum_system}
\frac{\partial\betaum}{\partial x_i} =  \frac{\Axi^{12} \frac{\betaum}{\Aum^{12}} - \Aum^{12}\frac{\betaxi}{\Axi^{12}}}{2(x_i-u_m)} +  \frac{g\Aum^{12}\Axi^{12}\left(  \left( \frac{\betaum}{\Aum^{12}}\right)^2- \left( \frac{\betaxi}{\Axi^{12}}\right)^2 \right)}{2t(x_i-u_m)}
\\
-\frac{\partial u_m}{\partial x_i} \sum_{\substack{a_j\in B\\a_j\ne u_m}}\left( \frac{\Aaj^{12} \frac{\betaum}{\Aum^{12}} - \Aum^{12}\frac{\betaaj}{\Aaj^{12}}}{2(a_j-u_m)} +  \frac{g\Aum^{12}\Aaj^{12}\left(  \left( \frac{\betaum}{\Aum^{12}}\right)^2- \left( \frac{\betaaj}{\Aaj^{12}}\right)^2 \right)}{2t(a_j-u_m)}\right)
\\
+\sum_{\substack{\alpha=1\\ \alpha\neq m}}^{g-1} \frac{\partial u_\alpha}{\partial x_i} \left(  \frac{\Aualpha^{12} \frac{\betaum}{\Aum^{12}} - \Aum^{12}\frac{\betaualpha}{\Aualpha^{12}}}{2(u_\alpha-u_m)} +  \frac{g\Aum^{12}\Aualpha^{12}\left(  \left( \frac{\betaum}{\Aum^{12}}\right)^2- \left( \frac{\betaualpha}{\Aualpha^{12}}\right)^2 \right)}{2t(u_\alpha-u_m)}\right)\,.
\end{multline}

Let us now differentiate $\betaum$ given in the form \eqref{betas}. We have
\begin{equation}
\label{dbetaum}
\frac{\partial \beta_{u_m}}{\partial x_i} = \frac{\betaum}{\Aum^{12}} \frac{\partial \Aum^{12}}{\partial x_i} + \Aum^{12} \frac{\partial}{\partial x_i}\left\{  \frac{1}{\varphi(Q_0)} \sum_{l=0}^{g-1} \frac{L_l}{l!(u_m-y_0)^{g-l}}   \right\} -  \Aum^{12}\frac{g }{2}\frac{\partial\Omega(P_\infty)}{\partial x_i}\,.
\end{equation}
We use \eqref{Aum_derivative} for the derivative of $\Aum^{12}$ and Proposition \ref{proposition_Omega_infty} for that of $\Omega(P_\infty)\,.$ It remains to calculate the derivative of the middle term, which we compute in a straightforward way:
\begin{multline*}
\frac{\partial}{\partial x_i}\left\{  \frac{1}{\varphi(Q_0)} \sum_{l=0}^{g-1} \frac{L_l}{l!(u_m-y_0)^{g-l}}   \right\}=
\sum_{l=0}^{g-1}\frac{1}{l!} \sum_{\substack{p_1+2p_2+\dots +lp_l=l \\ p_1,\dots, p_l\geq 0}}\frac{C_{p_1\dots p_l} \Sigma_1^{p_1}\dots\Sigma_l^{p_l}}{\varphi(Q_0)(u_m-y_0)^{g-l}}\times
\\
\times\left[\sum_{k=1}^l \frac{kp_k}{\Sigma_k} \left(  \Sigma_{k+1}\frac{\partial y_0}{\partial x_i}  - \frac{1}{(x_i-y_0)^{k+1}}  - \sum_{\alpha=1}^{g-1} \frac{1}{(u_\alpha-y_0)^{k+1}} \frac{\partial u_\alpha}{\partial x_i} \right)
 -\frac{g-l}{u_m-y_0}\left(\frac{\partial u_m}{\partial x_i} - \frac{\partial y_0}{\partial x_i}  \right) \right.
 \\
\left. + \frac{1}{2} \left( -\Sigma_1 \frac{\partial y_0}{\partial x_i} + \frac{1}{x_i-y_0} + \sum_{\alpha=1}^{g-1}
\frac{1}{u_\alpha-y_0}\frac{\partial u_\alpha}{\partial x_i}\right)
\right]\,.
\end{multline*}
Now, putting all these together and using \eqref{dy0-2} for  $\partial y_0/\partial x_i$, we  obtain the derivative of $\betaum$ which we need to compare to the formula \eqref{dbetaum_system} imposed by the Schlesinger system on $\betaum.$ Thus, we need to prove that the difference between \eqref{dbetaum} and  \eqref{dbetaum_system}  vanishes. Denote by $T_{u_m}$ the difference between \eqref{dbetaum} and  \eqref{dbetaum_system} after dividing by the common factor of $\Aum^{12}\,$:
\begin{multline*}
T_{u_m} = \frac{\frac{\betaxi}{\Axi^{12}} - \frac{\betaum}{\Aum^{12}}}{2(x_i-u_m)} + \frac{g\Axi^{12}}{2t(x_i-u_m)} \left( \frac{\betaum}{\Aum^{12}} - \frac{\betaxi}{\Axi^{12}}\right)^2
+\frac{g}{2}\frac{\Axi^{12}\prod_{\alpha=1}^{g-1}(x_i-u_\alpha)}{t(x_i-y_0)^{g-1}}
\\
-\frac{\partial u_m}{\partial x_i} \left[ \sum_{\substack{a_j\in B\\ a_j\neq u_m}} \left( \frac{\frac{\betaaj}{\Aaj^{12}} - \frac{\betaum}{\Aum^{12}}}{2(a_j-u_m)} + \frac{g\Aaj^{12}}{2t(a_j-u_m)} \left( \frac{\betaum}{\Aum^{12}} - \frac{\betaaj}{\Aaj^{12}}\right)^2 \right)\right.
\\
\left.+\sum_{l=0}^{g-1} \sum_{\substack{p_1+2p_2+\dots +lp_l=l \\ p_1,\dots, p_l\geq 0}}\frac{C_{p_1\dots p_l} \Sigma_1^{p_1}\dots\Sigma_l^{p_l}}{l!\varphi(Q_0)(u_m-y_0)^{g-l}}\left( \frac{g-l}{u_m-y_0} + \sum_{k=1}^l \frac{kp_k}{\Sigma_k} \left[  \frac{\Sigma_{k+1}\Aum^{12}}{t\varphi(Q_0)(u_m-y_0)^g} + \frac{1}{(u_m-y_0)^{k+1}}\right] \right.\right.
\\
\left.\left. + \frac{(g-l)\Aum^{12}}{t\varphi(Q_0)(u_m-y_0)^{g+1}} - \frac{\Aum^{12}\Sigma_1}{2t\varphi(Q_0)(u_m-y_0)^g} - \frac{1}{2(u_m-y_0)}\right)\right]
\\
+\sum_{\substack{\alpha=1\\ \alpha\neq m}}^{g-1}\frac{\partial u_\alpha}{\partial x_i} \left[  \frac{\frac{\betaualpha}{\Aualpha^{12}} - \frac{\betaum}{\Aum^{12}}}{2(u_\alpha-u_m)} + \frac{g\Aualpha^{12}}{2t(u_\alpha-u_m)} \left( \frac{\betaum}{\Aum^{12}} - \frac{\betaualpha}{\Aualpha^{12}}\right)^2\right.
\\
\left. - \sum_{l=0}^{g-1} \sum_{\substack{p_1+2p_2+\dots +lp_l=l \\ p_1,\dots, p_l\geq 0}}\frac{C_{p_1\dots p_l} \Sigma_1^{p_1}\dots\Sigma_l^{p_l}}{l!\varphi(Q_0)(u_m-y_0)^{g-l}}\left( \sum_{k=1}^l \frac{kp_k}{\Sigma_k} \left[  \frac{\Sigma_{k+1}\Aualpha^{12}}{t\varphi(Q_0)(u_\alpha-y_0)^g} + \frac{1}{(u_\alpha-y_0)^{k+1}}\right] \right.\right.
\\
\left.\left. +\frac{g-l}{u_m-y_0} \frac{\Aualpha^{12}}{t\varphi(Q_0)(u_\alpha-y_0)^g} - \frac{\Aualpha^{12}\Sigma_1}{2t\varphi(Q_0)(u_\alpha-y_0)^g} - \frac{1}{2(u_\alpha-y_0)}\right)\right]
\\
-\sum_{l=0}^{g-1} \sum_{\substack{p_1+2p_2+\dots +lp_l=l \\ p_1,\dots, p_l\geq 0}}\frac{C_{p_1\dots p_l} \Sigma_1^{p_1}\dots\Sigma_l^{p_l}}{l!\varphi(Q_0)(u_m-y_0)^{g-l}}\left( \sum_{k=1}^l \frac{kp_k}{\Sigma_k}\left[\frac{\Sigma_{k+1}\Axi^{12}}{t\varphi(Q_0)(x_i-y_0)^g} + \frac{1}{(x_i-y_0)^{k+1}}\right]\right.
\\
\left. +\frac{g-l}{u_m-y_0} \frac{\Axi^{12}}{t\varphi(Q_0)(x_i-y_0)^g} - \frac{\Axi^{12}\Sigma_1}{2t\varphi(Q_0)(x_i-y_0)^g} - \frac{1}{2(x_i-y_0)}\right)\,.
\end{multline*}
We are gong to prove that $T_{u_m}=0$. Now we apply Lemma \ref{lemma_A} twice to reduce the above expression for $T_{u_m}$:  first with $a_i=x_i$ and $a_j=u_m$ to cancel the terms that do not have a factor of $A^{12}$ nor a derivative of any $u_k$, and then with $a_i=u_\alpha$ and $a_j=u_m$ to cancel some terms that multiply $\partial u_\alpha/\partial x_i\,.$ We then apply Lemma \ref{lemma_B} to group three terms multiplying $\partial u_m/\partial x_i\,.$ Corollary \ref{lemma_C} can be used to cancel some terms multiplying $\partial u_m/\partial x_i$ that have a factor of $\Aum^{12}$, then similarly in the factor of each $\partial u_\alpha/\partial x_i$ for the terms that have a factor of $\Aualpha^{12}$, and finally Corollary \ref{lemma_C} is used analogously in the last two lines for the terms having $\Axi^{12}$ as a factor. As a result of these cancellations, we have
\begin{multline*}
\frac{T_{u_m}t}{g} =   \frac{\Axi^{12}}{2(x_i-u_m)} \left( \frac{\betaum}{\Aum^{12}} - \frac{\betaxi}{\Axi^{12}}\right)^2
- \frac{\Axi^{12}}{\varphi^2(Q_0)(x_i-y_0)^g} \sum_{l=0}^{g} \frac{L_l}{l!(u_m-y_0)^{g-l+1}}
\\
-\frac{\partial u_m}{\partial x_i} \left[ \sum_{\substack{a_j\in B\\ a_j\neq u_m}}   \frac{\Aaj^{12}}{2(a_j-u_m)} \left( \frac{\betaum}{\Aum^{12}} - \frac{\betaaj}{\Aaj^{12}}\right)^2 +\frac{\Aum^{12}}{\varphi^2(Q_0)}\sum_{l=0}^{g} \frac{L_l}{l!(u_m-y_0)^{2g-l+1}}
\right.
\\ \left.+\sum_{l=0}^g\frac{tL_l}{l!\varphi(Q_0)(u_m-y_0)^{g-l+1}}\right]
+\frac{\Axi^{12}\prod_{\alpha=1}^{g-1}(x_i-u_\alpha)}{2(x_i-y_0)^{g-1}}
\\
+\sum_{\substack{\alpha=1\\ \alpha\neq m}}^{g-1}\frac{\partial u_\alpha}{\partial x_i} \left[   \frac{\Aualpha^{12}}{2(u_\alpha-u_m)} \left( \frac{\betaum}{\Aum^{12}} - \frac{\betaualpha}{\Aualpha^{12}}\right)^2 -  \frac{\Aualpha^{12}}{\varphi^2(Q_0)(u_\alpha-y_0)^g}\sum_{l=0}^{g} \frac{L_l}{l!(u_m-y_0)^{g-l+1}}\right]\,.
\end{multline*}

Our next step is to substitute \eqref{umder_A} for derivatives $\partial u_\alpha/\partial x_i$ with $\alpha\neq m$ and use \eqref{rat1} from Lemma \ref{lemma_rational1} in the obtained expression. This yields
\begin{multline}
\label{M10}
\frac{tT_{u_m}}{g} =   \frac{\Axi^{12}}{2(x_i-u_m)} \left( \frac{\betaum}{\Aum^{12}} - \frac{\betaxi}{\Axi^{12}}\right)^2
-\frac {A_{x_i}^{12}\prod _{\alpha} (x_i-u_{\alpha})}{2(x_i-y_0)^{g-1}}\sum_{\substack{\alpha=1\\ \alpha\neq m}}^{g-1}\frac {(u_\alpha-y_0)^{g-1}\left( \frac{\betaum}{\Aum^{12}} - \frac{\betaualpha}{\Aualpha^{12}}\right)^2}{(u_\alpha-u_m)(x_i-u_{\alpha})\prod _{\beta \neq \alpha}(u_\alpha-u_{\beta})}
\\
-\frac{\partial u_m}{\partial x_i} \left[ \sum_{\substack{a_j\in B\\ a_j\neq u_m}}   \frac{\Aaj^{12}}{2(a_j-u_m)} \left( \frac{\betaum}{\Aum^{12}} - \frac{\betaaj}{\Aaj^{12}}\right)^2 +\frac{\Aum^{12}}{\varphi^2(Q_0)}\sum_{l=0}^{g} \frac{L_l}{l!(u_m-y_0)^{2g-l+1}}
\right.
\\
\left.
+\sum_{l=0}^g\frac{tL_l}{l!\varphi(Q_0)(u_m-y_0)^{g-l+1}}\right]
-\frac {A_{x_i}^{12}\prod _{\alpha} (x_i-u_{\alpha})}{\varphi^2(Q_0)(x_i-y_0)^{g-1}} \left(
\frac {1 }{(x_i-y_0)\prod _\alpha(y_0-u_\alpha)}
\right.
\\
\left.
+\frac {1}{(x_i-u_m)(u_m-y_0)\prod _{\alpha \neq m}(u_m-u_{\alpha})}\right)\sum_{l=0}^{g} \frac{L_l}{l!(u_m-y_0)^{g-l+1}}
+\frac{\Axi^{12}\prod_{\alpha}(x_i-u_\alpha)}{2(x_i-y_0)^{g-1}}
\,.
\end{multline}
Let us now look at the  factor of $\frac{\partial u_m}{\partial x_i}\,.$ The sum over the branch points in this factor can be evaluated with the help of Lemma \ref{lemma_residues}.

\begin{proposition}
\label{proposition_factor_of_dumdxi}
The sum of two terms multiplying $-\frac{\partial u_m}{\partial x_i}$ in the right-hand side expression of \eqref{M10} for $\frac{tT_{u_m}}{g}\,$ is:
\begin{equation*}
 \sum_{\substack{a_j\in B\\ a_j\neq u_m}}   \frac{\Aaj^{12}}{2(a_j-u_m)} \left( \frac{\betaum}{\Aum^{12}}
- \frac{\betaaj}{\Aaj^{12}}\right)^2
+\sum_{l=0}^g\frac{tL_l}{l!\varphi(Q_0)(u_m-y_0)^{g-l+1}}
=\frac{2\Aum^{12}}{\varphi^2(P_{u_m})(u_m-y_0)^{2g}}\,.
\end{equation*}
\end{proposition}
{\it Proof.}
Rewriting the square as in \eqref{square}, we have for the sum over the branch points
\begin{multline*}
\sum_{\substack{a_j\in B\\ a_j\neq u_m}}   \frac{\Aaj^{12}}{2(a_j-u_m)} \left( \frac{\betaum}{\Aum^{12}}
- \frac{\betaaj}{\Aaj^{12}}\right)^2
\\
= \frac{1}{2\varphi^2(Q_0)}
\sum_{l_1, l_2=0}^{g-1} \frac{L_{l_1}L_{l_1}}{l_1! l_2!} \sum_{k=0}^{g-l_2-1}\left( \frac{\sum_{\substack{a_j\in B\\ a_j\neq u_m}} \frac{\Aaj^{12}}{(a_j-y_0)^{k+1}}}{(u_m-y_0)^{2g-l_1-l_2-k}}
- \frac{\sum_{\substack{a_j\in B\\ a_j\neq u_m}} \frac{\Aaj^{12}}{(a_j-y_0)^{g-l_1+k+1}}}{(u_m-y_0)^{g-l_2-k}}
\right)\,.
\end{multline*}
Now let us split the sums according to the value of the exponents of the factor $(a_j-y_0)$ in a way to allow for applying the identities from Lemma \ref{lemma_residues}.
After a lengthy but straightforward calculation, using, in particular,
the following change of order of summation:
\begin{equation*}
\sum_{k=0}^{g-1}\sum_{l_1=0}^{k}\sum_{l_2=0}^{g-k-1}F(l_1, l_2)
=\sum_{l_1=0}^{g-1}\sum_{l_2=0}^{g-l_1-1}(g-l_1-l_2)F(l_1, l_2)
\end{equation*}
and
\begin{equation*}
\sum_{k=0}^{g-1}\sum_{l_1=k+2}^{g-1}\sum_{l_2=0}^{g-k-1}F(l_1, l_2)
=
\sum_{l_1=2}^{g-1}\sum_{l_2=0}^{g-1}(g-l_2)F(l_1, l_2) - \sum_{l_1=2}^{g-1}\sum_{l_2=0}^{g-l_1}(g-l_1-l_2+1)F(l_1, l_2)\,,
\end{equation*}
 this yields
\begin{multline}
\label{temp_factor_of_dumdxi}
\sum_{\substack{a_j\in B\\ a_j\neq u_m}}   \frac{\Aaj^{12}}{2(a_j-u_m)} \left( \frac{\betaum}{\Aum^{12}}
- \frac{\betaaj}{\Aaj^{12}}\right)^2
+\sum_{l=0}^g\frac{tL_l}{l!\varphi(Q_0)(u_m-y_0)^{g-l+1}}
\\
=
\frac{t (u_m-y_0)^{1-g}}{2\varphi(Q_0)}
\left[
\sum_{k=0}^{g-1}\sum_{l_1=0}^{k}\sum_{l_2=0}^{g-k-1}\frac{L_{l_1}L_{l_2}}{l_1!l_2!}  \frac{\underset{P=Q_0}{\rm res}\frac{\Omega(P)\varphi(P)}{(u-y_0)^{k-l_1+1}du}}{\varphi(Q_0)(u_m-y_0)^{1-l_2-k}}
 +
\sum_{l_1=0}^{g}\sum_{l_2=0}^{g-l_1}\frac{L_{l_1}L_{l_2}}{l_1!l_2!(u_m-y_0)^{2-l_2-l_1}}
\right]
\\
=:\frac{t (u_m-y_0)^{1-g}}{2\varphi(Q_0)}S.
\end{multline}
The last equality should be seen as a definition of the notation $S$.
Let us now study the term containing the residue. Writing the differential $\Omega$ as in \eqref{Omega-delta} with explicit expressions \eqref{v} for the basis of holomorphic differentials $v_i$ and using \eqref{Ll} for derivatives of $\varphi(Q_0)$ with respect to $y_0$, we have for the residue with $1\leq s\leq g$:
\begin{multline}
\label{residue}
\underset{P=Q_0}{\rm res}\frac{\Omega(P)\varphi(P)}{(u-y_0)^sdu} = R_1(s)+R_2(s)+R_3(s):=
\frac{1}{s!\varphi(Q_0)} \frac{\partial^s \varphi^2(Q_0)}{\partial y_0^s}
\\
+
\frac{1}{(s-2)!} \sum_{i=1}^{g-1} \delta_i
\frac{\frac{\partial^{s-2}}{\partial y_0^{s-2}}\{\varphi^2(Q_0) \prod_{\alpha\neq i}(y_0-u_\alpha)\}}{\varphi(P_{u_i})(u_i-y_0)\prod_{\alpha\neq i}(u_i-u_\alpha)} H(s- 2)
\\
+\frac{\delta_g}{(s-1)!} \frac{\frac{\partial^{s-1}}{\partial y_0^{s-1}}\{\varphi^2(Q_0) \prod_{\alpha}(y_0-u_\alpha)\}}{\varphi(Q_0)\prod_{\alpha}(y_0-u_\alpha)}\,,
\end{multline}
where $H$ is the Heaviside step function, that is the middle term is only present for $s\geq 2$.

Now considering the right hand side of \eqref{temp_factor_of_dumdxi}, we single out the term with $l_1=l_2=0$ in the last double sum and change the order of summation to obtain
\begin{multline*}
S:=\sum_{k=0}^{g-1}\sum_{l_1=0}^{k}\sum_{l_2=0}^{g-k-1}\frac{L_{l_1}L_{l_2}}{l_1!l_2!}  \frac{\underset{P=Q_0}{\rm res}\frac{\Omega(P)\varphi(P)}{(u-y_0)^{k-l_1+1}du}}{\varphi(Q_0)(u_m-y_0)^{1-l_2-k}}
 +
\sum_{l_1=0}^{g}\sum_{l_2=0}^{g-l_1}\frac{L_{l_1}L_{l_2}}{l_1!l_2!(u_m-y_0)^{2-l_2-l_1}}
\\
= \frac{1}{(u_m-y_0)^2} + \sum_{r=0}^{g-1} (u_m-y_0)^{r-1} \left(
\sum_{l_1+l_2=r+1}\frac{L_{l_1}L_{l_2}}{l_1!l_2!} + \frac{1}{\varphi(Q_0)} \sum_{s=0}^r\sum_{l_1+l_2=r-s}\frac{L_{l_1}L_{l_2}}{l_1!l_2!}
\underset{P=Q_0}{\rm res}\frac{\Omega(P)\varphi(P)}{(u-y_0)^{s+1}du}
\right).
\end{multline*}
Due to Proposition \ref{prop_Ll} and the Leibniz rule \eqref{L-Leibniz},  the sum over $l_1+l_2=r+1$ is equal to $\frac{\varphi^2(Q_0)}{(r+1)!} \partial^{r+1}_{y_0} \varphi^{-2}(Q_0)$  and similarly for the other sum over the values of $l_1+l_2\,:$
\begin{multline*}
S
= \frac{1}{(u_m-y_0)^2}
+ \sum_{r=0}^{g-1} (u_m-y_0)^{r-1} \!\!\left(
\frac{\varphi^2(Q_0)}{(r+1)!} \frac{\partial^{r+1}}{\partial{y_0^{r+1}}} \left\{\frac{1}{\varphi^2(Q_0)}\right\}
\right.
\\
\left.
 +  \sum_{s=0}^r\frac{\varphi(Q_0)}{(r-s)!} \frac{\partial^{r-s}}{\partial{y_0^{r-s}}} \left\{\frac{1}{\varphi^2(Q_0)}\right\}
\underset{P=Q_0}{\rm res}\frac{\Omega(P)\varphi(P)}{(u-y_0)^{s+1}du}
\right).
\end{multline*}
When using the expression \eqref{residue} for the residue, consider the contribution of the first term $R_1(s+1)$ of \eqref{residue} to the quantity $S$. After changing the summation variable to $s'=s+1$ this contribution simplifies,
\begin{equation*}
\sum_{s'=1}^{r+1}\frac{1}{(r-s'+1)!s'!} \frac{\partial^{r-s'+1}}{\partial{y_0^{r-s'+1}}} \left\{\frac{1}{\varphi^2(Q_0)}\right\} \frac{\partial^{s'} \varphi^2(Q_0)}{\partial y_0^{s'}}
= - \frac{1}{(r+1)!} \frac{\partial^{r+1}}{\partial{y_0^{r+1}}} \left\{\frac{1}{\varphi^2(Q_0)}\right\},
\end{equation*}
and cancels against the same term in $S$. We thus have
\begin{equation}
\label{S}
S
= \frac{1}{(u_m-y_0)^2} + \sum_{r=0}^{g-1} (u_m-y_0)^{r-1} \sum_{s=0}^r\frac{\varphi(Q_0)}{(r-s)!} \frac{\partial^{r-s}}{\partial{y_0^{r-s}}} \left\{\frac{1}{\varphi^2(Q_0)}\right\}
(R_2(s+1)+R_3(s+1))
\end{equation}
with $R_2, R_3$ defined in \eqref{residue}. The contribution of $R_3$ amounts to
\begin{equation*}
\frac{\delta_g}{\prod_{\alpha=1}^{g-1}(y_0-u_\alpha)}\sum_{r=0}^{g-1} \frac{(u_m-y_0)^{r-1}}{r!} \frac{\partial^{r}}{\partial{y_0^{r}}} \prod_{\alpha=1}^{g-1}(y_0-u_\alpha),
\end{equation*}
which vanishes due to Lemma \ref{lemma_E}.

Consider the contribution of $R_2(s+1)$ to $S$ \eqref{S}.
Introducing the summation index $s'=s-1$ and using the Leibniz rule, we have for this contribution
\begin{multline*}
\sum_{i=1}^{g-1} \frac{\delta_i\,\varphi(Q_0)}{\varphi(P_{u_i})(u_i-y_0)\prod_{\alpha\neq i}(u_i-u_\alpha)} \sum_{r=0}^{g-1} \sum_{s'=0}^{r-1}\frac{(u_m-y_0)^{r-1} }{(r-s'-1)!s'!} \frac{\partial^{r-s'-1}}{\partial{y_0^{r-s'-1}}} \left\{\frac{1}{\varphi^2(Q_0)}\right\}
\frac{\partial^{s'}}{\partial y_0^{s'}}\{\varphi^2(Q_0) \prod_{\alpha\neq i}(y_0-u_\alpha)\}
\\
=
\sum_{i=1}^{g-1} \frac{\delta_i\,\varphi(Q_0)}{\varphi(P_{u_i})(u_i-y_0)\prod_{\alpha\neq i}(u_i-u_\alpha)} \sum_{r=1}^{g-1} \frac{(u_m-y_0)^{r-1} }{(r-1)!} \frac{\partial^{r-1}}{\partial{y_0^{r-1}}}
 \prod_{\alpha\neq i}(y_0-u_\alpha).
\end{multline*}
Due to Lemma \ref{lemma_E}, the last sum is zero when $m\neq i$. When $m=i$, due to Lemma \ref{lemma_F}, the last sum is equal to $\prod_{\alpha\neq m}(u_m-u_\alpha)$. Thus the contribution of $R_2(s+1)$ to $S$ reduces to
\begin{equation*}
 \frac{\delta_m\,\varphi(Q_0)}{\varphi(P_{u_m})(u_m-y_0)}
\end{equation*}
and we have for $S$ \eqref{S}
\begin{equation*}
S
= \frac{1}{(u_m-y_0)^2}  + \frac{\delta_m\,\varphi(Q_0)}{\varphi(P_{u_m})(u_m-y_0)}
= \frac{\varphi(Q_0)\Omega(P_{u_m})}{\varphi(P_{u_m})(u_m-y_0)}
\end{equation*}
where the last equality is obtained by evaluating $\Omega(P)$ in the form \eqref{Omega-delta} at $P=P_m$ using the defining relations $v_i(P_{u_j}) = \delta_{ij}$  for the basis of holomorphic differentials $v_i\,.$
Using this result for $S$ in \eqref{temp_factor_of_dumdxi} and definition \eqref{A12_aj_g} of $\Aum^{12}$, we prove the proposition.
$\Box$

Let us now look at the first two lines in \eqref{M10}.

\begin{proposition}
\label{proposition_M10wavyline}
 We have
\begin{multline*}
\frac{1}{x_i-u_m} \left( \frac{\betaum}{\Aum^{12}} - \frac{\betaxi}{\Axi^{12}}\right)^2
-\frac {\prod _\alpha (x_i-u_{\alpha})}{(x_i-y_0)^{g-1}}\sum_{\substack{\alpha=1\\ \alpha\neq m}}^{g-1}\frac {(u_\alpha-y_0)^{g-1}\left( \frac{\betaum}{\Aum^{12}} - \frac{\betaualpha}{\Aualpha^{12}}\right)^2}{(u_\alpha-u_m)(x_i-u_{\alpha})\prod _{\beta \neq \alpha}(u_\alpha-u_{\beta})}
\\
=
\frac{\prod _\alpha (x_i-u_{\alpha})}{(x_i-y_0)^{g-1}}
\left(
\frac{2}{(x_i-y_0)\varphi^2(Q_0)\prod_\alpha(y_0-u_\alpha)}
\sum_{l=0}^{g} \frac{L_l}{l!(u_m-y_0)^{g-l+1}}
\right.
\\
\left.
-1-\frac{4}{\varphi^2(P_{u_m})(u_m-y_0)^{g+1}(x_i-u_m)\prod _{\alpha \neq m}(u_m-u_\alpha)} \right)\,.
\end{multline*}
\end{proposition}

{\it Proof.} Denoting the left hand side by $F$,  and using \eqref{difference_betas_1}, we have
\begin{multline*}
F=
\frac{1}{(x_i-u_m)\varphi^2(Q_0)} \left( \sum_{l_1, l_2=0}^{g-1} \frac{L_{l_1}L_{l_2}}{l_1!l_2!(u_m-y_0)^{2g-l_1-l_2}} + \sum_{l_1, l_2=0}^{g-1} \frac{L_{l_1}L_{l_2}}{l_1!l_2!(x_i-y_0)^{2g-l_1-l_2}}
\right.
\\
\left.
-2\sum_{l_1=0}^{g-1} \frac{L_{l_1}}{l_1!(u_m-y_0)^{g-l_1}}\sum_{l_2=0}^{g-1} \frac{L_{l_2}}{l_2!(x_i-y_0)^{g-l_2}}
\right)
\\
-\frac{\prod _{\alpha} (x_i-u_{\alpha})}{(x_i-u_m)(x_i-y_0)^{g-1}\varphi^2(Q_0)} \left[
\sum_{l_1, l_2=0}^{g-1} \frac{L_{l_1}L_{l_2}}{l_1!l_2!(u_m-y_0)^{2g-l_1-l_2}} \sum_{\alpha\neq m} \frac{(u_\alpha-y_0)^{g-1}}{\prod _{\beta \neq \alpha}(u_\alpha-u_{\beta})} \left( \frac{1}{u_\alpha-u_m}+\frac{1}{x_i-u_\alpha}\right)
\right.
\\
\left.
-2\sum_{l_1=0}^{g-1} \frac{L_{l_1}}{l_1!(u_m-y_0)^{g-l_1}}\sum_{l_2=0}^{g-1} \frac{L_{l_2}}{l_2!} \sum_{\alpha\neq m} \frac{(u_\alpha-y_0)^{l_2-1}}{\prod _{\beta \neq \alpha}(u_\alpha-u_{\beta})} \left( \frac{1}{u_\alpha-u_m}+\frac{1}{x_i-u_\alpha}\right)
\right.
\\
\left.
+ \sum_{N=0}^{g-1}\frac{L_N(2\Sigma_1, \dots, 2\Sigma_N)}{N!}\sum_{\alpha\neq m} \frac{ \frac{1}{u_\alpha-u_m}+\frac{1}{x_i-u_\alpha}}{(u_\alpha-y_0)^{g-N+1}\prod _{\beta \neq \alpha}(u_\alpha-u_{\beta})}
+\sum_{\substack{l_1, l_2=0\\ l_1+l_2=g}}^{g-1} \frac{L_{l_1}L_{l_2}}{l_1!l_2!} \sum_{\alpha\neq m} \frac{ \frac{1}{(u_\alpha-y_0)(u_\alpha-u_m)}}{\prod _{\beta \neq \alpha}(u_\alpha-u_{\beta})}
\right.
\\
\left.
+\sum_{\substack{l_1, l_2=0\\ l_1+l_2=g}}^{g-1} \frac{L_{l_1}L_{l_2}}{l_1!l_2!} \sum_{\alpha\neq m} \frac{ \frac{1}{x_i-y_0}\left(\frac{1}{u_\alpha-y_0}+\frac{1}{x_i-u_\alpha}\right)}{\prod _{\beta \neq \alpha}(u_\alpha-u_{\beta})}
+\sum_{\substack{l_1, l_2=0\\ l_1+l_2> g}}^{g-1} \frac{L_{l_1}L_{l_2}}{l_1!l_2!} \sum_{\alpha\neq m} \frac{(u_\alpha-y_0)^{l_1+l_2-g-1}\left( \frac{1}{u_\alpha-u_m}+\frac{1}{x_i-u_\alpha}\right)}{\prod _{\beta \neq \alpha}(u_\alpha-u_{\beta})}
\right].
\end{multline*}
Now we apply Lemmas \ref{lemma_rational1}, \ref{lemma_rational2} et \ref{lemma_rational3} from Section \ref{sect_rational} to the terms containing sums over $\alpha\neq m$. The only term which is not covered directly by the lemmas, can be represented as a derivative as follows:
\begin{equation*}
\sum_{\alpha\neq m} \frac{ \frac{1}{x_i-u_\alpha}}{(u_\alpha-y_0)^{g-N+1}\prod _{\beta \neq \alpha}(u_\alpha-u_{\beta})}
=
\frac{1}{(g-N)!}
\frac{\partial^{g-N}}{\partial y_0^{g-N}}\left\{
\sum_{\alpha\neq m} \frac{ \frac{1}{x_i-u_\alpha}}{(u_\alpha-y_0)\prod _{\beta \neq \alpha}(u_\alpha-u_{\beta})}
\right\}.
\end{equation*}
A lengthy but straightforward application of Lemmas \ref{lemma_rational1}, \ref{lemma_rational2} et \ref{lemma_rational3} results in the following expression for $F$:
\begin{multline}
\label{Ftemp}
F=
\frac{\prod _{\alpha} (x_i-u_{\alpha})}{\prod_{\alpha}(y_0-u_\alpha)(x_i-y_0)^{g}\varphi^2(Q_0)} \left[
2\sum_{l_1=0}^{g-1} \frac{L_{l_1}}{l_1!(u_m-y_0)^{g-l_1+1}}
-\frac{1}{(u_m-y_0)}\sum_{\substack{l_1, l_2=0\\ l_1+l_2=g}}^{g-1}  \frac{L_{l_1}L_{l_2}}{l_1!l_2!}
   \right]
\\
-\frac{\prod _{\alpha} (x_i-u_{\alpha})}{(x_i-u_m)(x_i-y_0)^{g-1}\varphi^2(Q_0)}
\sum_{N=0}^{g-1}\frac{L_N(2\Sigma_1, \dots, 2\Sigma_N)}{N!\,(g-N)!} \frac{\partial^{g-N}}{\partial y_0^{g-N}}\left\{ \frac{\frac{1}{u_m-y_0} - \frac{1}{x_i-y_0} }{\prod_{\alpha}(y_0-u_\alpha)}   \right\}
\,.
\end{multline}
Now consider the sum over $N$.  Adding and subtracting the term with $N=g$ and using Proposition \ref{prop_Ll}, we can apply the Leibniz rule \eqref{L-Leibniz} to rewrite this sum as follows
\begin{multline*}
\sum_{N=0}^{g-1}\frac{L_N(2\Sigma_1, \dots, 2\Sigma_N)}{N!\,(g-N)!} \frac{\partial^{g-N}}{\partial y_0^{g-N}}\left\{ \frac{\frac{1}{u_m-y_0} - \frac{1}{x_i-y_0} }{\prod_{\alpha}(y_0-u_\alpha)}   \right\}
 = -
 \frac{L_g(2\Sigma_1, \dots, 2\Sigma_g)\left(\frac{1}{u_m-y_0} - \frac{1}{x_i-y_0} \right)}{g!\prod_{\alpha}(y_0-u_\alpha)}
 \\
+
\frac{ \varphi^2(Q_0) }{g!} \frac{\partial^g}{\partial y_0^g} \left\{  \frac{\frac{1}{u_m-y_0} - \frac{1}{x_i-y_0}}{\varphi^2(Q_0)\prod_{\alpha}(y_0-u_\alpha)} \right\}
 \end{multline*}
Given that $\varphi^{-2}(Q_0)=\prod_{a_j}(y_0-a_j)$ we have for the term under the derivative
\begin{multline*}
 \frac{\frac{1}{u_m-y_0} - \frac{1}{x_i-y_0}}{\varphi^2(Q_0)\prod_{\alpha}(y_0-u_\alpha)}
 = (u_m-x_i)\frac{\prod_{a_j\notin\{u_\alpha, x_i\}}(y_0-a_j)}{u_m-y_0}
 \\
 =  -(u_m-x_i)\prod_{a_j\notin\{u_\alpha, x_i, x_j\}}(y_0-a_j)+(u_m-x_i)(u_m-x_j)\frac{\prod_{a_j\notin\{u_\alpha, x_i, x_j\}}(y_0-a_j)}{u_m-y_0}
 \,,
 \end{multline*}
where $\{u_\alpha\}$ denotes the set of all dependent branch points, with $1\leq \alpha\leq g-1$ and to obtain the last equality we represented the factor $y_0-x_j$ with some $j\neq i$ as $y_0-u_m+u_m-x_j\,.$ Repeating the same trick for each factor $y_0-a_j$ in the numerator, we obtain
\begin{equation*}
\frac{\partial^g}{\partial y_0^g} \left\{  \frac{\frac{1}{u_m-y_0} - \frac{1}{x_i-y_0}}{\varphi^2(Q_0)\prod_{\alpha}(y_0-u_\alpha)} \right\}
=
\prod_{a_j\notin\{u_\alpha\}}(u_m-a_j) \frac{\partial^g}{\partial y_0^g} \left\{ \frac{1}{u_m-y_0}  \right\}-g!(u_m-x_i)\,.
\end{equation*}
Finally, using the fact that $4\varphi^{-2}(P_{u_m})=\prod_{a_j\neq u_m}(u_m-a_j)$, we have for the sum over $N$
\begin{multline*}
\sum_{N=0}^{g-1}\frac{L_N(2\Sigma_1, \dots, 2\Sigma_N)}{N!\,(g-N)!} \frac{\partial^{g-N}}{\partial y_0^{g-N}}\left\{ \frac{\frac{1}{u_m-y_0} - \frac{1}{x_i-y_0} }{\prod_{\alpha}(y_0-u_\alpha)}   \right\}
 =- \varphi^2(Q_0) (u_m-x_i)
 \\
- \frac{(x_i-u_m)L_g(2\Sigma_1, \dots, 2\Sigma_g)}{g!(x_i-y_0)(u_m-y_0)\prod_{\alpha}(y_0-u_\alpha)}
+
 \frac{4\varphi^2(Q_0)}{\varphi^2(P_{u_m})(u_m-y_0)^{g+1}\prod_{\alpha\neq m}(u_m-u_\alpha)}
  \,.
 \end{multline*}
Note that the middle term cancels against another term in the first line of \eqref{Ftemp} after using the Leibniz rule to rewrite
the sum over $l_1, l_2$ ranging from $0$ to $g-1$ and such that $l_1+l_2=g$ as $(L_g(2\Sigma_1, \dots, 2\Sigma_g) - 2L_g)/g!$. The term with $2L_g/g!$ is then absorbed by the first sum in \eqref{Ftemp}. Putting these elements together, we prove the proposition.
$\Box$

Plugging the results of Propositions \ref{proposition_factor_of_dumdxi} and \ref{proposition_M10wavyline} into \eqref{M10}, we see that $T_{u_m}=0$ which finishes the proof of the fact that quantities $\Aaj^{11}, \Aaj^{12}$ and $\Aaj^{21}$ given by Theorem \ref{thm_main_g} satisfy equations \eqref{Aum11}.

\subsection{Proof for (21)-components}

The fact that functions $\Aaj^{21}$ defined in Theorem \ref{thm_main_g} satisfy equations \eqref{Aaj21} and \eqref{Aum21} can be obtained as a corollary of the preceding results. Namely, from the definition of $\Aaj^{21}$ we get for any $a_j\in B$
\begin{equation*}
\frac{\partial \Aaj^{21}}{\partial x_i} = -2\frac{\Aaj^{11}}{\Aaj^{12}} \frac{\partial \Aaj^{11}}{\partial x_i} - \frac{\Aaj^{21}}{\Aaj^{12}} \frac{\partial \Aaj^{12}}{\partial x_i}\,.
\end{equation*}
The rest is a simple calculation using the fact that $\Aaj^{12}$ and $\Aaj^{11}$ satisfy  \eqref{Aaj12}, \eqref{Aum12} and \eqref{Aaj11}, \eqref{Aum11}, respectively.

\subsection{Sum of $\Aaj$ is a constant diagonal matrix}

To finish the proof of Theorem \ref{thm_main_g}, it remains to show that the sum of matrices $\Aaj$ over all $a_j\in B$ is constant. Let us denote this sum by $\Ainfty$:
\begin{equation*}
\Ainfty:=\sum_{a_j\in B} A_{a_j}\,.
\end{equation*}
We already know from \eqref{sum_A12} that the $12$ component $\Ainfty^{12}$ of this matrix vanishes. It remains to find $\Ainfty^{11}$ and $\Ainfty^{21}$.
\begin{lemma}
\label{beta_sum}
For the quantities $\betaaj$ defined by \eqref{betas}, we have
\begin{equation*}
\sum_{a_j\in B}  \betaaj = -t\,.
\end{equation*}
\end{lemma}
{\it Proof.}
This is a straightforward calculation using \eqref{res0} and \eqref{res1} from Lemma \ref{lemma_residues}.
$\Box$

From Lemma \ref{beta_sum}, given definition of $\Aaj^{11}$ from Theorem \ref{thm_main_g}, we get the sum of $\Aaj^{11}$ over all branch points to be $ -1/4\,.$
\begin{lemma}
\label{sum21}
For the entries $\Aaj^{21}$ defined in Theorem \ref{thm_main_g}, we have
\begin{equation*}
\sum_{a_j\in B}  \Aaj^{21} = 0\,.
\end{equation*}
\end{lemma}
{\it Proof.}
Given the definition of $\Aaj^{21}$ from Theorem \ref{thm_main_g}, we can write the sum as
\begin{equation*}
\sum_{a_j\in B}  \Aaj^{21} =  - \frac{g}{4t} \sum_{a_j\in B} \frac{\betaaj}{\Aaj^{12}}  - \frac{g^2}{4t^2} \sum_{a_j\in B} \Aaj^{12} \left( \frac{\betaaj}{\Aaj^{12}}  \right)^2\,.
\end{equation*}
Substituting \eqref{betas} for $\betaaj$ and using the sums \eqref{res0}, \eqref{res1} over the branch points from Lemma \ref{lemma_residues}, we have
\begin{multline*}
\sum_{a_j\in B}  \Aaj^{21} =  - \frac{g}{4t}  \frac{1}{\varphi(Q_0)} \sum_{l=0}^{g-1} \frac{L_l}{l!}\sum_{a_j\in B}\frac{1}{(a_j-y_0)^{g-l}} + \frac{g^2(2g+1)}{8t}\Omega(P_\infty)
\\
  - \frac{g^2}{4t^2\varphi^2(Q_0)} \sum_{a_j\in B} \Aaj^{12}  \sum_{l_1,l_2=0}^{g-1} \frac{L_{l_1}L_{l_2}}{l_1!l_2!(a_j-y_0)^{2g-l_1-l_2}}
    - \frac{g^3}{4t}\Omega(P_\infty)      \,.
\end{multline*}
Rewriting the sum over $l_1, l_2$ as a sum over $N=l_1+l_2$, we see that the contribution of the terms with $N>g$ is zero due to \eqref{res0}. Singling out the term with $g=N$ and using \eqref{res1} for it, we have have
\begin{multline}
\label{sum-temp}
\sum_{a_j\in B}  \Aaj^{21} =  - \frac{g}{4t}  \frac{1}{\varphi(Q_0)} \sum_{l=0}^{g-1} \frac{L_l}{l!}\sum_{a_j\in B}\frac{1}{(a_j-y_0)^{g-l}} + \frac{g^2}{8t}\Omega(P_\infty)
\\
  - \frac{g^2}{4t^2\varphi^2(Q_0)}  \sum_{N=0}^{g-1} \sum_{l=0}^N \frac{L_lL_{N-l}}{l!(N-l)!}  \sum_{a_j\in B} \frac{\Aaj^{12}}{(a_j-y_0)^{2g-N}}
    + \frac{g^2}{4t\varphi(Q_0)}   \sum_{l=1}^{g-1} \frac{L_lL_{g-l}}{l!(g-l)!}     \,.
\end{multline}
Let us now work with the sum over $N$ in \eqref{sum-temp}. Note that the sum over $l$ in this term is equal to the $N$-fold derivative of $\varphi^{-2}(Q_0)/N!$ due to the Leibniz rule and derivative \eqref{Ll} of $1/\varphi(Q_0)$ from Proposition \ref{prop_Ll}.  We replace the sum over the branch points using \eqref{res2} from Lemma \ref{lemma_residues} and then calculate the residue from \eqref{res2} as in \eqref{residue}:
\begin{multline*}
\frac{1}{\varphi^2(Q_0)}  \sum_{N=0}^{g-1} \sum_{l=0}^N \frac{L_lL_{N-l}}{l!(N-l)!}  \sum_{a_j\in B} \frac{\Aaj^{12}}{(a_j-y_0)^{2g-N}}=  -t  \sum_{N=0}^{g-1} \frac{1}{N!} \frac{\partial^N}{\partial y^N} \left\{ \frac{1}{\varphi^2(Q_0)}\right\} \,\underset{P=Q_0}{\rm res}\frac{\Omega(P)\varphi(P)}{(u-y_0)^{g-N}du}
  \\
  =-\frac{t}{\varphi(Q_0)}  \sum_{N=0}^{g-1} \frac{1}{N!} \frac{\partial^N}{\partial y^N} \left\{ \frac{1}{\varphi^2(Q_0)}\right\}
  \frac{1}{(g-N)!} \frac{\partial^{g-N} \varphi^2(Q_0)}{\partial y_0^{g-N}}
\\
-t  \sum_{N=0}^{g-2} \frac{1}{N!} \frac{\partial^N}{\partial y^N} \left\{ \frac{1}{\varphi^2(Q_0)}\right\}
\frac{1}{(g-N-2)!} \sum_{i=1}^{g-1} \delta_i
\frac{\frac{\partial^{g-N-2}}{\partial y_0^{g-N-2}}\{\varphi^2(Q_0) \prod_{\alpha\neq i}(y_0-u_\alpha)\}}{\varphi(P_{u_i})(u_i-y_0)\prod_{\alpha\neq i}(u_i-u_\alpha)}
\\
-t  \sum_{N=0}^{g-1} \frac{1}{N!} \frac{\partial^N}{\partial y^N} \left\{ \frac{1}{\varphi^2(Q_0)}\right\}
\frac{\delta_g}{(g-N-1)!} \frac{\frac{\partial^{g-N-1}}{\partial y_0^{g-N-1}}\{\varphi^2(Q_0) \prod_{\alpha}(y_0-u_\alpha)\}}{\varphi(Q_0)\prod_{\alpha}(y_0-u_\alpha)} \,.
\end{multline*}
Now we apply the Leibniz rule in each of the three terms. Note that the first some completed with the term $N=g$ gives the $g$-fold derivative of $1$ and thus vanishes. In the second sum, we have the $(g-2)$-fold $y_0$-derivative of $\prod_{\alpha\neq i}(y_0-u_\alpha)/(g-2)!$ which is a polynomial on degree $g-2$ in $y_0$, thus the second sum over $N$ is equal to $1$. Similarly, the third sum over $N$ gives the $(g-1)$-fold $y_0$-derivative of a polynomial of degree $g-1$ in $y_0$. Thus we have
\begin{multline}
\label{sum-temp-1}
  \frac{1}{\varphi^2(Q_0)}\sum_{N=0}^{g-1} \sum_{l=0}^N \frac{L_lL_{N-l}}{l!(N-l)!}  \sum_{a_j\in B} \frac{\Aaj^{12}}{(a_j-y_0)^{2g-N}}
  \\
  =
  \frac{t\varphi(Q_0)}{g!} \frac{\partial^g}{\partial y^g} \left\{ \frac{1}{\varphi^2(Q_0)}\right\}
-
 \sum_{i=1}^{g-1}
\frac{t\delta_i}{\varphi(P_{u_i})(u_i-y_0)\prod_{\alpha\neq i}(u_i-u_\alpha)}
- \frac{t\delta_g}{\varphi(Q_0)\prod_{\alpha}(y_0-u_\alpha)} \,.
\end{multline}
Now, from the representation \eqref{Omega-delta} of our differential $\Omega$, evaluating it at $P_\infty$ with respect to the local parameter $u^{-1/2}$, we see that $\Omega(P_\infty)=\sum_{i=1}^{g}\delta_i v_i(P_\infty)$. Computing $v_i(\Pinfty)$ from definitions \eqref{v} and \eqref{vg} of the holomorphic differentials $v_i$, we see that the last two terms in \eqref{sum-temp-1} give exactly $t\Omega(\Pinfty)/2$, that is we get the following relation
\begin{equation*}
 \frac{1}{\varphi^2(Q_0)} \sum_{N=0}^{g-1} \sum_{l=0}^N \frac{L_lL_{N-l}}{l!(N-l)!}  \sum_{a_j\in B} \frac{\Aaj^{12}}{(a_j-y_0)^{2g-N}}
  =
  \frac{t\varphi(Q_0)}{g!} \frac{\partial^g}{\partial y^g} \left\{ \frac{1}{\varphi^2(Q_0)}\right\} + \frac{t}{2}\Omega(\Pinfty)\,.
  \end{equation*}
Plugging this back into \eqref{sum-temp}, we obtain
\begin{equation*}
\sum_{a_j\in B}  \Aaj^{21} =  - \frac{g}{4t}  \frac{1}{\varphi(Q_0)} \sum_{l=0}^{g-1} \frac{L_l}{l!}\sum_{a_j\in B}\frac{1}{(a_j-y_0)^{g-l}}
  - \frac{g\varphi(Q_0)}{4t(g-1)!}   \frac{\partial^g}{\partial y^g} \left\{ \frac{1}{\varphi^2(Q_0)}\right\}
    + \frac{g^2}{4t\varphi(Q_0)}   \sum_{l=1}^{g-1} \frac{L_lL_{g-l}}{l!(g-l)!}     \,.
\end{equation*}
Adding and subtracting the terms with $l=0$ and $l=g$ in the last sum and using again \eqref{Ll} from Proposition \ref{prop_Ll},  we can convert the completed sum by the Leibniz rule into $\varphi^2(Q_0)/g!$ times the $g$-fold derivative of ${1}/{\varphi^2(Q_0)}$ with respect to $y_0$, which cancels the middle term in the above line. Thus we have
\begin{equation*}
\sum_{a_j\in B}  \Aaj^{21} =  - \frac{g}{4t}  \frac{1}{\varphi(Q_0)} \sum_{l=0}^{g-1} \frac{L_l}{l!}\sum_{a_j\in B}\frac{1}{(a_j-y_0)^{g-l}}
    - \frac{g^2}{2t\varphi(Q_0)}  \frac{L_g}{g!}     \,.
\end{equation*}
Using again \eqref{Ll} from Proposition \ref{prop_Ll} and then the Leibniz rule, we can see that this is zero rewriting the sum over $l$ as follows
\begin{multline*}
\sum_{a_j\in B}  \Aaj^{21} =  - \frac{g}{4t}  \sum_{l=0}^{g-1} \frac{1}{l!}\frac{\partial^l}{\partial y_0^l}\left\{\frac{1}{\varphi(Q_0)}\right\}\frac{1}{({g-1-l})!}\frac{\partial^{g-1-l}}{\partial y_0^{g-1-l}}\left\{\sum_{a_j\in B}\frac{1}{a_j-y_0}  \right\}
    - \frac{g^2}{2t\varphi(Q_0)}  \frac{L_g}{g!}
\\
=- \frac{g}{4t(g-1)!}   \frac{\partial^{g-1}}{\partial y_0^{g-1}}\left\{\frac{1}{\varphi(Q_0)}\sum_{a_j\in B}\frac{1}{a_j-y_0}  \right\}
    - \frac{g}{2t\varphi(Q_0)}  \frac{L_g}{(g-1)!}\,.
\end{multline*}
It remains to rewrite $L_g$ as derivative using Proposition \ref{prop_Ll} again :
\begin{equation*}
\frac{L_g}{\varphi(Q_0)}=\frac{\partial^{g}}{\partial y_0^{g}}\left\{\frac{1}{\varphi(Q_0)} \right\}
=\frac{\partial^{g}}{\partial y_0^{g}}\left\{\sqrt{\prod_{a_j\in B} (y_0-a_j)} \right\}
= -\frac{1}{2}
\frac{\partial^{g-1}}{\partial y_0^{g-1}}\left\{\frac{1}{\varphi(Q_0)}\sum_{a_j\in B}\frac{1}{a_j-y_0}  \right\}
\,.
\end{equation*}
where we used \eqref{phi_Q0} for $\varphi(Q_0)$.
$\Box$

Summarizing, we have
\begin{equation*}
\Ainfty:=\left(\begin{array}{cc}-\frac{1}{4} & 0 \\ 0 & \frac{1}{4}\end{array}\right)\,.
\end{equation*}

This concludes the proof of Theorem \ref{thm_main_g}.

$\Box$

\section[Poncelet theorem and billiards. Combs. Inheritence problem.]{Links to Poncelet theorem and billiards within quadrics.\\ Rectification on combs. Inheritance problem}
\label{sect_Poncelet}

In \cite{Hitchin1}, Hitchin related a construction of algebraic solutions  to PVI $(\frac{1}{8}, -\frac{1}{8}, \frac{1}{8}, \frac{3}{8})$ to the classical Cayley condition characterizing pairs of conics in the plane with the property that there exist closed polygonal lines with prescribed number of sides, inscribed in one of the
conics and circumscribed about the other. Such closed polygonal lines are called Poncelet polygons. In the case when the conics are confocal, the Poncelet polygons are periodic billiard trajectories, where one of the conics is the billiard boundary and another one is the caustic of billiard trajectories. All the trajectories within the same boundary conic have the same period if they are periodic. More generally, the trajectories, with the same boundary and the caustic, have the same rotation number, which is rational in the periodic cases and irrational
otherwise, see e.g. \cite{DR2011, DR2019b}. As explained in \cite{DS2019, DS2021}, in the context of elliptical billiards in the plane, the Hitchin construction can be related to a simultaneous deformation of both the boundary conic and the caustic in a way that corresponding billiard trajectories
remain periodic with the same period. Then, \cite{DS2019, DS2021} went further and generalized Hitchin's consideration to non-periodic cases:  a general solution to the above Painlev\'e equation was linked to a  deformation of
both the boundary conic and the caustic in a way that the rotation number, which is not necessarily rational any more, remains preserved.

\smallskip

Elliptical billiards in the plane have higher-dimensional generalizations to billiards within quadrics in $\mathbb R^d$, see e.g. \cite{DR2011} and references therein. In \cite{DS2017a} a generalization of Hitchin's idea was presented giving solutions to non-constrained Schlesinger systems, which  were linked  there to  the so-called billiard ordered games (\cite{DR2006})  in $\mathbb R^d$  within $d-1$ confocal quadrics as boundaries and  other $d-1$ confocal quadrics as caustics.

\smallskip

A natural question remained open:

\begin{itemize}

\item  What would be an analogue of Hitchin's construction of solutions of isomonodromic deformations which would be related to a
usual billiard with the boundary consisting of {\it one single}  quadric  in $\mathbb R^d$ for $d>2$?

\end{itemize}

The  aim of the following Section \ref{sec_billiard} is to demonstrate that this  question has an answer in the framework  of isoharmonic deformations, constrained variations of Jacobi inversion, and constrained Schlesinger systems.

\smallskip

Let us observe that a geometric construction of four-periodic families of Poncelet polygons inscribed in a given circle with caustics belonging to a confocal family of conics was recently presented in \cite{DR2021}. The period-preserving flows in nonlinear Schr\"odinger systems were studied in \cite{GS}.

\medskip

In Section \ref{sec_SC} we construct a rectification of isoharmonic deformations in comb regions,
by applying classical ideas of Marchenko-Ostrovsky \cite{MO1975}. We conclude this Section and the paper in Section \ref{sec:inher} with a brief discussion on the inheritance problem and the injectivity of the frequency map.

\smallskip

\subsection{Applications to billiards in  $d$-dimensional space}
\label{sec_billiard}
Let us now consider billiards in domains in $d$-dimensional space bounded by the ellipsoid $\mathcal C_0$ from the confocal
family of quadrics:
\begin{equation}\label{eq:confocal1}
\mathcal{C}_{\lambda}\ :\ \frac{x_1^2}{b_1-\lambda}+\dots
+\frac{x_d^2}{b_d-\lambda}=1, \quad\lambda\in\mathbb{R},
\end{equation}
with $0<b_1<b_2<\dots < b_d$ being constants. These constants $b_j$ are going to be varied later on to induce the dynamics which is an object of our study. Any point of the space $\mathbb R^d$   belongs to exactly  $d$ quadrics from the family
(\ref{eq:confocal1}). The parameters $\lambda_1,\dots,$ $\lambda_d$ of these quadrics are called {\it the Jacobi
coordinates} of the point in  $\mathbb R^d$  of  Cartesian coordinates $(x_1,\dots, x_d)$ (see \cite{DR2019b} and references therein), with the convention: $\lambda_1<\lambda_2<\dots<\lambda_d$. The billiard flow has $d-1$ quadrics from the confocal family (\ref{eq:confocal1}) as {\it caustics},  they correspond to some values $\alpha_1, \dots, \alpha_{d-1}$ of  the parameter $\lambda$ from (\ref{eq:confocal1}). The billiard domain $\Gamma$   consists of all points of the billiard trajectories within $\mathcal{C}_0$ with caustics $\mathcal{C}_{\alpha_1}$, \dots, $\mathcal{C}_{\alpha_{d-1}}$. In Jacobi coordinates, $\Gamma$ is just a product of $d$ segments in $[0,+\infty)$; note that    $\hat \Delta_{2d-1}(\lambda)$,  defined in (\ref{eq:hyperelliptic}), is positive  on $\Gamma$. Along any billiard trajectory each of the Jacobi coordinates changes monotonically within  a certain  segment until it reaches an endpoint. Once a Jacobi coordinate reaches a segment endpoint, it reverses the direction of its motion.
Endpoints of the segments correspond to: i) reflection off of the boundary ellipsoid ($\lambda=0$); ii) touching a caustic ($\lambda=\alpha_k, k=1,\dots, d-1$); iii) crossing a coordinate hyper-plane ($\lambda=b_j, j=1, \dots, d$). For a given periodic trajectory of period $n$, denote {\it the winding numbers}: $m_0=n,\; m_1, \dots, m_{d-1}$,
where $m_i$ for $i=1,\dots, d-1$ is equal to the number of turns the Jacobi coordinate $\lambda_i$ makes in its interval during one period of the billiard trajectory. The billiard dynamics can be analyzed on the Jacobian of the hyperelliptic curve:
\begin{equation}\label{eq:hyperelliptic}
\mu^2=\hat \Delta_{2d-1}(\lambda)=(\lambda-b_1)\cdots(\lambda-b_d)(\lambda-\alpha_1)\cdots
(\lambda-\alpha_{d-1}).
\end{equation}

We will use the notation $\{e_1,\dots, e_{2d-1}\}=\{b_1,\dots, b_d,\alpha_1,\dots,\alpha_{d-1}\}$ with the ordering condition
$e_1<e_2<\dots<e_{2d-1}$ and $c_{2d}=0, c_j=1/e_j$.

\begin{lemma}[\cite{Audin}]\label{lemma:audin}
If $\alpha_1<\alpha_2< \dots<\alpha_{d-1}$, then
$\alpha_j\in\{e_{2j-1},e_{2j}\}$, for $1\le j\le d-1$.
\end{lemma}

\smallskip

In particular, Lemma \ref{lemma:audin} implies $e_{2d-1}=b_d$.

\smallskip

Following \cite{DR2019b} we use  the polynomial $\Delta_{2d}(s)=s\prod_{j=1}^{2d-1}\left(s-c_j\right).$
The connection between ellipsoidal billiards in dimension $d$ and Chebyshev extremal polynomials on the systems of $d$ intervals was fully established in \cite{DR2019b}. Exploiting this connection all three  Ramirez-Ros conjectures \cite{RR2014} were proved there and all sets of caustics which generate periodic trajectories of  a given period $n$ of billiards within a given ellipsoid in $d$-dimensional space were classified.  In addition, the injectivity properties of the billiard frequency maps which were proved  in Theorem 13 of  \cite{DR2019b}
 play an important role in the current setting, as it was explained in Section \ref{sect_isoharmonic}.  It was proved in \cite{DR2019b} that the generalized Cayley condition $C(n,d)$ from \cite{DR2006}   describing $n$ - periodic trajectories in $d$-dimensional space with caustics $\alpha_1,\dots,\alpha_{d-1}$ is satisfied if and only if there exist a pair of real polynomials
$\mathcal{P}_n$, $\mathcal{Q}_{n-d}$ of degrees $n$ and $n-d$ respectively such that the Pell equation holds:
\begin{equation}\label{eq:Pell}
\mathcal P_n^2(s)-{\Delta}_{2d}(s)\mathcal {Q}_{n-d}^2(s)=1.
\end{equation}

The polynomial $\mathcal P_n$ is a rescaled  extremal polynomial of the Chebyshev type on the system of $d$ intervals
$E=[c_{2d},c_{2d-1}]\cup[c_{2d-2},c_{2d-3}]\cup\dots\cup[c_2,c_1]$.
If we denote the signature by $(\tau_1,\dots,\tau_d)$, where  $\tau_j$ is  the number of zeroes of $\mathcal Q_{n-d}$ in $(c_{2j},c_{2j-1})$, then the following relations were established \cite{DR2019b} for the winding numbers:
\begin{equation}\label{eq:signature}
m_{j}=m_{j+1}+\tau_{j}+1, \quad 1\le j\le d.
\end{equation}
 Thus, we see that the degree of a Chebyshev polynomial and its signature determine the winding numbers and vice-versa.
There are relations
between the winding numbers and the frequencies $f_s$ as defined in \eqref{eq:freq}, see \cite{DR2019b}: $f_s=m_{d-s}/m_0$.

\smallskip
 The aim of this Section is to answer the following question:

\begin{itemize}

\item {\it Consider  a billiard system  with a given  ellipsoid in $\mathbb R^d$  as a boundary  and   a given  set of $d-1$  quadrics confocal with the boundary as the caustics,  which are tangent to an $n$-periodic trajectory. How to describe the dynamics of the caustics under the variation of semi-axes of the boundary ellipsoid, such that the generated billiard trajectories remain periodic with the same period and the same winding numbers?}

\end{itemize}

\smallskip

According to
Lemma \ref{lemma:audin}, we see that in the billiard setting for the first interval we have $c_{2d}=0, c_{2d-1}=1/b_{d}$. Moreover, for the remaining $d-1$ intervals $[c_{2j}, c_{2j-1}],$ with $j=1,\dots, d-1,$ exactly one endpoint of each of these intervals is $1/b_j$, i.e. a reciprocal of the square of a semi-axis, and another endpoint is the reciprocal of the parameter $\alpha_j$ of a caustic $\mathcal C_{\alpha_j}$. We then define maps $\sigma:\{1, 2,\dots, d-1\}\rightarrow \{\ell, r\}$ which are consistent with
the billiard dynamics by the rule that $\sigma (j)=\ell$ if the caustic parameter is the left endpoint of the $j$-th interval, and $\sigma (j)=r$ otherwise.

\smallskip

Let us consider one example.

\begin{example}
\label{ex:P5}
{\rm
Let us consider the Chebyshev polynomial $\mathcal P_5(s)$ of degree $5$ which is the extremal polynomial over three intervals  $[c_6, c_5]\cup [c_4, c_3]\cup [c_2, c_1]$, see Fig. \ref{fig:p542}. It appeared in \cite{DR2019b} in the study of $5$-periodic billiard trajectories in
$\mathbb R^3$.  Corresponding billiard trajectories have
 the winding numbers $(5,4,2)$. Thus, the signature is $(0,1,1)$. There are two geometrically distinct situations possible in this situation: i) $(c_3,c_4)=(1/b_2,1/\alpha_2)$; and ii) $(c_4,c_3)=(1/b_2,1/\alpha_2)$. The first corresponds to $\sigma_1: \{1, 2\}\rightarrow \{\ell, r\}$, defined by $\sigma_1(1)=r,\, \sigma_1(2)=\ell$. The second corresponds to $\sigma_2: \{1, 2\}\rightarrow \{\ell, r\}$, defined by $\sigma_2(1)=\ell,\, \sigma_2(2)=\ell$. Let us consider one of these options, say i). We are interested in the dynamics of $c_1$ and $c_4$ as the reciprocal values of the caustic parameters $\alpha_1$ and $\alpha_2$ when $c_2$ and $c_3$ (corresponding to the reciprocal values of $b_3$ and $b_2$ respectively) get deformed, provided that the Pell equation remain satisfied during the entire process and the signature $(0,1,1)$ preserved.
}
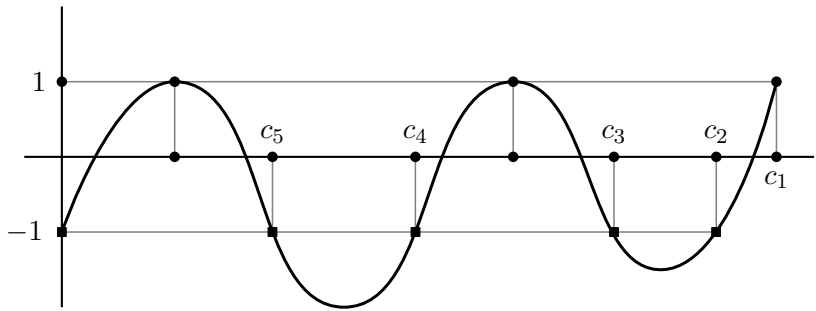
\begin{figure}[h]
\centering
\begin{pspicture}(-0.5,-2.5)(10.5,2.5)

\psline[linecolor=gray,linewidth=0.02](0,1)(9.5,1)
\psline[linecolor=gray,linewidth=0.02](0,-1)(8.7,-1)
\psline[linecolor=gray,linewidth=0.02](1.5,1)(1.5,0)
\psline[linecolor=gray,linewidth=0.02](2.8,-1)(2.8,0)
\psline[linecolor=gray,linewidth=0.02](4.7,-1)(4.7,0)
\psline[linecolor=gray,linewidth=0.02](6,1)(6,0)
\psline[linecolor=gray,linewidth=0.02](7.34,-1)(7.34,0)
\psline[linecolor=gray,linewidth=0.02](8.7,-1)(8.7,0)
\psline[linecolor=gray,linewidth=0.02](9.5,1)(9.5,0)

\psline(-0.5,0)(10,0)
\psline(0,-2)(0,2)

\pscurve[linewidth=0.04](0,-1)(1.5,1)(3.75,-2)(6,1)(8,-1.5)(9.5,1)

\pscircle*(0,1){0.07}
\rput(-0.3,1){$1$}

\psdots[linecolor=black,fillcolor=black,dotstyle=square,dotsize=0.14](0,-1)
\rput(-0.5,-1){$-1$}

\pscircle*(1.5,1){0.07}
\pscircle*(6,1){0.07}
\pscircle*(9.5,1){0.07}


\psdots[linecolor=black,fillcolor=black,dotstyle=square,dotsize=0.14]
(2.8,-1)(4.7,-1)(7.34,-1)(8.7,-1)

\pscircle*(1.5,0){0.07}

\pscircle*(2.8,0){0.07}
\rput(2.8,0.3){$c_{5}$}

\pscircle*(4.7,0){0.07}
\rput(4.7,0.3){$c_{4}$}

\pscircle*(6,0){0.07}

\pscircle*(7.34,0){0.07}
\rput(7.34,0.3){$c_{3}$}

\pscircle*(8.7,0){0.07}
\rput(8.7,0.3){$c_{2}$}

\pscircle*(9.5,0){0.07}
\rput(9.5,-0.3){$c_{1}$}

\end{pspicture} 
\caption{From \cite{DR2019b}: the graph of $\mathcal P_5(s)$. The endpoints are $c_1=1/\alpha_1$, $c_2=1/b_3$, $\{c_3,c_4\}=\{1/b_2,1/\alpha_2\}$, $c_{5}=1/b_1$, $c_6=0$. The signature is $(0,1,1)$. Corresponding billiard trajectories have
the winding numbers $(5,4,2)$.}\label{fig:p542}
\end{figure}

\end{example}

\smallskip
Now, we can formulate answers to the questions posed above in this section.  Let us recall that we defined the Chebyshev dynamics
at the beginning of Section \ref{sect_Chebyshev}.

\smallskip

\begin{theorem}\label{thm:billiard}  Let an ellipsoid in $\mathbb R^d$  be elected as the boundary of a billiard system together with a set of $d-1$ caustics $\mathcal C_{\alpha_j}$, the quadrics confocal with the boundary which are tangent to an $n$-periodic trajectory with given winding numbers $(m_0=n, m_1,\dots, m_{d-1})$. If the semi-axes of the boundary ellipsoid vary,  then assuming that the winding numbers remain preserved, the caustics $\mathcal C_{\alpha_j}$ will deform in a way that the reciprocals $u_j=1/\alpha_j$ of their confocal parameters  obey the Chebyshev dynamics as dependent variables (i.e. dependent end-points of intervals), while $x_j = 1/b_j, \, j=1, \dots, d-1,$ the reciprocal of the squares of semi-axes, are the independent variables (i.e. independently varying end-points of intervals).
\end{theorem}

\smallskip

Thus, the dynamics of the caustics of periodic billiard trajectories generates isoharmonic deformations if the semi-axes of the boundary ellipsoid vary and the winding numbers are preserved. Moreover, these isoharmonic deformations are isoequilibrium.
The families of deformations obtained in Theorem \ref{thm_main_g} with $\c_1$ being rational, $\c_2=0$ correspond to the dynamics of the caustics of periodic billiard trajectories when the semi-axes of the boundary ellipsoid vary and the winding numbers are preserved  in the same way as the algebraic solutions of PVI$(1/8; -1/8; 1/8; 3/8)$ corresponded to deformations of Poncelet polygons in Hitchin's work \cite{Hitchin1}, see also \cite{DS2019, DS2021}. More generally, in the case of non-rational $\c_1$ with $\c_2=0$, they relate to  the dynamics of the caustics under the variation of semi-axes of the boundary ellipsoid, provided that the generated billiard trajectories preserve their frequencies.

\bigskip

\subsection{Comb regions and rectification of the deformations}
\label{sec_SC}

Following the idea which goes back to Marchenko-Ostrovsky \cite{MO1975}, let us consider the Schwarz--Christoffel map $\theta$ generated by the Green function. This is a conformal map from the upper half-plane $\mathbb H=\{z|\Im z> 0\}$ to a comb region $\mathcal C$, a vertical semi-strip with vertical slits  mapping   the point at infinity to the point at infinity.  A  comb region $\mathcal C$ is defined as follows:
$$
\mathcal C=\{w|\, \Im w>0, 0\le \Re w\le 1\}\setminus\bigcup_j^{d-1}\{w|\, \Re w = q_j, 0\le \Im w\le h_j\}
$$
with some real values $q_1,\dots, q_{d}$ and $h_1, \dots, h_d$.

\begin{figure}[h]
\centering
\includegraphics[width=10cm,height=4cm]{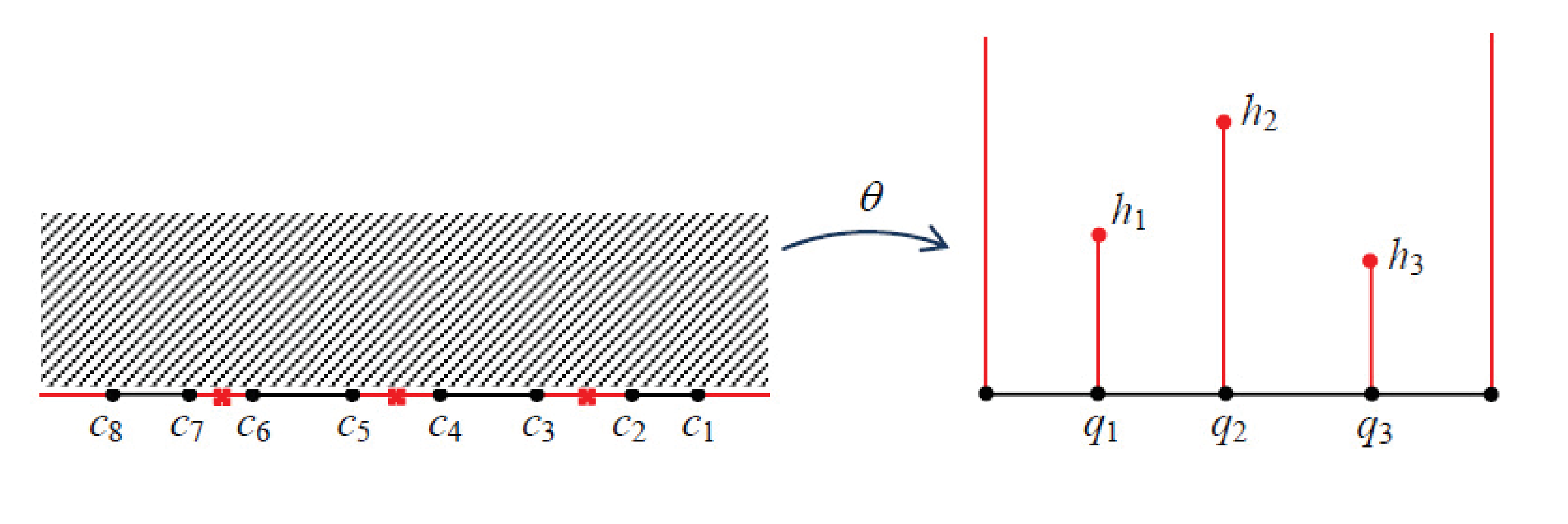}
\caption{Map $\theta$ for $d=4$, $g=3$; the crosses in $(c_7, c_6), (c_5, c_4), (c_3, c_2)$ denote the zeros of $\eta$ and map to $h_1, h_2, h_3$.}
\label{fig:Comb1}
\end{figure}

Using the notation of Section \ref{sect_isoharmonic}, we consider the set $E$ as the union of $d$ intervals, the frequency map $F$ given by \eqref{eq:freq} with the components  denoted by $(f_1, \dots, f_{d-1})$, and the Green function of $E^c$ with the pole at infinity
$$
G_{E}(z, \infty)=\int_{c_{2d}}^z\eta.
$$
Then
\begin{equation}
\label{theta}
\theta=-iG_E(z, \infty):\mathbb H\to \mathcal C,
\end{equation}
where $\mathcal C$ is defined with $q_j=f_j$, $j=1\dots, d-1$. The points  $(q_j, h_j)$ are the $\theta$-images  of the zeros of the differential $\eta$ from the gap intervals $(c_{2j-1}, c_{2j})$ (see  Remark \ref{rem:eta} and  Fig. \ref{fig:Comb1}).

The Chebyshev supports, which is the same as $n$-regular sets, are characterized by  $q_j$'s being integers  (up to a common factor), $j=1,\dots, d-1$  (see e.g. \cite{SY1992}). This corresponds to the conditions $\c_1\in {\mathbb Q}^g, \c_2=0$. In such a case, according to \cite{Bogatyrev2012}, we have
\begin{equation}
	\label{OmegaB}
	\Omega_{\mathbb A}(p) = \frac{1}{n}\frac{d{\mathbb A}(p)}{{\mathbb A}(p)}=\frac{1}{n}\frac{{\mathcal P_n}'(\lambda)d\lambda}{{\mathcal Q_{n-2}}(\lambda) \sqrt{\Delta_{2d}(\lambda)}} = \frac{\prod_j^{d-1}(\lambda-\hat{\gamma_j})d\lambda}{\sqrt{\Delta_{2d}(\lambda)}},
	\end{equation}
where $\mathbb A$ is the Akhiezer function and $\Omega_{\mathbb A}$  is the differential defined in \eqref{eq:Akhiezer}, which coincides with $\eta$ in this case, as was pointed out in Section \ref{sect_isoharmonic}.
Here $\hat{\gamma_j}\in (c_{2j}, c_{2j-1})$ are the gap critical points of the polynomial ${\mathcal P_n}$. The images  $(q_j, h_j)=\theta(\hat {\gamma_j})$ are the tips of the vertical slits at $q_j$, $j=1,\dots, d-1$.

\begin{example} {\rm
This is a continuation of Example \ref{ex:P5}. The image of the map $\theta$ is presented in Fig. \ref{fig:Comb5}.
Here $d=3, \;g=2$.  The values $q_1=2/5$ and $q_2=4/5$ are obtained from the data of the Chebyshev polynomial  $\mathcal P_5$: its degree being $5$, and the number of critical points in the intervals, namely one critical point in each of the intervals $(c_6, c_5)$ and $(c_4, c_3)$, and none in the interval $(c_2, c_1).$  The  $\theta$-images of the gap critical points of the polynomial $\mathcal P_5,$ those in  the gap intervals $(c_5, c_4), (c_3, c_2),$ are the points $(q_1, h_1)$ and $(q_2,h_2)$.   The winding numbers are $(m_0, m_1, m_2)=(5, 4, 2)$ and the frequencies are $m_2/m_1=f_1=q_1=2/5$ and $m_1/m_0=f_2=q_2=4/5$.
}

\begin{figure}[h]
\centering
\includegraphics[width=10cm,height=4cm]{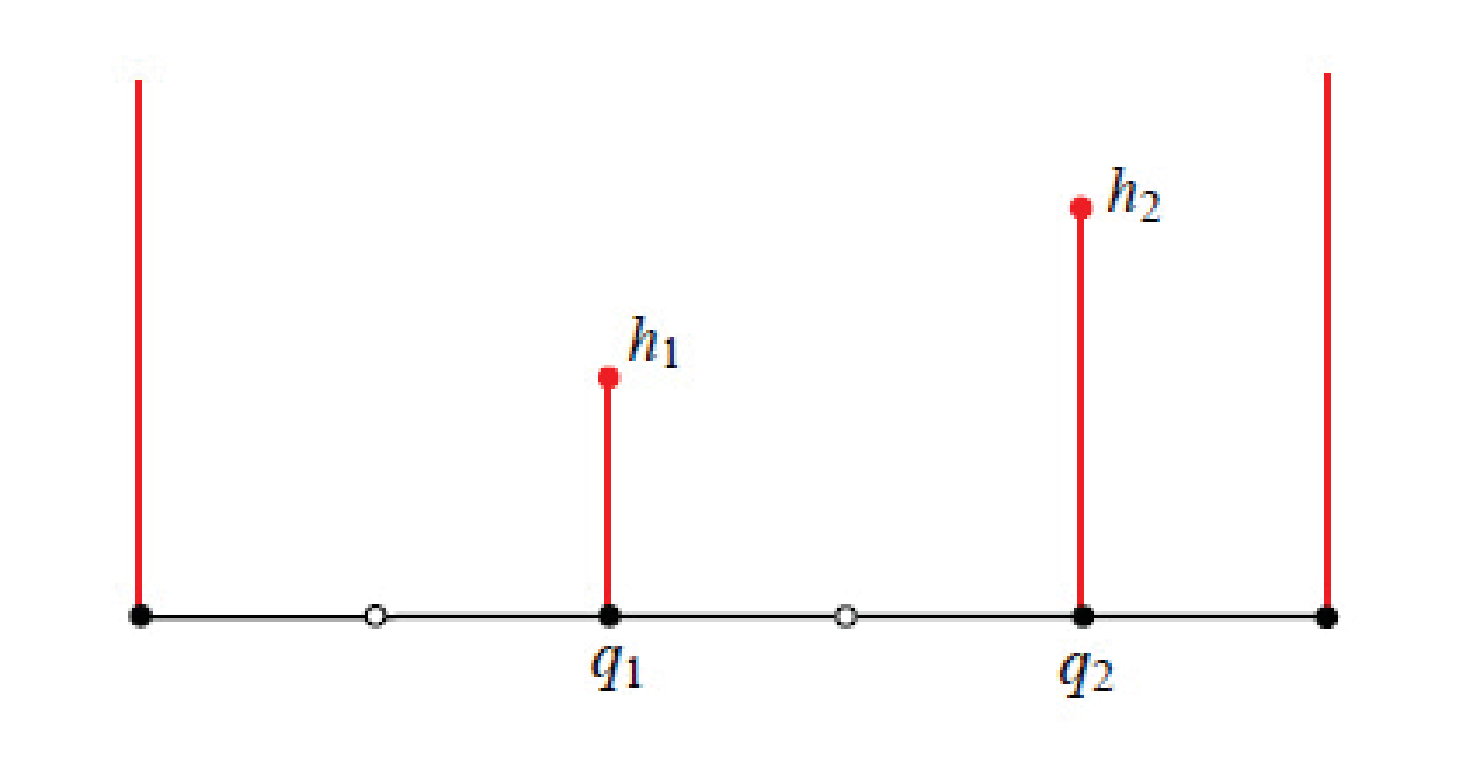}
\caption{The comb region for $\mathcal P_5$ from Example \ref{ex:P5}}
\label{fig:Comb5}
\end{figure}

\end{example}

\bigskip

We conclude with the following statement about explicit rectification of isoharmonic deformations.

\smallskip

\begin{theorem}\label{thm:SC} Consider the map $\theta$ \eqref{theta} corresponding to  the Green function of the complement of the union $E$ of $d$ real intervals, from the upper half-plane to the comb region defined by points $q_j, h_j$, $j=1, \dots, d$.  Any isoharmonic deformation of $(E^c, \infty)$ keeps the points $q_j$ unchanged, i.e. the base of the comb is invariant under the isoharmonic deformations. The deformation
applies only vertically by varying $h_j$, the $\theta$-images of the critical points of the Green function, along the vertical rays based at $q_j$, for each $j=1,\dots, d-1$.
\end{theorem}

\subsection{The inheritance problem and the injectivity of the frequency map}\label{sec:inher}

One of the problems of the 1991 Mikl\'os Schweitzer Mathematical Contest organized
by the J\'anos Bolyai Mathematical Society in Hungary, with a slight generalization was \cite{Totik2006, Totik2009}:

{\it
To divide an inheritance, $d>2$ siblings  turn to a judge. Secretly however, each of them
bribes the judge. What a given sibling inherits depends continuously and monotonically
on the bribes: it   increases monotonically with respect to the value of  their own bribe and  decreases monotonically with respect to the value of the bribe
of everybody else. Show that if the eldest sibling does not give too much  to the judge,
then the others can choose their bribes so that the decision will be fair, i.e., each of them
gets an equal share.}

In \cite{Totik2006, Totik2009}, this problem, its variations, and more precise formulations were discussed in the framework of  the so-called {\it monotone systems}:
denote by ${\bf t}=(t_1, \dots, t_d)$ the set of bribes, with $t_j$ being the bribe of the $j$-th sibling and $h_k(\bf t)$ denotes the share
obtained by the $k$-th sibling with the bribes being equal to $\bf t$.  Then these parameters satisfy  (i) each $h_k$ is a continuous function; (ii) each $h_k$ is  strictly increasing in $t_k$ and strictly  decreasing in $t_j$ for all $j\ne k$; (iii) the sum is constant: $\sum_{j=1}^dh_j({\bf t})=1$.

There are various incarnations of the monotone systems, see \cite{Totik2006}. Here we are particularly interested in the following  situation, see \cite{Totik2009}.  Using the notation from Section \ref{sect_isoharmonic}, let
$$
E=[c_{2d},c_{2d-1}]\cup[c_{2d-2},c_{2d-3}]\cup\dots\cup[c_2,c_1] \,\quad\text{with}\quad  c_{2d}<c_{2d-1}<\dots<c_1.
$$
Assume that each $t_j$ is nonnegative and smaller than $\hat t=\min_k(c_{2k-2}-c_{2k-1})/2, \, k=1, \dots, d-1,$  and
$$
E({\bf t})=[c_{2d},c_{2d-1}+t_d]\cup[c_{2d-2},c_{2d-3}+t_{d-1}]\cup\dots\cup[c_2,c_1+t_1] \,\quad\text{with}\quad  c_{2d}<c_{2d-1}<\dots<c_1.
$$
Then
$$
h_j({\bf t}):={\mathcal M}_{E(\bf t)}[c_{2j}, c_{2j-1}+t_j],
$$
form a monotone system.  Given the existence of a trivial solution ${\bf t}={\bf 0}$,  the positive answer to the inheritance problem may seem in contradiction with the injectivity of the frequency
map defined in \eqref{eq:freq}.   However, there is no contradiction because in Section \ref{sect_isoharmonic} (and in the rest of this paper)
we assume that {\it one interval remains fixed}, which, in the language of the inheritance problem, means that we assume that for one $j\in\{1, 2, \dots, d\}$ the corresponding $t_j$ is zero. We come to the following:

\begin{proposition} If one of the intervals does not participate in the bribing scheme, i.e. if there exists $j$ such that $t_j=0$ then the inheritance problem does not have any fair solution, apart from the trivial one when all other intervals also do not participate in the bribing, i.e.
when ${\bf t}={\bf 0}$.
\end{proposition}

\bigskip

\bigskip

{\bf Acknowledgements.} We thank Patrick Labelle for very helpful discussions concerning the polynomials $L_l$ and Jean-Philippe Burelle for his  very  useful comments upon reading parts of the manuscript.  V.D. acknowledges with gratitude the Simons Foundation grant no. 854861 and the support from the University of Texas at Dallas, MISANU, the Serbian Ministry of Education, Science, and Technological Development, the Science Fund of Serbia. V.S. gratefully acknowledges
support from the Natural Sciences and Engineering Research Council of Canada through a Discovery grant and from the University of Sherbrooke.

\end{document}